\documentclass{amsart}
\usepackage{amssymb,amsmath,amsfonts,euscript,mathrsfs,amsthm}
\usepackage[dvips]{graphicx}

\usepackage[dvips]{color}

\input xy
\xyoption{all}
\CompileMatrices

\def\risom{\overset{\sim}{\rightarrow}}

\DeclareMathOperator{\Id}{Id}

\def\cA{{\mathscr A}}
\def\cB{{\mathscr B}}

\def\cD{{\mathscr D}}
\def\cE{{\mathscr E}}
\def\cF{{\mathscr F}}

\def\cH{{\mathscr H}}

\def\cK{{\mathscr K}}
\def\cJ{{\mathscr J}}
\def\cL{{\mathscr L}}
\def\cM{{\mathscr M}}
\def\cN{{\mathscr N}}
\def\cO{{\mathscr O}}

\def\cT{{\mathscr T}}

\def\cV{{\mathscr V}}

\def\cY{{\mathscr Y}}

\def\fg{{\mathfrak g}}

\def\fm{{\mathfrak m}}

\def\fl{{\mathfrak l}}

\def\Hom{\operatorname{Hom}}

\def\End{\operatorname{End}}

\def\Spec{\operatorname{Spec}}

\def\res{\operatorname{res}}

\def\vac{|0\rangle}

\newtheorem{thm}{Theorem}[section]
\numberwithin{equation}{thm}

\newtheorem{prop}[thm]{Proposition}
\newtheorem{lem}[thm]{Lemma}

\theoremstyle{definition}
\newtheorem{defn}[thm]{Definition}

\newtheorem{ex}[thm]{Example}

\theoremstyle{remark}

\newtheorem{rem}[thm]{Remark}

\newtheorem{nolabel}[thm]{}

\DeclareMathOperator{\ber}{Ber}
\DeclareMathOperator{\lie}{Lie}
\DeclareMathOperator{\sdet}{sdet}
\newcommand{\ys}{Y}
\newcommand{\cys}{\mathscr{Y}}
\begin{document}
\title{SUSY Vertex algebras and supercurves}
\author{Reimundo Heluani}
\address{Department of Mathematics, MIT \\ Cambridge, MA 02139 \\ USA}
\email{heluani@math.mit.edu}
\date{}
\begin{abstract}This article is a continuation of 
	\cite{heluani3}.  Given a strongly conformal SUSY vertex algebra $V$ and a
supercurve
$X$, we construct a vector bundle $\cV^r_X$ on $X$,
the fiber of which, is isomorphic to $V$. Moreover, the state-field
correspondence of $V$ canonically gives rise to (local) sections of these vector
bundles. We also define \emph{chiral algebras} on any supercurve $X$, and show
that the vector bundle $\cV^r_X$, corresponding to a SUSY vertex algebra,
carries the structure of a
chiral algebra. 
\end{abstract}
\maketitle
\tableofcontents
\section{Introduction}
\begin{nolabel}
	Vertex algebras were introduced about 20 years ago by Borcherds
\cite{borcherds}. They provide a rigorous definition of the chiral part of
2-dimensional conformal field theory, intensively studied by physicists. Since
then,
they have had important applications to string theory and conformal field
theory, and to mathematics, by providing tools to study the most interesting
representations of infinite dimensional Lie algebras. Since their appearance, they
have been extensively studied in many papers and books (for the latter we
refer to \cite{flm}, \cite{frenkelhuang}, \cite{kac:vertex}, \cite{huang}, \cite{frenkelzvi},
\cite{beilinsondrinfeld}).

Vertex algebras also appeared in algebraic geometry as \emph{Factorization Algebras}
on complex curves
\cite{beilinsondrinfeld}, \cite{frenkelzvi}. In the last five years, numerous
applications of this deep connection between factorization algebras and vertex
algebras have been exploited, notably in the study of the \emph{moduli spaces} (of
curves, vector bundles, principal bundles, etc) arising in algebraic
geometry. There are also connections between the theory of vertex
algebras and the 
 geometric Langlands conjecture \cite[ch.
17]{frenkelzvi}. Vertex algebras have also given new invariants of manifolds
\cite{kapranov}, \cite{malikov}  and applications to mirror symmetry \cite{borisov}.

Even though these approaches have been successful in formalizing 2-dimensional conformal
field theories, it has been known for some time to physicists, that in order to describe supersymmetric 
theories, similar objects should be defined on \emph{supercurves} instead
of simply curves (cf. \cite{dolgikh}, \cite{cohn}, \cite{banks}). With this
motivation, mathematicians have studied in detail the supergeometry of
manifolds, and in particular supercurves (cf. \cite{deligne2}, \cite{manin3},
\cite{manin2}, \cite{vaintrob1} among others).

The purpose of this article is to generalize the above objects to describe chiral
algebras over supercurves. To accomplish this, we will use 
the \emph{supersymmetric (SUSY) vertex algebras} defined in 
\cite{heluani3}. Roughly speaking, SUSY vertex algebras are vertex algebras with a \emph{state-field
correspondence} that includes the odd coordinates of a supercurve as formal
parameters, that is, to any vector $a$ in a SUSY vertex algebra, we associate
 a \emph{superfield}
\begin{equation*}
	\ys(a,z,\theta^1, \dots, \theta^N),
\end{equation*}
such that structural properties, similar to
those of ordinary vertex
algebras, hold.  

Given a SUSY vertex algebra $V$ and a supercurve $X$, we want to assign a vector bundle
$\cV$ over $X$ in such a way that the fiber at a point $x \in X$ is identified
with $V$. Moreover, we would like $\ys$ to canonically define sections of this
vector bundle (more precisely, its restricted dual). Here we find the first
difference with the classical theory, namely, supercurves come in different
flavors: general $1|n$ dimensional supercurves and superconformal $1|n$ supercurves. The latter are to the former what holomorphic curves are to compact
connected 2-manifolds. In 
\cite{heluani3}, two different versions of
SUSY vertex algebras are defined, one which localizes to give vector bundles on
a general $1|n$-dimensional supercurve (called $N_W = n$ SUSY vertex algebras) and
another which gives vector bundles on superconformal supercurves (called
$N_K=n$ SUSY vertex algebras).

There are several relations between these different SUSY vertex algebras. As a
basic example, let us consider the cases with low odd dimensions. Roughly
speaking, a general $N=1$ supercurve is the data of a curve $X$ and a line
bundle $\cL$ over it, sections of this line bundle are considered to be the
\emph{values of a coordinate in the odd direction}. Similarly, an (oriented) superconformal
$N=2$ supercurve consists of a curve $X$ and two line bundles $\cL$ and $\cH$
over it such that $\cL \otimes \cH$ is the canonical bundle $\omega$ of $X$.
It follows that an $N=1$ supercurve gives rise canonically to another $N=1$
supercurve (interchanging $\cL$ with $\omega\otimes \cL^{-1}$) and to a
superconformal $N=2$ supercurve (by taking $\cH = \omega \otimes \cL^{-1}$). On
the algebraic side, any (conformal) $N_W=1$ SUSY vertex algebra gives rise to a
(conformal) $N_K=2$ SUSY vertex
algebra (this corresponds to the isomorphism between the superconformal Lie
algebras $K(1|2)$ and $W(1|2)$) and both of them correspond to vertex algebras with
$N=2$ superconformal structure. It follows that any such vertex algebra gives
vector bundles in both $N=1$ supercurves and in the corresponding
superconformal $N=2$ supercurve. These three vector bundles are intimately
related as we will see in section \ref{sec:7}. 

As in the ordinary vertex algebra case, the vector bundles we construct (more precisely
quotients of them) are extensions of (powers of) the Berezinian bundle of $X$ (a super analog
of the canonical bundle). The algebraic properties of $V$ reflect in geometric properties of $\cV$ as in
the ordinary vertex algebra case. We obtain thus superprojective structures, affine
structures, global differential operators, etc. as splittings of these
extensions. In particular, the state-field correspondence itself gives such
splittings (locally).  
	\label{no:intro.1}
\end{nolabel}
\begin{nolabel}
	After constructing these vector bundles, it is natural to ask if they
	carry the structure of a \emph{chiral algebra} on a supercurve. It is shown that
	the usual definitions carry over to the super case with minor
	difficulties, and that the vector bundles obtained from $V$ are indeed
	chiral algebras. This allows us to define the coinvariants and conformal
	blocks of a SUSY vertex algebra in a coordinate independent way as in
	\cite{frenkelzvi}. 
	\label{no:intro.2}
\end{nolabel}
\begin{nolabel}
	The organization of this article is as follows: In section
	\ref{sec:intro} we recall some well known notions about vertex algebras
	and supercurves. We also summarize here some results on the structure theory of SUSY
	vertex algebras.

	In section \ref{sec:vector_bundles_section} we construct a vector bundle
	with a flat connection associated to an $N_W=n$ SUSY vertex algebra, over
	any $N=n$ supercurve. We also construct vector bundles associated to
	$N_K=n$ SUSY vertex algebras over oriented superconformal $N=n$
         supercurves.
	In this section we follow closely
	\cite{frenkelzvi} since most proofs are straightforward adaptations of the
	those in
	the ordinary vertex algebra case. The new phenomena can be found in section
	\ref{sec:7} where we analyze in detail examples with supersymmetry.
	
	  In section
	\ref{sub:group_aut_o} we define the groups $\mathrm{Aut} \cO$ of changes of
	coordinates and the $\mathrm{Aut} \cO$-torsor $\mathrm{Aut}_X$ for a 
	supercurve. In section \ref{sec:vector_bundles} we construct the vector bundles
	themselves and their sections. In particular we show that the
	state-field correspondence for a SUSY vertex algebra is a section of the
	dual of the corresponding vector bundle. In section \ref{sec:7} we
	compute explicitly some examples of vector bundles over supercurves of low odd
	dimension. 

	In section \ref{sec:chiral} we define chiral algebras over supercurves
	and we prove that the vector bundles constructed from SUSY vertex
	algebras are examples of chiral algebras. We also define
	the spaces of coinvariants in a coordinate independent way. 

	In appendix \ref{ap:gl1_1_rep} we give a brief description of a family of
	representations of the Lie algebra $\fg\fl(1|1)$ and their realizations as
	fibers of certain natural vector bundles over $N=1$ supercurves. 
	\label{no:intro_organization}
\end{nolabel}
\textbf{Acknowledgments.} The author would like to thank Victor G. Kac for reading
the manuscript, encouraging him, and many useful discussions. He would also like to
thank David Ben-Zvi, for very useful discussions.
\numberwithin{thm}{subsection}
\section{Preliminaries}\label{sec:intro}
\subsection{SUSY vertex algebras} \label{sub:super_recolection}
In this section we collect some results and examples of SUSY vertex algebras from
\cite{heluani3}. For the general theory of vertex algebras, the reader is referred
to the book \cite{kac:vertex} and \cite{kacdesole2} for an excellent exposition.
\begin{nolabel}
	Let $N$ be a non-negative integer. We will denote $Z=(z,\theta^1,\dots,\theta^N)$ where
	$z$ is an even indeterminate, $\theta^i$'s are odd anticommuting indeterminates
	commuting with $z$. For $J= (j_1, \dots, j_s)$ ordered subset of $\{1,\dots, N\}$, and
	$j \in \mathbb{Z}$, we
	will denote  
	\begin{equation*}
		\theta^J = \theta^{j_1} \dots \theta^{j_s}, \quad Z^{j|J} =  \theta^J z^j,
	\end{equation*}
	and we will denote by $N\setminus J$ the ordered complement of $J$ in $\{1, \dots,
	N\}$. For two disjoint subsets $I,\,J \subset \{1, \dots, N\}$ define
	$\sigma(I,J) = \pm 1$ by
	$\theta^I \theta^J = \sigma(I,J) \theta^{I \cup J}$,
	and $\sigma(J) = \sigma(J, N\setminus J)$. Finally, define $e_i = \{i\}$.

	Given a vector superspace $V$, we will denote by $V[ [Z]]$ (resp. $V( (Z))$) the space
	of formal power series (resp. formal Laurent series) in $Z$ with values in $V$, namely,
	formal sums of the form 
	\begin{equation*}
		\sum_{j \geq 0, J} Z^{j|J} v_{j|J} \quad \left( \text{ resp. } \sum_{j \geq
		N_0, J} Z^{j|J} v_{j|J} \right), 
	\end{equation*}
	where $N_0$ is some integer number and $v_{j|J} \in V$. Finally, we will denote by
	$V[Z, Z^{-1}]$ the space of Laurent polynomials with coefficients in $V$, namely,
	elements of $V( (Z))$ which are finite sums.
	
	Let $\cH_W$ (resp. $\cH_K$) be the associative superalgebra generated by
	an even element
	$T$ and $N$ odd elements $S^i$ subject to the relations
	\begin{equation*}
		{[}T, S^i] = 0, \quad [S^i, S^j] = 0 \text{ (resp. }2\delta_{i,j}T \text{)}.
	\end{equation*}
	
	An $N_W=N$ (resp $N_K=N$) SUSY vertex algebra $(V, \vac, Y)$ is the data of a
	$\cH_W$-module (resp. $\cH_K$-module)
	 $V$ (the space of states), an even vector $\vac \in V$ (the vacuum
	 vector) and a parity preserving $\mathbb{C}$-bilinear product with values in Laurent
	series over $V$:
	\begin{equation*}
		V \otimes V \rightarrow V( (Z)), \quad a \otimes b \mapsto Y(a,Z)b = \sum_{j
		\in \mathbb{Z}, J} Z^{-1-j|N\setminus J} a_{(j|J)}b,
	\end{equation*}
	subject to the following axioms $(a,\,b \in V)$:
	\begin{itemize}
		\item (vacuum axioms) $Y(a,Z)\vac|_{Z=0} = a$, $T\vac = S^i
			\vac = 0$, for $i = 1, \dots, N$, 
		\item (translation invariance) $[T, Y(a,z)] = \partial_z Y(a,Z)$, $[S^i,
			Y(a,Z)] = \partial_{\theta^i} Y(a,Z)$ (resp.
			$(\partial_{\theta^i} - 
			\theta^i \partial_z) Y(a,Z)$),
		\item (locality) $(z-w)^n [Y(a,Z), Y(b,W)] = 0$ for some $n \in \mathbb{Z}_+$.
	\end{itemize}
	\label{no:basic}
\end{nolabel}
\begin{rem}
	Note that when $N=0$, this definition agrees with the usual definition of vertex
	algebra as in \cite{kac:vertex}
	\label{rem:usual_defn}
\end{rem}
\begin{nolabel}
	Denote $\Lambda = (\lambda, \chi^1, \dots, \chi^N)$, where $\lambda$ is an even
	indeterminate and $\chi^i$'s are odd indeterminates, subject to the relations:
	\begin{equation*}
		{[}\lambda, \chi^i] = 0, \quad [\chi^i, \chi^j] = 0 \; (\text{ resp. } -
		2\delta_{i,j}\lambda),
	\end{equation*}
	and write $Z \Lambda = z \lambda + \sum_{i=1}^N \theta^i \chi^i$. 
	
	Let $V$ be an $N_W=N$ (resp. $N_K = N$) SUSY vertex algebra. For $a,\,b \in V$ we
	define
	\begin{equation}
		[a_\Lambda b] = \res_Z e^{Z \Lambda} Y(a,Z)b,
		\label{basic.2.2}
	\end{equation}
	where $\res_Z$ stands for the coefficient to the right of $Z^{-1|N}$. We note that the
	right hand side of (\ref{basic.2.2}) is a finite sum of monomials in $\Lambda$ times
	elements of $V$ (cf. \cite{heluani3}). This operation is called the $\Lambda$-bracket.
	As in the usual vertex algebra case, it encodes the singular part of the
	\emph{operator product expansion} (OPE) in $V$. 

	Define the \emph{normally ordered product} $: :$ as a $\mathbb{C}$-bilinear
	product on $V$:
	\begin{equation*}
		V \otimes V \rightarrow V, \qquad a \otimes b \mapsto :a b: =
		a_{(-1|N)}b.
	\end{equation*}
	The action of $\cH_W$ (resp. $\cH_K$) on $V$ is by derivations of both, the
	$\Lambda$-bracket, and the normally ordered product. As proved in
	\cite{heluani3}, these two operations encode all the structure of the SUSY
	vertex algebra $V$.
	\label{no:basic.2}
\end{nolabel}
\begin{ex} (Virasoro) This is an ordinary vertex algebra generated by one even
	field $L$ satisfying
	\begin{equation}
		{[}L_\lambda L] = (T + 2\lambda) L + \frac{\lambda^3}{12} c, 	
		\label{basic.ex.2.0}
	\end{equation}
	where $c \in \mathbb{C}$ is the \emph{central charge}. Expanding this field as
	\begin{equation*}
		L(z) = \sum_{n \in \mathbb{Z}} z^{-2-n} L_n,
		\label{basic.ex.2.0a}
	\end{equation*}
	we obtain that the operators $L_n$ satisfy the commutation relations of the Virasoro
	algebra of central charge $c$, namely:
	\begin{equation*}
		{[}L_m, L_n] = (m - n) L_{m+n} + \delta_{m, -n} \frac{m^3 - m}{12} c.
		\label{basic.ex.2.0b}
	\end{equation*}
	A field $L$ in an ordinary vertex algebra $V$ 
	satisfying (\ref{basic.ex.2.0}) will be called a \emph{Virasoro field} of central
	charge $c$. 
	\label{ex:virasoro_example}
\end{ex}
\begin{defn}
	Let $V$ be an ordinary vertex algebra with a Virasoro field $L$. We will say that a vector $a \in V$
	that satisfies $[L_\lambda a] = (T + \Delta \lambda) a + O(\lambda^2)$ has
	\emph{conformal weight} $\Delta$. If moreover, $a$ satisfies $[L_\lambda a] = (T +
	\Delta \lambda)a$, we will say that $a$ is \emph{primary}. 
	\label{defn:primary_def}
\end{defn}
\begin{ex} (Neveu-Schwarz) This vertex algebra is generated by a Virasoro field as in Example
	\ref{ex:virasoro_example} and an odd field $G$, primary of conformal weight $3/2$.
	The remaining $\lambda$-bracket is given by:
	\begin{equation*}
		[G_\lambda G] = 2 L + \frac{\lambda^2}{3} c.
	\end{equation*}
	If we expand the corresponding fields as: 
	\begin{equation*}
			L(z) = \sum_{n \in \mathbb{Z}} L_n z^{-2-n}, \quad
			G(z) = \sum_{n \in 1/2 + \mathbb{Z}} G_n z^{-3/2
			- n},
	\end{equation*}
	then the coefficients of such expansions satisfy the following
	commutation relations: 
	\begin{equation}
		\begin{gathered}
			{[}L_m, L_n] = (m-n) L_{m+n} + \delta_{m, -n}
			\frac{m^3 - m}{12} c, \\
			[G_m, L_n] = \left( m - \frac{n}{2} \right)
			G_{m+n}, \quad
			[G_m, G_n] = 2 L_{m+n} + \delta_{m,-n}
			\frac{m^2 - 1/4}{3} c.
		\end{gathered}
		\label{eq:n=1fouriercommutations}
	\end{equation}
	\label{ex:ns.example}
\end{ex}
\begin{rem}
	Let $V$ be a vertex algebra with an $N=1$ superconformal vector $\tau$ (cf.
	\cite[definition 5.9]{kac:vertex}). Namely, the Fourier modes of the fields
	\begin{equation*}
		\begin{aligned}
			G(z) &= Y(\tau,z) \sum_{n \in 1/2 + \mathbb{Z}} G_n
			z^{-n-3/2}, \\
			L(z) &= \frac{1}{2} Y(G_{-1/2} \tau,z) = \sum_{n \in
			\mathbb{Z}} L_n z^{-2-n},
		\end{aligned}
	\end{equation*}
	satisfy the relations (\ref{eq:n=1fouriercommutations}) of a Neveu-Schwarz
	algebra for some $c \in \mathbb{C}$, $L_{-1} (= G_{-1/2}^2 ) = T$ and the
	operator $L_0$ is diagonalizable with eigenvalues bounded below. Then $V$
	carries a structure of an $N_K=1$ SUSY vertex algebra with $S = G_{-1/2}$ and the
	superfields are defined as:
	\begin{equation}
		\ys(a,z,\theta) = Y(a,z) + \theta Y(G_{-1/2}a,z).
		\label{eq:nimepreguntes}
	\end{equation}
	\label{rem:quecosache}
\end{rem}
\begin{ex} \cite[ex. 5.9a]{kac:vertex} \label{ex:2.16}
	Let $B_1$ be the vertex algebra generated by an even vector (free boson)
	 $\alpha$ and an odd vector (free fermion) $\varphi$, namely:
	\begin{equation*}
			{[}\alpha_\lambda \alpha] = \lambda, \quad
			{[}\varphi_\lambda \varphi] = 1, \quad
			{[}\alpha_\lambda \varphi] = 0.
		\label{}
	\end{equation*}
	Then $B_1$ is a (simple) vertex algebra with a family of $N=1$
	superconformal vectors 
	\begin{equation*}
		\tau = (\alpha_{(-1)} \varphi_{(-1)} + m
		\varphi_{(-2)})\vac, \qquad m \in \mathbb{C},
		\label{}
	\end{equation*}
	of central charge $c = \tfrac{3}{2} - 3 m^2$. 
	\label{ex:boson_fermion_ex}
\end{ex}
\begin{ex} \cite[Thm 5.10]{kac:vertex} \label{ex:n=2}
	The $N=2$ vertex algebra is generated by a Virasoro field $L$ of central
	charge $c$, an
	even field $J$, primary of conformal weight $1$, and two odd fields
	$G^\pm$, primary of conformal weight $3/2$. The remaining
	$\lambda$-brackets are:
	\begin{equation*}
	\begin{gathered}
			{[}J_\lambda J] = \frac{c}{3} \lambda, \quad 
			{[}G^\pm_\lambda G^\pm] = 0, \quad
			{[}J_\lambda G^\pm] = \pm G^\pm,\\
			{[}G^+_\lambda G^-] = L + \frac{1}{2} \partial J
			+ \lambda J + \frac{c}{6} \lambda^2.
		\label{}
	\end{gathered}
	\end{equation*}
	This vertex algebra contains an $N=1$ superconformal vector:
	\begin{equation*}
		\tau = G^+_{(-1)} \vac + G^-_{(-1)} \vac.
		\label{}
	\end{equation*}
	Also, this vertex algebra admits a $\mathbb{Z}/2 \mathbb{Z} \times
	\mathbb{C}^*$ family of automorphisms. The generator of
	$\mathbb{Z}/2 \mathbb{Z}$ is given by $L \mapsto L$, $J \mapsto
	-J$ and $G^\pm \mapsto G^\mp$. The $\mathbb{C}^*$ family is given
	by $G^+ \mapsto \mu G^+$ and $G^- \mapsto  \mu^{-1} G^-$.
	Applying these automorphisms, we get a family of $N=1$
	superconformal
	structures. 

	By expanding the corresponding fields:
	\begin{equation*}
		\begin{gathered}
			L(z) = \sum_{n \in \mathbb{Z}} L_n z^{-2-n}, \\
			G^\pm(z) = \sum_{n \in 1/2 + \mathbb{Z}} G^\pm_n
			z^{-3/2 
			-n}, \quad
			J(z) = \sum_{n \in \mathbb{Z}} J_n z^{-1-n},
		\end{gathered}
		\label{}
	\end{equation*}
	we get the commutation relations of the Virasoro operators $L_n$, and the following
	remaining commutation relations
	\begin{equation*}
		\begin{gathered}
			\begin{aligned}
			{[}J_m, J_n] &= \frac{m}{3} \delta_{m,-n} c, &  
			{[}J_m, G^\pm_n] &= \pm G^\pm_{m+n}, \\
			[G^\pm_m, L_n] &= \left( m -\frac{n}{2}  \right)
			G^\pm_{m+n}, &
			[L_m, J_n] &= - n J_{m+n}, 
		\end{aligned} \\
			{[}G^+_m, G^-_n] = L_{m+n} +
			\frac{m-n}{2} J_{m+n} + \frac{c}{6} \left( m^2 -
			\frac{1}{4}
			\right) \delta_{m,-n}.			
		\end{gathered}
		\label{}
	\end{equation*}
	
	Sometimes it is convenient to introduce a different set of generating
	fields for this vertex algebra. We define $\tilde{L}(z) = L(z) - 1/2 \partial_z
	J(z)$. This is a Virasoro field with central charge zero,
	namely:
	\begin{equation*}
		{[}\tilde{L}_\lambda \tilde{L}] = (\partial + 2\lambda)
		\tilde{L}.
		\label{}
	\end{equation*}
	With respect to this Virasoro element, $G^+$ is
	primary of conformal weight $2$ and $G^-$ is primary of conformal
	weight $1$; $J$ has conformal weight $1$ but is no longer
	a primary field. To summarize the commutation relations, we write:
	\begin{equation}
		\begin{gathered}
			Q(z) = G^+(z) = \sum_{n \in \mathbb{Z}} Q_n
			z^{-2-n}, \\
			H(z) = G^-(z) = \sum_{n \in \mathbb{Z}} H_n
			z^{-1-n}, \quad
			\tilde{L}(z) = \sum_{n \in \mathbb{Z}} T_n
			z^{-2-n}. 
		\end{gathered}
		\label{eq:n=2expansion}
	\end{equation}
		The corresponding $\lambda$-brackets of these fields are
		given by:
	\begin{equation}
		\begin{gathered}
		\begin{aligned}
			{[}\tilde{L}_\lambda \tilde{L}] &= (\partial + 2\lambda)
			\tilde{L}, &
			{[}\tilde{L}_\lambda J] &= (\partial + \lambda)J -
			\frac{\lambda^2}{6}c, \\
			{[}\tilde{L}_\lambda Q] &= (\partial + 2\lambda) Q, &
			{[}\tilde{L}_\lambda H] &= (\partial + \lambda) H,
		\end{aligned} \\
			[H_\lambda Q] = \tilde{L} - \lambda J +
			\frac{c}{6} \lambda^2. 
		\end{gathered}
		\label{eq:agregado.w.n=2}
	\end{equation}
	The commutation relations of the coefficients in
	(\ref{eq:n=2expansion}) are:
	\begin{equation}
	\begin{aligned}
			{[}T_m, T_n] &= (m-n) T_{m+n}, & [Q_m,Q_n] &= [H_m, H_n] = 0, \\ 
			{[} T_m, H_n] &= - n H_{m+n}, &
			[T_m , J_n] &= -n J_{m+n} - m (m+1)
			\frac{c}{12}\delta_{m, -n},\\  
			{[} T_m, Q_n] &= (m-n) Q_{m+n}, &
			[H_m, Q_n] &= T_{m+n} - m J_{m+n} + m
			(m-1)\frac{c}{6} \delta_{m, -n}\\
	\end{aligned}
	\label{eq:n=2commutations_new_fields}
	\end{equation}

	Finally, defining $G^{(1)} = G^+ + G^{-}$ and $G^{(2)} = i (G^+ - G^-)$, we
	obtain another set of generators for this vertex algebra.  We note
	that, with respect to $L$, the fields $G^{(i)}$ are primary of
	conformal weight $3/2$, and $J$ is primary of conformal weight $1$. The other
	commutation relations between the generating fields $L,\,J,\, G^{(i)}$ are 
	\begin{equation*}
		\begin{aligned}
			{[}{G^{(i)}}_\lambda G^{(i)}] &= 2 L +
			\frac{c \lambda^2}{3}, &
			[{G^{(1)}}_\lambda G^{(2)}] &= -i \left( \partial
			+ 2\lambda
			\right) J, \\
			[J_\lambda G^{(1)}] &= - i G^{(2)}, &
			[J_\lambda G^{(2)}] &= i G^{(1)},
		\end{aligned}
	\end{equation*}
	or, equivalently:
	\begin{equation}
		\begin{gathered}
		\begin{aligned}
			{[}G^{(i)}_m, G^{(i)}_n] &= 2 L_{m+n} + \left( m^2 -
		\frac{1}{4}
		\right) \frac{c}{3} \delta_{m, -n}, \\
		{[}G^{(1)}_m, G^{(2)}_n] &= i \left( n - m \right) J_{m+n},
	\end{aligned} \\
		[J_{m}, G^{(1)}_n] = -i G^{(2)}_{m+n}, \quad
		[J_m, G^{(2)}_n] = i G^{(1)}_{m+n}.
		\end{gathered}
	\label{eq:n=2superconf}
	\end{equation}
	\label{ex:n=2_defn_ex}
\end{ex}
\begin{rem}
	As in the $N=1$ case, given an $N=2$ superconformal vertex algebra, namely a vertex
	algebra with a vector $j$ and two operators $S^1, S^2$ satisfying
	\begin{equation*}
		{[T}, S^i] = 0, \qquad [S^i, S^j] = 2 \delta_{i,j} T,
	\end{equation*}
	and such that the corresponding fields:
	\begin{equation*}
		\begin{aligned}
			J(z) &= -i Y(j, z), \qquad
			L(z) = \frac{1}{2} Y(S^2 S^1 j, z),\\
			G^{(1)}(z) &\equiv G^+(z) + G^-(z) = - Y(S^2 j, z), \\
			G^{(2)}(z) & \equiv i \left(G^+(z) - G^-(z)\right) = Y(S^1 j,
			z),\\
		\end{aligned}
	\end{equation*}
	satisfy the $\lambda$-brackets of Example
	\ref{ex:n=2_defn_ex}, $L_{-1} = T$, $G^{(i)}_{-1/2} = S^i$, and $L_0$ is diagonalizable with eigenvalues bounded
	below, we obtain an $N_K = 2$ SUSY vertex algebra by letting 
	$\ys(a,Z) = Y(a,z) + \theta^1 Y(S^1 a,z) + \theta^2 Y(S^2a, z) +  \theta^2
		\theta^1 Y(S^1 S^2 a, z)$.
	
		Similarly, given a vertex algebra with two vectors $j, h$ and an
	odd operator $S$ such that $[T,S] = 0$, $S^2 = 0$ and the
	associated fields:
	\begin{equation*}
		\begin{aligned}
			J(z) &= - Y(j,z), & H(z) &= Y(h, z), \\
			Q(z) &= Y(S j, z), & \tilde{L}(z) &= Y(S h, z) -
			 \partial_z J(z),
		\end{aligned}
	\end{equation*}
	satisfy the commutation relations
	(\ref{eq:agregado.w.n=2}), $T_{-1} = T$, $Q_{-1} = S$, $J_0$ is
	diagonalizable, and $T_0$ is diagonalizable with eigenvalues
	bounded below,
	 we obtain an $N_W = 1$ SUSY
	vertex algebra
	by letting:
	$\ys(a,Z) = Y(a,z) + \theta Y(S a, z)$.
	\label{rem:n=2.construction}
\end{rem}
\begin{ex} \cite[ex. 5.9d]{kac:vertex}\label{ex:2.23} 
	Consider the vertex algebra generated by a pair of free charged
	bosons $\alpha^\pm$ and a pair of free charged fermions
	$\varphi^\pm$ where the only non-trivial commutation relations are:
	\begin{equation*}
			{[}{\alpha^\pm}_\lambda \alpha^\mp] = \lambda,  \qquad
			[{\varphi^\pm}_\lambda \varphi^\mp] = 1.
		\label{}
	\end{equation*}
	This vertex algebra contains the following family of $N=2$ vertex
	subalgebras:
	\begin{equation*}
		\begin{gathered}
			G^\pm = :\alpha^\pm \varphi^\pm: \pm m \partial
			\varphi^\pm, \quad  J = :\varphi^+ \varphi^-: -m
			(\alpha^+ + \alpha^-), \qquad m \in \mathbb{C},\\
			L = :\alpha^+ \alpha^-: + \frac{1}{2} : \partial
			\varphi^+ \varphi^-: +  \frac{1}{2} :\partial \varphi^-
			\varphi^+: - \frac{m}{2} \partial (\alpha^+ -
			\alpha^-).
		\end{gathered}
		\label{}
	\end{equation*}
	The vector $j = i J_{-1}\vac$ provides this vertex algebra
	with the structure of an
	$N_K=2$ SUSY vertex algebra, by letting $T=L_{-1}$ and $S^i
	= G^{(i)}_{-1/2}$, where $G^{(1)} = G^+ + G^-$ and $G^{(2)} = i( G^+ - G^-)$. 
\end{ex}
\begin{ex} ($W_N$ series) Now we define an $N_W = N$ SUSY vertex algebra for each non-negative
	integer $N$. 
	When $N=0$, $W_0$ is the Virasoro vertex algebra of central charge $c$.
	
	$W_1$ is the $N_W=1$ SUSY vertex algebra generated\footnote{See \cite{heluani3} for the
	definition of \emph{generating fields} of a SUSY vertex algebra.} by an odd superfield
	$L$ and an even
	superfield $G$ satisfying:
	\begin{equation}
		\begin{aligned}
			{[}L_\Lambda L] &= (T + 2\lambda) L, \qquad [Q_\Lambda Q] = SQ +
			\frac{\lambda \chi}{3} c, \\
			[L_\Lambda Q] &= (T+\lambda)Q - \chi L + \frac{\lambda^2}{6} c,
		\end{aligned}
		\label{basic.ex.2.1}
	\end{equation}
	where $c \in \mathbb{C}$ is the \emph{central charge}. Expanding these
	superfields as:
	\begin{equation*}
		Q(Z) = - J(z) + \theta G^+(z), \qquad L(Z) = G^-(z) + \theta \left( L(z) +
		\frac{1}{2} \partial_z J(z)
		\right),
		\label{basic.ex.2.2}
	\end{equation*}
	we find that the fields $L, J, G^\pm$ generate an $N=2$ vertex algebra of central
	charge $c$ as in example \ref{ex:n=2_defn_ex}.  

	$W_2$ is the $N_W=2$ SUSY vertex algebra generated by an even superfield $L$ and two
	odd superfields $Q^1$ and $Q^2$ satisfying:
	\begin{equation}
		\begin{aligned}
		{[}L_\Lambda L] &= (T+2\lambda) L, & [{Q^i}_\Lambda Q^i] &= S^i Q^i, \\
		{[Q^1}_\Lambda Q^2] &= (S^1 + \chi^1) Q^2  - \chi^2 Q^1 +
		\frac{\lambda}{6}c,
		& [L_\Lambda Q^i] &= (T + \lambda) Q^i + \chi^i L.
	\end{aligned}
		\label{eq:basic.ex.2.4}
	\end{equation}

	Finally, for $N \geq 3$ we let $W_N$ be the $N_W = N$ SUSY vertex algebra generated by
	a superfield $L$ of parity $N \mod 2$ and $N$ superfields $Q^i$, $i=1, \dots, N$ of
	parity $N + 1 \mod 2$, satisfying:
	\begin{equation}
		\begin{aligned}
			{[L}_\Lambda L] &= (T + 2\lambda) L, \qquad {[Q^i}_\Lambda Q^j] = (S^i
			+ \chi^i) Q^j - \chi^j Q^i, \\
			[L_\Lambda Q^i] &= (T+ \lambda) Q^i + (-1)^N \chi^i L.
		\end{aligned}
		\label{eq:basic.ex.2.5}
	\end{equation}

	It is proved in \cite{heluani3} that the Lie superalgebra $W(1|N)$  of derivations on
	$\mathbb{C}[Z, Z^{-1}]$ acts on $W_N$. We let $W(1|N)_-$ (resp. $W(1|N)_<$) be the Lie
	subalgebras of regular vector fields (resp. regular vector fields vanishing at the
	origin). 
\end{ex}
\begin{defn}
	An $N_W=N$ SUSY vertex algebra $V$ is called \emph{conformal}  if there exists
	$N+1$ vectors $\nu, \tau^1,\dots, \tau^N$ in $V$ such that their associated superfields
	$L(Z) = \ys(\nu, Z)$ and $Q^i(Z) = \ys(\tau^i, Z)$ satisfy (\ref{eq:basic.ex.2.5}) for $N
	\geq 3$, (\ref{eq:basic.ex.2.4}) for $N=2$, (\ref{basic.ex.2.1}) for $N=1$ or
	(\ref{basic.ex.2.0}) for $N=0$, 
        and moreover:
	\begin{itemize}
		\item $\nu_{(0|0)} = T$, 
		 $\tau^i_{(0|0)} =  S^i$,
		\item The operator $\nu_{(1|0)}$ acts diagonally with eigenvalues bounded 
			below and with finite dimensional eigenspaces.
	\end{itemize}
	
	If moreover, 
	 the action of $W(1|N)_<$ on $V$ can be
	 exponentiated to the group of automorphisms of the $1|N$ dimensional superdisk
	 $D^{1|N}$ (see \ref{no:disks} for a definition), we will say that
	$V$ is \emph{strongly conformal}. 
	This amounts to the following extra condition:
	\begin{itemize}
		\item The operators $\nu_{(1|0)}$ and $\sum_{i=1}^N \sigma(e_i)
			\tau^i_{(0|e_i)}$ have integer
			eigenvalues.
	\end{itemize}
	\label{defn:w.conformal.definition}
\end{defn}
\begin{ex}[$K_N$ series] 
	For $N \leq 3$, let $K_N$ be the $N_K=N$ SUSY vertex algebra generated by one
	superfield $G$ of parity $N \mod 2$ 
	satisfying:
	\begin{equation}
		[G_\Lambda G] = \left( 2 T + (4-N) \lambda + \sum_{i = 1}^N
		\chi^i S^i \right) G + \frac{\lambda^{3-N}\chi^N}{3} c,
		\label{eq:ex.k_n.1}
	\end{equation}
	where $c \in \mathbb{C}$ is called the \emph{central charge}. Let $K_4$ be the $N_K=4$
	SUSY vertex algebra generated by an even superfield $G$ satisfying:
	\begin{equation}
		{[G}_\Lambda G] = \left( 2 T + \sum_{i=1}^4 S^i \chi^i \right) G + \lambda c.
		\label{eq:ex.k_n.2.a}
	\end{equation}

	In the case $N=1$, if we expand the corresponding superfield 
	 as
	\begin{equation*}
		G(z,\theta) = G(z) + 2 \theta L(z),
	\end{equation*}
	we find that the fields $G(z)$ and $L(z)$ generate a Neveu-Schwarz vertex
	algebra of central charge $c$ as in Example \ref{ex:ns.example}. 
	
	When $N=2$ expanding the corresponding superfield as:
	\begin{equation*}
		G(z, \theta^1,\theta^2) = \sqrt{-1} J(z) +\theta^1
		G^{(2)}(z) - \theta^2  G^{(1)}(z) + 2 \theta^1 \theta^2 L(z),
	\end{equation*}
	where $G^{(1)} = G^+ + G^-$ and $G^{(2)} = i (G^+ - G^-)$, we find that the corresponding fields $J, L, G^\pm$ satisfy the
	commutation relations of the $N=2$ vertex algebra as in Example
	\ref{ex:n=2_defn_ex}.

	It was proved in \cite{heluani3} that the Lie superalgebra $K(1|N)$ of vector fields on the
	$1|N$-dimensional superdisk $D^{1|N}$ preserving the differential 1-form 
	\begin{equation*}
		\omega = dz + \sum_{i=1}^N \theta^i d\theta^i,
		\label{ex.k.1.agrego.mas}
	\end{equation*}
	up to multiplication by a function, acts on $K_N$. We let $K(1|N)_-$ (resp. $K(1|N)_<$)
	be the Lie subalgebra of regular vector fields (resp. regular vector fields vanishing
	at the origin). 
	\label{ex:k.n.series}
\end{ex}
\begin{defn}
	Let $N \leq 4$, an $N_K=N$ SUSY vertex algebra $V$ is called \emph{conformal} if there exists
	a vector $\tau \in V$ (called the conformal vector) such that the corresponding field
	$G(Z) = \ys(\tau, Z)$ satisfies (\ref{eq:ex.k_n.1}) for $N \leq 3$ or
	(\ref{eq:ex.k_n.2.a}) for $N=4$, 
	and moreover
	\begin{itemize}
		\item $\tau_{(0|0)} = 2 T$, 
		 $\tau_{(0|e_i)} = \sigma(N\setminus e_i, e_i) S^i$,
		\item the operator $\tau_{(1|0)}$ acts diagonally with eigenvalues bounded
			 below and finite dimensional eigenspaces.
	\end{itemize}

	If moreover,  the
	representation of $K(1|N)_<$ can be exponentiated to the group of
	automorphisms of the disk $D^{1|N}$ preserving the differential form
	$\omega$ 
	up to multiplication by a function, 
	we will say that $V$ is \emph{strongly conformal}. 
This amounts to the extra condition \begin{itemize} \item the operator $\tau_{(1|0)}$ has
			integer eigenvalues, and if $N=2$, the operator
			$\sqrt{-1} \tau_{(0|N)}$ has integer eigenvalues. \end{itemize}

	\label{defn:k.conf.vertex.def}
\end{defn}
\begin{ex} (boson-fermion system) Let $B_1$ be the strongly conformal $N_K=1$ SUSY vertex algebra
	generated  by one odd superfield $\Psi$ satisfying:
	\begin{equation*}
		[\Psi_\Lambda \Psi] = \chi.
		\label{b_1.def}
	\end{equation*}
	Expanding the corresponding superfield as:
	\begin{equation*}
		\Psi(Z) = \varphi(z) + \theta \alpha(z),
		\label{b_1.expansion}
	\end{equation*}
	we find that the ordinary fields $\varphi$ and $\alpha$ generate the well known
	\emph{boson-fermion system} as in Example \ref{ex:boson_fermion_ex}.
	\label{ex:boson_fermion_susyk}
\end{ex}
\begin{nolabel} Now we summarize some basic results in the structure theory of SUSY
	vertex algebras.
The proofs can be found in \cite{heluani3}.  We will denote $\nabla = (T, S^1,
\dots, S^N)$, and for each $(j|J)$ as
above, we define:
\begin{equation*}
	\nabla^{j|J} = T^j S^J, \qquad \nabla^{(j|J)} =
	\frac{(-1)^{\frac{J(J+1)}{2}}}{j!} \nabla^{j|J}.
\end{equation*}
Let $Z = (z,\theta^1, \dots, \theta^N)$ be as before, and $W = (w, \zeta^1, \dots,
\zeta^N)$ be such that $w$ is an even indeterminate, commuting with $z$,
$\theta^i$'s and $\zeta^i$'s, and $\zeta^i$ are odd anticommuting indeterminates,
commuting with $z$, $w$ and anticommuting with $\theta^i$'s. In the $N_W=N$
 case we will write:
\begin{equation}
	Z - W = (z-w, \theta^1 - \zeta^1, \dots, \theta^N - \zeta^N),
	\label{eq:notation.bas.2}
\end{equation}
and in the $N_K=N$ case:
\begin{equation}
	Z-W = \left(z-w-\sum_{i=1}^N \theta^i \zeta^i, \theta^1 - \zeta^1, \dots,
	\theta^N - \zeta^N \right).
	\label{eq:notation.bas.3}
\end{equation}
We define the formal super delta-function to be:
\begin{equation*}
	\delta(Z,W) = \left( i_{z,w} - i_{w,z} \right) (Z-W)^{-1|N},
\end{equation*}
where $i_{z,w}$ denotes the \emph{expansion in the domain} $|z| > |w|$. We note
that this definition is independent of the definition of $Z-W$ as
(\ref{eq:notation.bas.2}) or (\ref{eq:notation.bas.3}). 

Put $Z \nabla = z T + \sum_i \theta^i S^i$. Similarly, in the $N_W=N$ (resp.
$N_K=N$) case, we put $D_W = (\partial_w, \partial_{\zeta^i})$ (resp. $D_W =
(\partial_w, \partial_{\zeta^i} + \zeta^i \partial_w)$). 
\end{nolabel}
\begin{prop}
	Let $V$ be a $N_W=N$ or $N_K = N$ SUSY vertex algebra. Then 
	\begin{multline*}
		[a_{(n|I)}, \ys(b,W)] = \sum_{(j|J), j\geq 0}
		(-1)^{JN + IN + IJ}  \sigma(J) \sigma(I) \times \\ \times
		 \bigl(D_W^{(j|J)} W^{n|I}\bigr)
		\ys\left( a_{(j|J)}b,W \right). 
	\end{multline*}
	If, moreover, $n \geq 0$, this becomes:
	\begin{equation*}
		[a_{(n|I)}, \ys(b,W)] = \ys(e^{-W\nabla}a_{(n|I)}e^{W\nabla}b,W).
	\end{equation*}	
	\label{prop:commutator_w.1}
\end{prop}
\begin{thm}
	In an $N_W=N$ (resp. $N_K = N$)  SUSY vertex algebra the following identities hold (see
	\cite{heluani3} for the
	definition of $(j|J)$-th product of superfields)
	\begin{enumerate}
		\item $\ys(a_{(j|J)}b,Z) = \sigma(J)
			\ys(a,Z)_{(j|J)}\ys(b,Z)$ (the $(j|J)$-th
			product identity), 
		\item $\ys(a_{(-1|N)}b,Z) = :\ys(a,Z)\ys(b,Z):$,
		\item $\ys(Ta, Z) = \partial_z \ys(a,Z)$,
		\item $\ys(S^ia, Z) = \partial_{\theta^i}\ys(a,Z)$ (resp.
			$(\partial_{\theta^i} + \theta^i \partial_z) Y(a,Z)$),
		\item we have the following OPE formula:
			\begin{equation}
				\begin{aligned}
				{[}\ys(a,Z), \ys(b,W)] &= \sum_{(j,J), j \geq 0}
				\sigma(J)
				(D_W^{(j|J)} \delta(Z,W) )
				\ys(a_{(j|J)}b,W) \\
				&= \sum_{(j|J), j \geq 0} (i_{z,w} - i_{w,z})
				(Z-W)^{-1-j|N\setminus
				J} \ys(a_{(j|J)} b, W),
			\end{aligned}
				\label{eq:ope_w.1}
			\end{equation}
			where the sum is finite and the operator $i_{z,w}$ denotes
			the expansion in the domain $|z| > |w|$.
	\end{enumerate}
	\label{thm:structure.2}
\end{thm}
\begin{lem} The following identity is true (note that $D_W$ and $Z-W$ here have
	different meanings in the $N_W=N$ and $N_K=N$ case):
	\begin{equation}
		D^{(j|J)}_W \delta(Z,W) =  \sigma(J, N \setminus
	J) (i_{z,w} - i_{w,z}) (Z-W)^{-1-j|N\setminus J}.
		\label{eq:k.deriv.1}
	\end{equation}
	\label{lem:k.deriv}
\end{lem}
\begin{thm}[Skew-symmetry]
	In a SUSY vertex algebra the following identity, called
	\emph{skew-symmetry}, holds
	\begin{equation}
		\ys(a,Z)b = (-1)^{ab}e^{Z \nabla} \ys(b,-Z)a
		\label{eq:skew-symmetry.w.1}
	\end{equation}
	\label{thm:skew-symmetry.w}
\end{thm}
\begin{thm}[Cousin property]
	For any SUSY vertex algebra $V$ and vectors $a, b, c \in V$,
	the three expressions:
	\begin{equation*}
		\begin{aligned}
			\ys(a,Z) \ys(b, W)c & \in V( (Z) ) ( ( W) ) \\
			(-1)^{ab} \ys(b,W) \ys(a, Z) c & \in V( (W)) ( (Z)) \\
			\ys\left( \ys(a,Z-W)b,W \right)c &\in V( (W) ) ( (Z-W))
		\end{aligned}
	\end{equation*}
	are the expansions, in the domains $|z| > |w|$, $|w| > |z|$ and $|w| > |w-z|$
	respectively, of the same element of 
	$V [ [ Z, W] ] [ z^{-1}, w^{-1}, (z-w)^{-1}]$.
	\label{cor:cousin.w}
\end{thm}
\begin{nolabel} Let $V$ be an $N_W = N$ (resp. $N_K=N$) SUSY vertex algebra. 
	It was proved in \cite{heluani3} that $\lie(V) = \tilde{V}/\tilde{\nabla}
	\tilde{V}$ is naturally a Lie algebra\footnote{Here we need to change the
	parity of $\tilde{V}/\tilde{\nabla}\tilde{V}$ if $N$ is odd, see
	\cite{heluani3} for details.}, where $\tilde{V} = V
	\otimes_{\mathbb{C}}
	\mathbb{C}[Z, Z^{-1}]$ and $\tilde{\nabla} \tilde{V}$ is the space spanned by vectors
	of the form:
	\begin{equation}
		\begin{split}
			&T a \otimes f(Z) + a \otimes \partial_z f(Z), \\
			&S^i a \otimes f(Z) + (-1)^{aN} a \otimes \partial_{\theta^i}
			f(Z), \\ &\text{ (resp.} S^i a \otimes f(Z) + (-1)^{aN} a
			\otimes (\partial_{\theta^i} + \theta^i \partial_z)
			f(Z) \text{)},
		\end{split}
		\label{eq:acabodeagregarte}
	\end{equation}
	for $a \in V$, $f(Z) \in \mathbb{C}[Z, Z^{-1}]$. 
	Let $\varphi : \lie (V) \rightarrow \End(V)$ be the linear map
	defined by 
	\begin{equation}
		a_{<n|I>} = a \otimes Z^{n|I} \mapsto (-1)^{aI}\sigma(I)
		a_{(n|I)}, \qquad a \in V.
		\label{eq:nuevamentealvicio}
	\end{equation}
	Similarly, we construct $V \otimes_{\mathbb{C}} \mathbb{C}( (Z))$ and consider its
	quotient $\lie'(V)$ by the vector space generated by vectors of the form
	(\ref{eq:acabodeagregarte}). Then (\ref{eq:nuevamentealvicio}) defines a map $\varphi' : \lie'(V) \rightarrow
	\End(V)$. 
\end{nolabel}
\begin{thm}
	The maps $\varphi$, and $\varphi'$ are Lie algebra homomorphisms.
	\label{thm:lie_fourier.vertex.1}
\end{thm}
\subsection{Supercurves}\label{sub:super_curves_recolection}
\begin{nolabel}
For a general introduction to the theory of supermanifolds and superschemes,
the reader should refer to \cite{manin2}. We will
 follow \cite{rabin1} for the theory of \emph{supercurves} over a
Grassmann algebra $\Lambda$. The deformation theory of superspaces and
sheaves over them can be found in \cite{vaintrob1}. The relations between
superconformal Lie algebras and the moduli spaces of supercurves was
stated in \cite{vaintrob2}. The reader may also find useful the notes
\cite{deligne2}.
\end{nolabel}
\begin{defn}
	A \emph{superspace} is a locally ringed space $(X, \cO_X)$ where
	$X$ is a topological space and $\cO_X$ is a
	sheaf of supercommutative rings. A \emph{morphism} of superspaces 
	is a graded morphism of locally ringed spaces. We will use $X$ to denote such
	a superspace when no confusion should arise. A
	\emph{superscheme} is a superspace such that $(X, \cO_{X,\bar{0}})$
	is a scheme, where from now on $\cO_{X,i}$ denotes the i-th
	graded part of $\cO_X$, $i = \bar{0}, \bar{1}$.
	\label{defn:super_scheme}
\end{defn}
\begin{nolabel}
	Given a superspace $(X,\cO_X)$ define $\cJ = \cO_{X,\bar{1}} +
	\cO_{X,\bar{1}}^2$. $\cJ$ is a sheaf of ideals in $(X,
	\cO_X)$, the corresponding subspace $(X, \cO_X/\cJ)$ will be
	denoted $(X_{\mathrm{rd}}, \cO_{X_{\mathrm{rd}}})$. 
\end{nolabel}
\begin{ex} \label{ex:specR}
	Let $R$ be a supercommutative ring, and let $J = R_1 + R_1^2$ be
	the ideal generated by $R_{\bar{1}}$ as above, then $(\Spec R, R)$ is a superscheme. Note
	that as topological spaces $\Spec R = \Spec R/J$ since every
	element in $J$ is nilpotent 
	(we consider only homogeneous ideals with respect to the $\mathbb{Z}/2
	\mathbb{Z}$-grading). 
\end{ex}
\begin{defn}[cf. \cite{manin2}]
	A \emph{supermanifold} is a superspace $(X,\cO_X)$ such that
	for every point $x \in X$ there exists an open neighborhood $U$
	of $x$ and a locally free sheaf $\cE$ of
	${\cO_{X_{\mathrm{rd}}}}|_U$-modules, of (purely odd) rank $0|q$ such
	that $(U, {\cO_X}|_U)$ is isomorphic to $(U_{\mathrm{rd}},
	S_{\cO_{X_{\mathrm{rd}}}}(\cE)|_U)$. Here $S(\cE)$ denotes the
	symmetric algebra of a  vector bundle.
\end{defn}
\begin{nolabel}
	An open sub-supermanifold of $(X,\cO_X)$ consist of an open
	subset $U \subset X$ and the restriction of the structure sheaf,
	namely $(U, \cO_X|_U)$. 
\end{nolabel}
\begin{nolabel}
	In the analytic setting, the situation is easier to describe. 
	The supermanifold $\mathbb{C}^{p|q}$ is the topological space
	$\mathbb{C}^p$ endowed with the sheaf of supercommutative algebras
	$\cO[\theta_1,\dots , \theta_q]$ where $\cO$ is
	the sheaf of germs of holomorphic functions on $\mathbb{C}^p$ and
	$\theta_i$ are odd anticommuting variables. A
	complex supermanifold is a topological space $|X|$ with a sheaf of
	supercommutative algebras $\cO_X$ locally isomorphic to
	$\mathbb{C}^{p|q}$. Morphisms of supermanifolds are continuous
	maps $\sigma:|X| \rightarrow |Y|$ together with morphisms of
	sheaves $\sigma^\sharp: \sigma^* \cO_Y \rightarrow \cO_X$. 
	
	Let $\Lambda = \mathbb{C}[\alpha_1, \dots, \alpha_n]$ be a
	Grassmann algebra. The $0|n$-dimensional superscheme $\Spec \Lambda$ has
	as underlying topological space a single point. We will work in
	the category of superschemes over $\Lambda$, namely super
	schemes $S$ together with a structure morphism $S \rightarrow
	\Spec \Lambda$. In the case when $S$ is a proper, smooth of
	relative dimension $1|q$ super-scheme, we say that $S$ is a $N=q$
	\emph{supercurve} (over $\Lambda$).
\end{nolabel}
\begin{defn}
	More explicitly (cf. \cite{rabin1}), a smooth compact connected
	complex supercurve over $\Lambda$ of dimension $1|N$ is a pair
	$(X, \cO_X)$, where $X$ is a topological space and $\cO_X$ is a
	sheaf of supercommutative $\Lambda$-algebras over $X$, 
	such that:
	\begin{enumerate}
		\item $(X, \cO_X^{\mathrm{red}})$ is a smooth compact
			connected algebraic curve. Here
			$\cO_X^{\mathrm{red}}$ is the reduced sheaf of
			$\mathbb{C}$-algebras on $X$ obtained by
			quotienting out the nilpotents in $\cO_X$. 
		\item For some open sets $U_\alpha \subset X$ and
			some linearly independent odd elements
			$\theta^i_\alpha$ of $\cO_X(U_\alpha)$ we have
				$\cO_X(U_\alpha) = \cO_X^{\mathrm{red}}
				\otimes \Lambda[\theta^1_\alpha, \dots ,
				\theta^N_\alpha]$.
	\end{enumerate}
	The $U_\alpha$ above are called \emph{coordinate neighborhoods}
	of $(X, \cO_X)$ and $Z_\alpha=(z_\alpha, \theta^1_\alpha, \dots ,
	\theta^N_\alpha)$ are called local coordinates for $(X, \cO_X)$
	if $z_\alpha$ (mod nilpotents) are local coordinates for $(X,
	\cO_X^{\mathrm{red}})$. On overlaps $U_\alpha \cap U_\beta$ we
	have:
	\begin{equation}
		z_\beta = F_{\beta \alpha} (z_\alpha, \theta^j_\alpha), \quad \theta^i_\beta
		= \Psi^i_{\beta \alpha} (z_\alpha,
		\theta^j_\alpha),
		\label{eqn:coordinate_change}
	\end{equation}
	where $F_{\beta \alpha}$ are even and $\Psi_{\beta
	\alpha}$ are odd. We will write such a change of coordinates as
	$Z_\beta = \rho_{\beta,\alpha}(Z_\alpha)$ with $\rho = (F,\Psi^i)$ where
	no confusion should arise.
\end{defn}
\begin{nolabel}
	A $\Lambda$-point of a supercurve $(X, \cO_X)$ is a morphism
	$\varphi:\Spec \Lambda \rightarrow (X, \cO_X)$ over $\Lambda$, namely
	the composition of $\varphi$ with the structure morphism $(X, \cO_X)
	\rightarrow \Spec \Lambda$ is the identity. Locally, a $\Lambda$
	point is given by specifying the images of the local coordinates
	under the even $\Lambda$-homomorphism $\varphi^\sharp :
	\cO_X(U_\alpha) \rightarrow \Lambda$. These local parameters
	$(p_\alpha = \varphi^\sharp(z_\alpha), \pi^i_\alpha =
	\varphi^\sharp(\theta^i_\alpha))$ transform as the coordinates do in
	(\ref{eqn:coordinate_change}).
\end{nolabel}
\begin{nolabel} \label{no:disks}
	The $N=q$ formal superdisk is an ind-superscheme as in the non-super
	situation, namely, let $R=\mathbb{C}[t,\theta^1,\dots, \theta^q]$ and
	let $\fm$ be the maximal ideal generated by $(t, \theta^1,\dots,
	\theta^q)$. We define the superschemes $D^{(n)} = \Spec
	R/\fm^{n+1}$ and we have embeddings $D^{(n+1)}
	\hookrightarrow D^{(n)}$. The formal disk is then 
	\begin{equation*}
		D = \varinjlim_{n \rightarrow \infty} D^{(n)}.
	\end{equation*}
	If we want to emphasize the dimensions of these disks we will
	denote them by $D^{1|q}.$ 
\end{nolabel}
\begin{nolabel}
	Vector bundles of rank $(p|q)$ over a supermanifold $(X, \cO_X)$ are locally
	free sheaves $\cE$ of $\cO_X$-modules over $X$, of rank $p|q$.
	That is, locally, $\cE$ is isomorphic to $\cO_X^p \oplus (\Pi
	\cO_X)^q$ where $\Pi$ is the parity change functor.

	An example is the \emph{tangent bundle} to a $p|q$-dimensional supermanifold $(X,
	\cO_X)$; it is a rank $p|q$ vector bundle. Its fiber at the
	point $x \in X$ is given as in the non-super case as the subset
	of morphisms in $\Hom(D^{(1)}, X)$ mapping the closed point in
	$D^{(1)}$ to $x$.  The \emph{cotangent bundle} $\Omega^1_X$ of
	$(X, \cO_X)$ is the dual of the tangent bundle.

	Another example is the \emph{Berezinian bundle} of a supermanifold $(X, \cO_X)$. We
	will define
	this bundle by giving local trivializations. Recall
	\cite[$\S 1.10$]{deligne2} that given a free module $L$ of finite
	type over a supercommutative algebra $A$, the superdeterminant is a
	homomorphism 
	\begin{equation*}
		\sdet: \mathrm{GL}(L) \rightarrow \mathrm{GL(1|0)} = A_0^\times,
		\label{}
	\end{equation*}
	defined in coordinates as follows: for a parity preserving
	automorphism $T$ of $A^{p|q}$ with matrix ${K \, L \choose M \,
	N}$ we put:
		\begin{equation*}
			\sdet (T) = \det(K - L N^{-1}M) \det(N)^{-1}.
			\label{}
		\end{equation*}
	With this definition we can now define the Berezinian of the
	module $L$ as the following $A$-module denoted $\ber(L)$. Let $\{e_1, \dots, e_{p+q}\}$
	be a basis of
	$L$ where the first $p$ elements are even and the last $q$ are
	odd. This basis defines a one-element basis of $\ber (L)$ denoted
	by $[e_1 \dots e_{p+q}]$ of parity $q \mod 2$. Given an automorphism
	$T$ of $L$ we put:
	\begin{equation*}
		[Te_1 \dots Te_{p+q}] = \sdet (T) [e_1 \dots e_{p+q}].
		\label{}
	\end{equation*}
	This makes $\ber (L)$ a well defined rank $1|0$ $A$-module when $q$ is even and a
	rank $0|1$ $A$-module when $q$ is odd. Now we can define the
	Berezinian bundle of $(X, \cO_X)$ as $\ber_X = \ber(\Omega_X^1)$. 
	
	The definition of coherent and quasi-coherent sheaves is exactly
	the same as in the non-super case, in particular for super
	manifolds it follows that the structure sheaf is coherent
	\cite{vaintrob1}.
	\label{no:berezinian.def}
\end{nolabel}
\begin{nolabel} \label{no:n=2superconformal}
	Given an $N=1$ supercurve $(X, \cO_X)$ and an extension of $\cO_X$ by an invertible
	sheaf $\cE$:  
	\begin{equation}
		0 \rightarrow \cO_X \rightarrow \hat{\cE}\rightarrow \cE
		\rightarrow 0,
		\label{eq:extensionz}
	\end{equation}
	we can construct an $N=2$ supercurve $(Y, \cO_Y)$ canonically. Its local coordinates
	are given by $(z_\alpha, \theta_\alpha, \rho_\alpha)$, where $(z_\alpha,
	\theta_\alpha)$ are local coordinates of $X$ and $\rho_\alpha$ are local sections of
	$\cE$. In each coordinate patch $U_\alpha$ we can construct the form $d z_\alpha - d\theta_\alpha
	\rho_\alpha$. We say that the $N=2$ supercurve $(Y, \cO_Y)$ is \emph{superconformal} if
	this form is globally defined up to multiplication by a function.
	
	This
	happens if on overlaps $U_\alpha \cap U_\beta$ we
	have (see (\ref{eqn:coordinate_change})):
	\begin{equation}
		\rho_\beta = \sdet \begin{pmatrix} \partial_z F &
			\partial_z \Psi \\ \partial_\theta F &
			\partial_\theta \Psi \end{pmatrix} \rho_\alpha +
			\frac{\partial_\theta F}{\partial_\theta \Psi}.
		\label{eqn:cocycle_extension}
	\end{equation}
	Here $\sdet$ is the \emph{superdeterminant} of an automorphism 
	defined above, which can be written as 
	\begin{equation*}
		\sdet \begin{pmatrix} \partial_z F &
		\partial_z \Psi \\ \partial_\theta F &
		\partial_\theta \Psi \end{pmatrix} =
		D \left( \frac{D F}{D \Psi} \right),
		\label{}
	\end{equation*}
	where $D = \partial_\theta + \theta \partial_z$. 
	Conversely, if (\ref{eqn:cocycle_extension}) is satisfied on
	overlaps, the cocycle condition is satisfied and we have an
	extension as in (\ref{eq:extensionz}). 
	
	Therefore to each $N=1$ supercurve $(X, \cO_X)$, we canonically associate a
	$N=2$ superconformal curve $(Y, \cO_Y)$. 
	
	From
	(\ref{eqn:cocycle_extension}) we see that we have an exact
	sequence of sheaves on $Y$:
	\begin{equation*}
		0 \rightarrow \cO_X \rightarrow \cO_Y \rightarrow \ber_X
		\rightarrow 0,
		\label{}
	\end{equation*}
	where $\ber_X$ is the Berezinian bundle on $(X, \cO_X)$. The last
	map $\hat{D}:\cO_Y \rightarrow \ber_X$ is given in the above local coordinates, by the
	differential 
	operator $\partial_{\rho_\alpha}$. 

	Introducing new coordinates
	\begin{equation*}
			\hat{z}_\alpha = z_\alpha - \theta_\alpha
			\rho_\alpha, \qquad
			\hat{\theta}_\alpha = \theta_\alpha, \qquad
			\hat{\rho}_\alpha = \rho_\alpha,
		\label{}
	\end{equation*}
	we obtain on overlaps $U_\alpha \cap U_\beta$:
	\begin{equation}
			\hat{z}_\beta = F(\hat{z}_\alpha,
			\hat{\rho}_\alpha)
			+ \frac{D F(\hat{z}_\alpha,\hat{\rho}_\alpha)}{D
			\Psi(\hat{z}_\alpha,
			\hat{\rho}_\alpha)}\Psi(\hat{z}_\alpha,
			\hat{\rho}_\alpha), \quad
			\hat{\rho}_\beta = \frac{D F(\hat{z}_\alpha,
			\hat{\rho}_\alpha)}{D \Psi(\hat{z}_\alpha,
			\hat{\rho}_\alpha)},
		\label{eqn:dual_transition}
	\end{equation}
	where $D = \partial_\theta + \theta \partial_z$ in local coordinates $(z, \theta,
	\rho)$ as above.  
	
	We see from (\ref{eqn:dual_transition}) that $\cO_Y$ contains the
	structure sheaf of 
	another $N=1$ supercurve $(\hat{X}, \cO_{\hat{X}})$, whose local
	coordinates are $(\hat{z}_\alpha, \hat{\rho}_\alpha)$. We call
	$(\hat{X}, \cO_{\hat{X}})$ the \emph{dual curve} of $(X, \cO_X)$. 

	Finally, we define an $N=1$ superconformal curve as an $N=1$ supercurve $(X,\cO_X)$
	which is self-dual. We see from
	(\ref{eqn:dual_transition})
	that the transition functions $F, \Psi$ must satisfy:
	\begin{equation}
		DF = \Psi D\Psi ,
		\label{eqn:super_conformal1}
	\end{equation}
	for $(X, \cO_X)$ to be superconformal. In this case the operator
	$D_\alpha = \partial_{\theta_\alpha} + \theta_\alpha
	\partial_{z_\alpha}$ transforms as:
	\begin{equation}
		D_\beta = (D \Psi)^{-1} D_\alpha,
		\label{eq:D_transforms_as}
	\end{equation}
	hence in this situation the supercurve $(X, \cO_X)$ carries a
	$0|1$-dimensional distribution $D$ such that $D^2$ is nowhere
	vanishing (since $D^2 = \partial_z$ in local coordinates). 
\end{nolabel}
\begin{rem} \label{rem:3.14} An equivalent definition of $N=1$ and $N=2$ superconformal curves
	was given by Manin \cite{manin3} (under the name SUSY curves). Let $X$ be a complex
	supermanifold of dimension $1|N$ ($N=1$ or $2$). When $N=1$ we
	say that a locally free direct subsheaf $\cT^1 \subset \cT_X$
	($\cT_X$ is the tangent sheaf of $X$) of rank $0|1$
	for which
	the Frobenius form
	\begin{equation*}
		(\cT^1)^{\otimes 2} \rightarrow \cT^0 := \cT_X /\cT^1, \quad
		t_1 \otimes t_2 \mapsto [t_1, t_2] \mod \cT^1,
		\label{}
	\end{equation*}
	is an isomorphism, is a SUSY structure on $X$. 

	When $N=2$, a SUSY structure consists of two locally free direct
	subsheaves $\cT', \, \cT''$ of $\cT_X$ of rank $0|1$ whose sum in
	$\cT_X$ is direct, they are integrable distributions and the
	Frobenius form 
	\begin{equation*}
		\cT' \otimes \cT'' \rightarrow \cT_X/ (\cT' \oplus \cT''), \quad
		t_1 \otimes t_2 \mapsto [t_1, t_2] \mod (\cT' \oplus
		\cT''),
		\label{}
	\end{equation*}
	is an isomorphism.

	Let  $(X, \cO_X)$ be an $N=1$ supercurve and $D_\alpha$ be a family of vector
	fields in
	$U_\alpha$, such that $D_\alpha$ and $D_\alpha^2$ form a basis for $\cT_X$ on
	$U_\alpha$ and $D_\alpha = G_{\alpha \beta} D_\beta$ on $U_\alpha \cap U_\beta$,
	where $G_{\alpha \beta}$ is a family of invertible even
	functions. The sheaf defined by $\cT^1|_{U_\alpha} = \cO_X D_\alpha$ is a SUSY
	structure in $(X, \cO_X)$ \cite{manin3}. In local coordinates as above,
	the vector fields
	$D_\alpha = \partial_{\theta_\alpha} + \theta_\alpha
	\partial_{z_\alpha}$ satisfy these conditions when $X$ is an $N=1$ superconformal
	curve (see (\ref{eq:D_transforms_as})).  

	The $N=2$ case is similar. Let $(X, \cO_X)$ be an $N=2$ supercurve and  $\left\{ D_\alpha^1,\,D_\alpha^2
	\right\}$ be a family of vector fields such that $D^i_\alpha,
	[D_\alpha^1, D_\alpha^2]$
	generate $\cT_X$ in $U_\alpha$ and, moreover, we have:
	\begin{xalignat*}{2}
		(D_\alpha^1)^2 &= f^1_\alpha D_\alpha^1, &
		(D_\alpha^2)^2 &= f_\alpha^2 D_\alpha^2,\\
		D_\alpha^1 &= F_{\alpha, \beta}^1 D_\beta^1, & D^2_\alpha
		&= F_{\alpha \beta}^2 D_\beta^2 \quad \text{on} \;
		U_\alpha \cap U_\beta, 
	\end{xalignat*}
	where $f^i_\alpha$ and $F_{\alpha \beta}^i$ are even
	functions. Putting $\cT'|_{U_\alpha} = \cO_X D_\alpha^1$ and $\cT''|_{U_\alpha} = \cO_X
	D_\alpha^2$ we obtain an $N=2$ superconformal structure on $(X,\cO_X)$. If the two distributions
	$\cT'$ and $\cT''$ can be distinguished globally, the $N=2$ superconformal curve is
	called orientable
	and a choice of one of these distributions is called its orientation.
	
	It is clear that the construction given in
	\ref{no:n=2superconformal} gives an oriented $N=2$ superconformal curve;
	conversely, given such a curve, we can consider the functor $X
	\rightarrow X/\cT'$ (recall that $\cT'$ is integrable therefore
	this quotient makes sense). The duality that was explained in
	\ref{no:n=2superconformal} corresponds to the duality $X/\cT'
	\leftrightarrow X/\cT''$. 
\end{rem}
\begin{nolabel}
	\label{no:3.15} Recall  \cite{rabin1} that a $\Lambda$-point of an $N=1$ supercurve
	$X$ transforms as an irreducible divisor of the dual curve $\hat{X}$.
	Indeed, an irreducible divisor of $X$ is given in local coordinates $(z_\alpha,
	\theta_\alpha)$ by
	expressions of the form $P_\alpha= z_\alpha - \hat{z}_\alpha - \theta_\alpha
	\rho_\alpha$. Two divisors $P_\alpha$ and $P_\beta$ are said to correspond
	to each other in the intersection $U_\alpha \cap U_\beta$ if in this
	intersection we have
	\begin{equation*}
		P_\beta(z_\beta, \theta_\beta) = P_\alpha(z_\alpha, \theta_\alpha)
		g(z_\alpha,\theta_\alpha)
	\end{equation*}
	for some even invertible function $g(z_\alpha, \theta_\alpha)$ (we
	consider Cartier divisors). It is easy to see that  
	the parameters $\hat{z}_\alpha, \rho_\alpha$ transform as in
	(\ref{eqn:dual_transition}), namely as the parameters of a $\Lambda$-point
	of $\hat{X}$. 
\end{nolabel}
\begin{nolabel}
	\label{no:3.16} We can define a theory of contour integration on an $N=1$
	superconformal curve 
	as in \cite{friedan}, \cite{mcarthur}, \cite{rogers}. We describe briefly a
	generalization to arbitrary $N=1$ supercurves due to Bergvelt and Rabin (cf.
	\cite{rabin1}). For simplicity, we will work  
	in the analytic category. Let us define a \emph{super contour} to be a
	triple $\Gamma = (\gamma, P, Q)$ consisting of an ordinary contour $\gamma$
	on the reduction $|X|$ and two Cartier divisors as in \ref{no:3.15} such
	that their reductions to $|X|$ are the endpoints of $\gamma$. If in local coordinates
	\begin{equation*}
		P=z - \hat{p} - \theta \hat{\pi}, \qquad Q=z - \hat{q} - \theta
		\hat{\xi},
	\end{equation*}
	then the
	corresponding $\Lambda$-points of the dual curve $\hat{X}$ are given 
	by $(\hat{p}, \hat{\pi})$ and $(\hat{q}, \hat{\xi})$. Let
	$z=\hat{p}_{\mathrm{rd}}$ and $z = \hat{q}_{\mathrm{rd}}$ be the equations
	for the reductions of these points, i.e. the endpoints for $\gamma$. 
	We define the integral of a section $\omega_\alpha = D \hat{f}_\alpha$ of
	the Berezinian sheaf of $X$ (Recall that $D : \cO_{\hat{X}}
	\twoheadrightarrow \ber{X}$) along $\Gamma$ by:
	\begin{equation*}
		\int_P^Q \omega = \int_P^Q D \hat{f} = \hat{f}(\hat{q}, \hat{\xi})
		- \hat{f}(\hat{p}, \hat{\pi}),
	\end{equation*}
	where we assume that the contour connecting $P$ and $Q$ lies in a single
	simply connected open set $U_\alpha$. If the contour traverses several open
	sets then we need to choose intermediate divisors on each overlap and we
	have to prove that the resulting integral is independent of these divisors.
	In what follows we will only need the integration in a sufficiently
	``small'' open set $U_\alpha$ (the formal disk around a point). 

	Dually, we can integrate sections of $\ber_{\hat{X}}$ along contours in
	$X$. Indeed, let $\gamma$ be a path in the topological space $|X|$ and two
	$\Lambda$-points $P,\, Q$ of $X$ whose reduced parts are the end-points of
	$\gamma$. Let $\hat{\omega} \in \ber_{\hat{X}}(U_\alpha)$ and suppose that
	$\gamma$ lies in a simply connected open $U_\alpha$. Then $\hat{\omega} =
	\hat{D} f$ for some function $f \in \cO_X(U_\alpha)$, and we put
	\begin{equation*}
		\int_P^Q \hat{\omega} = f(Q) - f(P).
	\end{equation*}

	As it is shown in \cite{rabin1}, this theory of integration can be
	understood in terms of a theory of contour integration on the corresponding
	$N=2$ superconformal curve (cf. \cite{cohn}). For this let $X$ and
	$\hat{X}$ be an $N=1$ supercurve and its dual, and let $Y$ be the
	corresponding $N=2$ superconformal curve. We have two short exact
	sequences of sheaves in $Y$:
	\begin{equation*}
		\begin{aligned}
			& 0 \rightarrow  \cO_X \rightarrow \cO_Y
			\xrightarrow{D^-}  \ber_{\hat{X}} \rightarrow
			0, \\
			& 0 \rightarrow  \cO_{\hat{X}} \rightarrow \cO_Y
			\xrightarrow{D^+}  \ber_X \rightarrow 0.
		\end{aligned}
	\end{equation*}
	We can define a sheaf operator on $\cO_Y^{\oplus 2}$ by the component-wise
	action of the differential operators $(D^-, D^+)$. It is shown in
	\cite{rabin1} that for $U$ a simply connected open in $|Y|=|X|$ and $(f,g)$
	a section of $\cO_Y^{\oplus 2}(U)$ such that $(D^-, D^+)(f,g)= 0$, 
	there exists a section $H \in \cO_Y(U)$, unique up to an additive constant,
	such that $(f,g) = (D^- H, D^+ H)$. Let $\cM$ be the subsheaf of
	$\cO_Y^{\oplus 2}$ consisting of closed sections $(f,g)$ as above. It
	follows that $\cM =  \ber_{X} \oplus  \ber_{\hat{X}}$. A
	super contour in $Y$ consists of a triple $(\gamma, P, Q)$ where $P$ and
	$Q$ are $\Lambda$-points of $Y$ such that their reduced points are the
	endpoints of $\gamma$. If $\gamma$ is supported on a simply connected open set
	$U$, then any section $\omega \in \cM(U)$ can be written as $(D^- H, D^+
	H)$ and we put:
	\begin{equation*}
		\int_P^Q \omega = H(Q) - H(P).
	\end{equation*}
	The extension to contours not lying in a single simply connected $U$ is
	straightforward but we will not need it. 	
\end{nolabel}
\begin{nolabel}
	We will define in general a superconformal $N=n$ supercurve to
	be a supercurve such that in some coordinate system $Z_\alpha =
	(z_\alpha, \theta^i_\alpha)$, the differential form 
	\begin{equation}
		\omega = d z_\alpha + \sum_i \theta^i_\alpha d
		\theta^i_\alpha
		\label{eq:holaquetal}
	\end{equation}
	is well defined up to multiplication by a function. It is easy to show that
	this definition agrees with the definition above in the $N=1$ and
	$N=2$ cases (cf. $\S$\ref{no:supern=1group} and
	$\S$\ref{no:n=2supergroup}).

	A set of coordinates $Z= (z,\theta^i)$ such that the form $\omega$ has
	the form (\ref{eq:holaquetal}) (up to multiplication by a function) will
	be called $SUSY$ coordinates (or \emph{coordinates compatible with the superconformal
	structure}).

	Let $(z,\theta)$ and $(z',\theta')$ be two local coordinates
	compatible with a (local) superconformal structure on an $N=1$
	supercurve $(X,\cO_X)$. Denote $D = \partial_\theta + \theta \partial_z$
	and $D'= \partial_{\theta'} + \theta' \partial_{z'}$. Let $G$ be the invertible function such that 
	$D = G D'$ (cf. (\ref{eq:D_transforms_as})).
	We define the \emph{Schwarzian derivative} of $(z',\theta')$ with
	respect to $(z,\theta)$ to be the (odd) function
	\begin{equation}
		\sigma(G) = \frac{D^3 G}{G} - 2 \frac{D G D^2
		G}{G^2}.
		\label{eq:schwarzian_1_def}
	\end{equation}
\end{nolabel}
\begin{defn}
	A \emph{superprojective} structure on an $N=1$ superconformal
	curve over 
	$\Lambda$ is a (maximal) atlas consisting of coordinates
	$(z_\alpha, \theta_\alpha)$ compatible with the superconformal
	structure,  and such that its
	transition functions are fractional linear transformations,
	 that is, changes of coordinates of the
	form:
	\begin{equation*}
			z' = \frac{a z + b + \alpha \theta}{c z + d + \beta
			\theta}, \qquad
			\theta' = \frac{\gamma z + \delta + e \theta}{c z + d +
			\beta \theta},
	\end{equation*}
	for some even constants $a$, $b$, $c$, $d$ and $e \in \Lambda$ and some odd constants
	$\alpha$, $\beta$, $\gamma$ and $\delta \in \Lambda$, such that
	\begin{equation*}
		\sdet
		\begin{pmatrix}
			a & b & \alpha \\ c & d & \beta \\ \gamma & \delta & e
		\end{pmatrix} = 1
	\end{equation*}
\end{defn}
\begin{prop}[\cite{manin2} Proposition 4.7] \label{prop:3.22}
	Let $(z,\theta)$ and $(z', \theta')$ be two local coordinates on
	$(X, \cO_X)$. The following statements are equivalent:
	\begin{enumerate}
		\item $(z,\theta)$ and $(z',\theta')$ are compatible with
			a common superconformal structure and $\sigma =
			0$.
		\item $(z,\theta)$ and $(z', \theta')$ define the same
			superprojective structure.
	\end{enumerate}
\end{prop}
\section{The associated vector bundles} \label{sec:vector_bundles_section}
\subsection{The groups $\mathrm{Aut} \cO$} \label{sub:group_aut_o}
\begin{nolabel}
	We start this section by describing the groups of changes of coordinates
	in the formal superdisk $D^{1|N}$. We analyze in detail their
	corresponding Lie superalgebras in the cases $N=1$ and $N=2$. We then
	define principal bundles for these groups over any smooth supercurve. 

	In this section, we let $\Lambda$ be a Grassman algebra over $\mathbb{C}$. We will work in
	the category of superschemes over $\Lambda$ unless explicitly stated. When
	we work with a
	supergroup $G$, we will be interested in its $\Lambda$-points. 
	\label{no:4.1.explanation}
\end{nolabel}
  Let $\mathrm{SSch}/k$
	be the category of superschemes over a field $k$ and let $\mathrm{Set}$
	be the category of sets.
	Fix a non negative integer $N$ and a separated superscheme $X$ of finite
	type over $k$ (cf.
	\ref{defn:super_scheme}). Let
	$D^{(m)}$ be as in
	\ref{no:disks} and $D^{1|N}$ be the formal superdisk.  Define a family of contravariant
	functors $F_m : \mathrm{SSch}/k \rightarrow \mathrm{Set}$
	\begin{equation*}
		F_m(Y) = \Hom_k(Y \times_k D^{(m)}, X).
		\label{}
	\end{equation*}
	The proof of the following proposition is standard:
\begin{prop} \label{no:jetschemes}
	The functors $F_m$ are representable by superschemes $X_m$.
\end{prop}
	Note in particular that $X_0 = X$, and when $N=1$ we see that $X_1$ is
	the total tangent space of $X$. 
	
	The embeddings $D^{(m)} \hookrightarrow D^{(m+1)}$ induce
	projections $X_{m+1} \rightarrow X_{m}$ and we define the
	\emph{Jet superscheme} of $X$ as
	\begin{equation*}
		JX = \varprojlim_{m \rightarrow \infty} X_m
		\label{}
	\end{equation*}

\begin{nolabel} \label{no:auto1} Let us analyze first the case $N=1$. 
	Consider the group of continuous (even) automorphisms of the
	topological commutative superalgebra $\Lambda[ [Z]]$, where 
	$Z=(z,\theta)$ are topological generators. Such an automorphism
	is
	given by a pair of power series
	\begin{equation*}
		\begin{aligned}
			z &\mapsto a_{1,0} z + a_{0,1} \theta + a_{1,1} z
			\theta + \dots \\
			\theta &\mapsto b_{1,0} z + b_{0,1} \theta +
			b_{1,1} z \theta + \dots,
		\end{aligned}
		\label{}
	\end{equation*}
	where the matrix ${a_{1,0} \, a_{0,1} \choose b_{1,0} \,
	b_{0,1}}$ is in $GL(1|1)$\footnote{Here and further, $GL(p|q)$ is the group
	of even
	automorphisms of a $p|q$ dimensional module over $\Lambda$.}. Denote this
	supergroup by
	$\mathrm{Aut} \cO^{1|1}$. In what follows we will analyze its
	$\mathbb{C}$-points. 

        This supergroup is a semidirect product of $GL(1|1)$ and a pro-unipotent
	super group, namely, the subgroup $\mathrm{Aut}_+\cO^{1|1}$ of automorphisms
	where ${ a_{1,0} \, 
	a_{0,1} \choose b_{1,0} \, b_{0,1}} = \Id$. In fact, 
	\begin{equation*}
		\mathrm{Aut}_+ \cO^{1|1} = \varprojlim_{n\rightarrow \infty} \Spec
		\mathbb{C}[a_{1,1}, b_{1,1}, a_{2,0}, b_{2,0}, \dots,
		a_{n,1}, b_{n,1}].
		\label{}
	\end{equation*}
	Let $\fm$ be the maximal ideal of $\mathbb{C}[Z]$
	generated by $(z, \theta)$. We have 
	\begin{equation*}
		\mathrm{Aut}_+ \cO^{1|1} = \varprojlim_{n \rightarrow \infty}
		\mathrm{Aut} (\mathbb{C}[Z]/\fm^n)
		\label{}.
	\end{equation*}
	Similarly for its Lie superalgebra $\mathrm{Der}_+ \cO^{1|1}$, we have
	\begin{equation*}
		\mathrm{Der}_+ \cO^{1|1} = \varprojlim_{n \rightarrow \infty}
		\mathrm{Der} (\mathbb{C}[Z]/\fm^n),
		\label{}
	\end{equation*}
	where for each $\mathbb{C}$-superalgebra $R$, we denote $\mathrm{Der} (R)$ the
	Lie superalgebra of derivations of $R$.
	The exponential map is an isomorphism at each step, giving an
	isomorphism $\mathrm{exp}: \mathrm{Der}_+ \cO^{1|1} \rightarrow
	\mathrm{Aut}_+ \cO^{1|1}$.

	The linearly compact Lie superalgebra $\mathrm{Der}_0 \cO^{1|1} = \mathrm{Lie}
	(\mathrm{Aut} \cO^{1|1})$ has the following topological basis:
	\begin{xalignat*}{2}
		z^n \partial_z  &\quad (n \geq 1), & z^n \partial_\theta &
		\quad (n \geq 1), \\
		z^n \theta \partial_z & \quad (n \geq 0), & z^n \theta
		\partial_\theta & (\quad n \geq 0),
	\end{xalignat*}
	or the following one $(n \geq 0)$:
	\begin{equation}
		\begin{aligned}
		T_n &=  -
		z^{n+1} \partial_z - (n+1) z^n \theta \partial_\theta, &
		J_n &= -z^n \theta \partial_\theta, \\
		Q_n &= -z^{n + 1} \partial_\theta, & H_n &=
		z^{n} \theta \partial_z.
	\end{aligned}
		\label{eq:4.2.7}
	\end{equation}
	These elements satisfy the
	commutation relations of the $N=2$ algebra
	(\ref{eq:n=2commutations_new_fields}) 
	for $n \geq 0$. In particular, we see that $\mathrm{Der}_0 \cO^{1|1}$ is
	the formal 
	completion of the Lie algebra $W(1|1)_<$.
	The Lie subalgebra $\mathrm{Der}_+ \cO$ is
	topologicaly generated by the same vectors with $n \geq 1$.
\end{nolabel}
\begin{nolabel} \label{no:supern=1group}
	We now turn our attention to the superconformal
	$N=1$ case. Consider the differential 
	form
	$\omega = d z + \theta d \theta$ on the formal superdisk $D^{1|1}$, and
	the supergroup $\mathrm{Aut}^\omega \cO^{1|1}$ of
	automorphisms of $D^{1|1}$ preserving this form, up to
	multiplication by a function. This is a
	subgroup of $\mathrm{Aut} \cO^{1|1}$ whose Lie superalgebra
	$\mathrm{Der}^\omega_0 \cO^{1|1}$ consist of
	derivations $X$ in $\mathrm{Der}_0 \cO^{1|1}$ such that $L_X \omega = f
	\omega$ for some formal power series $f$ (here $L_X$ denotes the Lie
	derivative). More explicitly, the linearly compact Lie superalgebra
	$\mathrm{Der}^\omega_0 \cO^{1|1}$
	 is topologically generated by 
	\begin{equation}
		\begin{aligned}
			L_n &= - \frac{n+1}{2} z^n \theta \partial_\theta
			- z^{n+1} \partial_z, \qquad n \in \mathbb{Z}_+ \\
			G_n &= - z^{n+1/2} (\partial_\theta - \theta
			\partial_z), \qquad n \in \frac{1}{2} + \mathbb{Z}_+.
		\end{aligned}
		\label{eq:neveu_shwarz_generators}
	\end{equation}
	These generators satisfy the commutation relation of
	the Neveu-Schwarz algebra as defined in (\ref{eq:n=1fouriercommutations}).
	In particular, we see that $\mathrm{Der}_0^\omega \cO^{1|1}$ is the formal
	completion of the Lie superalgebra $K(1|1)_<$.
	
	An automorphism of the formal superdisk is
	determined by two power series $F(Z),
	\Psi(Z)$ which are the images of the generators $Z=(z,
	\theta)$. Under this transformation we have (recall
	$\partial_\theta$ is an odd derivation)
	\begin{equation*}
		\begin{aligned}
			dz + \theta d\theta &\mapsto \partial_z F dz -
			\partial_\theta F d\theta + \Psi (\partial_z\Psi
			dz + \partial_\theta \Psi d\theta) \\
		       	 &= (\partial_z F + \Psi \partial_z \Psi )
			dz - (\partial_\theta F - \Psi \partial_\theta
			\Psi) d\theta,
		\end{aligned}
		\label{}
	\end{equation*}
	therefore we get that, in order for $\omega$ to be preserved up to
	multiplication by a function, we need
		$(\partial_\theta F - \Psi \partial_\theta \Psi) =
			-\theta (\partial_z F +\Psi \partial_z\Psi)$,
	and this is equivalent to (\ref{eqn:super_conformal1}). 
\end{nolabel}
\begin{nolabel} \label{no:n=2supergroup}
	Finally we turn our attention to the (oriented) superconformal $N=2$ case.
	For this we consider the differential form $\omega= dz + \theta^1 d\theta^1 + \theta^2
	d \theta^2$ on the formal superdisk $D^{1|2}$. We want to analyze the group
	of automorphisms of $D^{1|2}$ preserving this form in the sense
	of the previous paragraph \ref{no:supern=1group}. Such an
	automorphism is determined by an even power series $F(Z)$
	and two odd power series $\Psi^1(Z)$ and $\Psi^2(Z)$, where $Z =
	(z,\theta^1, \theta^2)$ are the coordinates on $D^{1|2}$. Under such a
	change of coordinates, the 
	differential form $\omega$ changes to:
	\begin{multline} \label{mult:formtransform}
		\partial_z F dz - \partial_{\theta^1} F d\theta^1 -
		\partial_{\theta^2} F d \theta^2 + \Psi^1 \left(
		\partial_z \Psi^1 dz + \partial_{\theta^1} \Psi^1 d\theta^1 +
		\partial_{\theta^2} \Psi^1 d\theta^2
		\right) + \\ + \Psi^2 \left( \partial_z \Psi^2 dz +
		\partial_{\theta^1} \Psi^2 d \theta^1 +
		\partial_{\theta^2} \Psi^2 d\theta^2 \right) = \\ =
		\left( \partial_z F + \Psi^1 \partial_z \Psi^1 + \Psi^2
		\partial_z \Psi^2 \right) d z + \left( -
		\partial_{\theta^1} F + \Psi^1 \partial_{\theta^1}\Psi^1
		+ \Psi^2 \partial_{\theta^1} \Psi^2\right) d\theta^1 + \\
		+ 
		\left( - \partial_{\theta^2} F + \Psi^1
		\partial_{\theta^2} \Psi^1 + \Psi^2 \partial_{\theta^2}
		\Psi^2 \right) d\theta^2.	
	\end{multline}
	Collecting terms, imposing that the form $\omega$ is
	preserved up to multiplication by a function, and defining the
	differential operators $D^i = \partial_{\theta^i} + \theta^i
	\partial_z$ we obtain that the automorphisms we are
	considering satisfy the equations:
	\begin{equation} \label{eq:allowedn=2}
		D^i F = \Psi^1 D^i \Psi^1 + \Psi^2 D^i \Psi^2, \qquad i =
		1,2.
	\end{equation}

	Note also that a particular case of (\ref{mult:formtransform})
	when $F= z - \tfrac{1}{2} \theta^1 \theta^2$, $\Psi^1 =
	\tfrac{i}{2} (\theta^2 - \theta^1)$ and $\Psi^2 = \tfrac{1}{2}
	(\theta^1 + \theta^2)$ transforms the form
	\begin{equation*}
	\omega \mapsto dz + \theta^2 d \theta^1 = dz -
	d\theta^1 \theta^2,
	\end{equation*}
	and the supergroup of automorphisms of $D^{1|2}$ preserving the
	latter form is the supergroup of changes of coordinates preserving an
	$N=2$
	superconformal structure as in \ref{no:n=2superconformal}. 

%
	The linearly compact Lie superalgebra $\mathrm{Der}^\omega_0 \cO^{1|2} =
	\mathrm{Lie}
	(\mathrm{Aut}^\omega \cO^{1|2})$ is topologicaly generated by:
	\begin{equation}
		\begin{aligned}
			L_n &= - z^{n+1} \partial_z -
			\frac{n+1}{2} z^n \left( \theta^1
			\partial_{\theta^1} + \theta^2
			\partial_{\theta^2} \right), \quad  n \in \mathbb{Z}_+\\
			G^{(2)}_{n} &= + z^{n+1/2} \left( \theta^2
			\partial_z - \partial_{\theta^2}
			\right) - \left( n + \frac{1}{2} \right)
			z^{n-1/2} \theta^1 \theta^2 \partial_{\theta^1}, \quad n
			\in \frac{1}{2} + \mathbb{Z}_+
			\\
			G^{(1)}_n &= + z^{n+1/2} \left( \theta^1
			\partial_z - \partial_{\theta^1} \right) + \left(
			n + \frac{1}{2} \right) z^{n-1/2} \theta^1
			\theta^2 \partial_{\theta^2}, \quad n \in
			\frac{1}{2} + \mathbb{Z}_+\\
			J_n &= -i z^n \left(\theta^2 \partial_{\theta^1} -
			\theta^1 \partial_{\theta^2} \right) \quad n \in
			\mathbb{Z}_+.
		\end{aligned}
		\label{eq:4.4.8}
	\end{equation}
	These operators satisfy
	the commutation relations of the $N=2$ generators as in
	(\ref{eq:n=2superconf})  for $n \geq
	0$. We see that the Lie superalgebra $\mathrm{Der}_0^\omega \cO^{1|2}$ is
	the formal completion of the Lie superalgebra $K(1|2)_<$.

	It is useful to consider \emph{complex} coordinates $\theta^\pm = \theta^1
	\pm i \theta^2$, and derivations $D^\pm = \tfrac{1}{2} (D^1 \pm i D^2)$. In
	the coordinates $(z,\theta^+, \theta^-)$, these derivations are expressed
	as:
	\begin{equation*}
		D^\pm = \partial_{\theta^\mp} + \frac{1}{2} \theta^\pm
		\partial_z.
	\end{equation*}
	If we change coordinates by $\rho = (F, \Psi^+, \Psi^-)$, with
	$\Psi^\pm = \Psi^1 \pm i \Psi^2$, the superconformal
	condition (\ref{eq:allowedn=2}) reads 
	\begin{equation}
		D^\pm F = \frac{1}{2} \Psi^+ D^\pm \Psi^- + \frac{1}{2} \Psi^-
		D^\pm \Psi^+.
		\label{eq:7.7.6.a}
	\end{equation}
	Therefore, under a change of coordinates $(z_\alpha,
	\theta^\pm_\alpha) \mapsto (z_\beta, \theta^\pm_\beta)$, the operators
	$D^\pm$ transform as
	\begin{equation}
		D^\pm_\alpha = (D^\pm_\alpha \Psi^-_{\beta,\alpha}) D^+_\beta +
		(D^\pm_\alpha \Psi^+_{\beta,\alpha}) D_\beta^-.
		\label{eq:7.7.7.a}
	\end{equation}
	In the following sections,  we will consider only \emph{oriented}
	superconformal $N=2$
	supercurves (cf. remark \ref{rem:3.14}), namely those for which there
	exists a coordinate atlas $(U_\alpha,z_\alpha, \theta^\pm_\alpha)$ such
	that on overlaps we have
	\cite{cohn}:
	\begin{equation}
		D^\pm_\alpha \Psi^\pm_{\beta,\alpha} = 0.
		\label{eq:7.7.8.a}
	\end{equation}
	In these coordinates, the topological generators of the Lie superalgebra
	$\mathrm{Der}_0^\omega \cO^{1|2}$ are
	expressed as:
	\begin{equation}
		\begin{aligned}
			L_n &= - z^{n+1} \partial_z - \frac{n+1}{2}z^n (\theta^+
			\partial_{\theta^+} + \theta^- \partial_{\theta^-}), \quad n
			\in \mathbb{Z}_+ \\
			J_n &= - z^n (\theta^+ \partial_{\theta^+} - \theta^-
			\partial_{\theta^-}), \quad n \in \mathbb{Z}_+\\
			G^\pm_n &=- z^{n+1/2} \left( \partial_{\theta^\pm} -
			\frac{1}{2} \theta^\mp \partial_z \right) -
			\frac{n +1/2}{2} z^{n-1/2} \theta^\pm \theta^\mp
			\partial_{\theta^\pm}, \quad  n \in \frac{1}{2} + \mathbb{Z}_+
		\end{aligned}
		\label{eq:7.7.9.a}
	\end{equation}
	where as before we have $G^\pm = \tfrac{1}{2}( G^{(1)} \mp i G^{(2)})$.

	Recall from \ref{no:n=2superconformal} that an oriented superconformal
	$N=2$ supercurve $(Y, \cO_Y)$
	projects onto two
	$N=1$ supercurves $X$ and its dual $\hat{X}$.
	Defining new coordinates $(u, \theta^+, \theta^-)$, where $u
	= z - \tfrac{1}{2} \theta^+ \theta^-$, we 
	see that equations (\ref{eq:7.7.6.a}), for a change of
	coordinates $\rho=(G = F + \tfrac{1}{2}\Psi^+ \Psi^-, \Psi^+, \Psi^-)$
	are expressed in these coordinates as:
	\begin{equation}
			D^- G = \Psi^- D^- \Psi^+, \qquad
			D^+ G = 0.
		\label{eq:7.8.3'.a}
	\end{equation}
	Moreover, the operators $D^\pm$ are expressed as:
	\begin{equation}
		D^+ = \partial_{\theta^-}, \qquad D^- = \partial_{\theta^+} +
		\theta^- \partial_u.
		\label{eq:7.8.4.a}
	\end{equation}
	Note that the coordinate $\theta^-$ does not appear in the
	transition functions for $u, \theta^+$, therefore these coordinates give the
	topological space $|Y|$ the structure of an $N=1$ supercurve. Let us call
	this curve $X$. Similarly, if we define $u' = z - \tfrac{1}{2} \theta^+
	\theta^-$ we obtain that $u', \,\theta^-$ defines the dual curve $(\hat{X},
	\cO_{\hat{X}})$.

	It follows from the above discussion, that given a change of coordinates
	$\rho = (G, \Psi^+) \in \mathrm{Aut} \cO^{1|1}$, we obtain uniquely a
	change of coordinates $\rho = (G, \Psi^+, \Psi^-) \in \mathrm{Aut}^\omega
	\cO^{1|2}$, where $\Psi^- = D^- G / D^- \Psi^+$. This map induces an
	isomorphism of supergroups from $\mathrm{Aut} \cO^{1|1}$ to the identity
	component of $\mathrm{Aut}^\omega \cO^{1|2}$.  
	 This isomorphism corresponds
	to the isomorphism of Lie superalgebras
	$K(1|2) \equiv W(1|1)$ (cf. \cite{kacleur}), and
	has a geometric counterpart (cf. \cite{vaintrob2})
	relating the moduli space of (oriented) superconformal $N=2$ supercurves
	and the moduli space of $N=1$ supercurves.
	
\end{nolabel}
\begin{rem} \label{rem:superconformaljets}
	Let $X$ be a superconformal $N=n$ supercurve.
	Then for some coordinate atlas $Z_\alpha = (z_\alpha,
	\theta^1_\alpha, \dots, \theta^n_\alpha)$, the form $\omega = d_z +
	\sum_{i=1}^n \theta^i d\theta^i$ is globally defined up to multiplication
	by a function.
	 Let $\omega$ be that form on the superdisk
	$D^{1|N}$ and on $D^{(m)}$ as well. Define the functors $F^\omega_m :
	\mathrm{SSch}/k \rightarrow \mathrm{Set}$ by:
	\begin{equation*}
		F_m(Y) = \Hom_k^\omega (Y \times_k D^{(m)}, X),
		\label{}
	\end{equation*}
	where $\Hom^\omega$ denotes the set of morphisms preserving the form
	$\omega$ (up to multiplication by a function). It follows in the same way
	as in Proposition \ref{no:jetschemes} that the functors $F_m$ are
	representable by superschemes $X_m^\omega$. This allows us to define
	the superscheme
	\begin{equation*}
		JX^\omega = \varprojlim_{m \rightarrow \infty} X^\omega_m,
		\label{}
	\end{equation*}
	parametrizing maps $D \rightarrow X$ preserving the
	superconformal structure.
\end{rem}
\begin{nolabel}
	Let $X$ be an $N=n$ supercurve and let $x \in X$. If $Z = (z,
	\theta^1, \dots, \theta^n)$ are local coordinates at $x$ and $\cO_x$ denotes the
	completion of the local ring at $x$, we have an isomorphism
	\begin{equation*}
		\cO_x \equiv \mathbb{C}[ [Z]],
		\label{}
	\end{equation*}
	where we should replace $\mathbb{C}$ by $\Lambda$ if $X$ is
	defined over $\Lambda$. For the purposes of this section it is
	enough to consider curves over $\mathbb{C}$, the relative case
	follows easily.
	Let $\mathrm{Aut}_x$ denote the set of local
	coordinates $Z=(z,\theta^i)$ at $x$. In the algebraic setting we
	mean by coordinates an \'{e}tale map $Z: X \rightarrow
	\mathbb{A}^{1|n}$. The set $\mathrm{Aut}_x$ is a torsor for the
	group $\mathrm{Aut} \cO^{1|n}$. The torsors $\mathrm{Aut}_x$ glue
	to form an $\mathrm{Aut} \cO^{1|n}$-torsor $\mathrm{Aut}_X$.
	Indeed $\mathrm{Aut_X}$ consists of pairs $(x,Z)$ where $x$ is a
	point in $X$ and $Z=(z,\theta^i)$ is a local coordinate at $x$. The
	action of $\mathrm{Aut} \cO^{1|n}$ on the fibers is by changes of
	coordinates. The torsor $\mathrm{Aut}_X$ may be described as an
	open subscheme of $JX$ consisting of jets of maps $D^{1|n}
	\rightarrow X$ such that their 1-jet is in $GL(1|n)$. Since we
	can cover $X$ by Zariski open subschemes $U_\alpha$ and \'{e}tale maps
	$f_\alpha: U_\alpha \rightarrow  \mathbb{A}^{1|n}$ we see that
	the $\mathrm{Aut} \cO^{1|n}$-torsor $\mathrm{Aut}_X$ is locally trivial
	in the Zariski topology (cf. \cite[$\S$5.4.2]{frenkelzvi}). 
\end{nolabel}
\begin{nolabel}
	Similarly, let $X$ be a (oriented) superconformal $N=n$ supercurve
 	and $x \in X$. Let $\mathrm{Aut}^\omega_x$ be the set of SUSY
	coordinates $Z$ at $x$ (that is, compatible with the superconformal
	structure). It follows that this set is an $\mathrm{Aut}^\omega
	\cO^{1|n}$-torsor. Moreover these torsors glue to form an
	$\mathrm{Aut}^\omega \cO^{1|n}$-torsor $\mathrm{Aut}^\omega_X \rightarrow X$. As
	in the previous paragraph, $\mathrm{Aut}^\omega_X$ is an open
	sub superscheme of $JX^\omega$ (cf. \ref{rem:superconformaljets})
	consisting of jets of maps $D^{1|n} \rightarrow X$ compatible with the
	superconformal structure and with invertible 1-jet. 
	\label{no:4.1.8}
\end{nolabel}
\begin{rem}
	Let $V$ a finite rank  $\mathrm{Aut}
	\cO$-module (resp. a finite rank $\mathrm{Aut}^\omega \cO$-module), and let
	$X$ be an $N=n$ supercurve (resp. a superconformal $N=n$ supercurve). We
	define a vector bundle on $X$ by
	\begin{equation*}
		\cV_X = \mathrm{Aut}_X \overset{\mathrm{Aut} \cO}{\times} V
		\qquad (\text{resp. } \mathrm{Aut}^\omega_X
		\overset{\mathrm{Aut}^\omega \cO}{\times} V),
		\label{}
	\end{equation*}
	consisting of pairs $(\tilde x, v)$ with $\tilde x$ in
	$\mathrm{Aut}_X$ (resp. $\mathrm{Aut}^\omega_X$) and $v \in V$
	with the identification $(\tilde x \cdot g, v) \sim (\tilde x, g
	\cdot v)$ for $g \in \mathrm{Aut} \cO$ (resp. $g \in
	\mathrm{Aut}^\omega \cO$ ). We call $\cV_X$ the
	$\mathrm{Aut}_X$ (resp. $\mathrm{Aut}^\omega_X$) \emph{twist of $V$}. 
	\label{rem:4.1.9}
\end{rem}

\subsection{Vector bundles, sections and connections}\label{sec:vector_bundles}
\begin{nolabel}
	In this section we construct vector bundles on supercurves
	associated with SUSY vertex algebras following
	\cite{frenkelzvi}. Briefly,  a strongly conformal 
	 $N_W=n$ SUSY vertex algebra is a module for the Harish-Chandra
	pair $(\mathrm{Der}
	\cO^{1|n}, \mathrm{Aut} \cO^{1|n})$, therefore we can apply the
	Beilinson-Bernstein
	localization construction \cite{BB} to get a vector bundle with a flat
	connection over any $N=n$ supercurve. Similarly, a strongly conformal $N_K=n$ 
	SUSY vertex algebra $(n \leq 4)$, is a module
	for the Harish-Chandra pair $(\mathrm{Der}^\omega \cO^{1|n},
	\mathrm{Aut}^\omega \cO^{1|n})$\footnote{From now on, we will abuse notation and denote
	by $\mathrm{Aut}^\omega \cO^{1|n}$ its identity component}, therefore we can construct vector
	bundles with flat connections over any oriented superconformal $N=n$
	curve. 

	As in \cite{frenkelzvi}, it turns out that the state-field
	correspondence in all these cases can be seen as a (local) section
	of the corresponding bundles. The corresponding \emph{change of
	coordinates} formula (a generalization of Huang's formula
	\cite{huang} ) is proved in this section. 
\end{nolabel}
\begin{nolabel}
	Let $V$ be a strongly conformal $N_W=n$ SUSY vertex algebra.
	Therefore we have $N+1$ vectors $\nu$ and $\tau^1, \dots, \tau^N$ such
	that their Fourier modes $\nu_{(m, I)}$ and $\tau^j_{(m,I)}$ with $m \geq
	0$ generate a Lie superalgebra
	isomorphic to $\mathrm{Der} \cO^{1|n}$.
	The derivation $\partial_z$ (corresponding to
	$\nu_{(0,0)}$) cannot be
	exponentiated to the group $\mathrm{Aut} \cO^{1|n}$ and the Lie
	superalgebra spanned by $\nu_{(m,I)}$, and $\tau^j_{(m,I)}$
	for $m \geq 1$ if $I \neq 0$ is isomorphic to $\mathrm{Der}_0
	\cO^{1|n}$. 

	In order to exponentiate the representation $V$ of
	$\mathrm{Der}_0 \cO^{1|n}$ to a representation of the group
	$\mathrm{Aut} \cO^{1|n}$ we note as before that this Lie
	algebra is a semidirect product of $\fg\fl (1|n)$ with the
	pro-nilpotent Lie subalgebra $\mathrm{Der}_+ \cO^{1|n}$. Namely, the
	subalgebra spanned by $z \partial_z$, $\theta^i
	\partial_{\theta^j}$, $z \partial_{\theta^i}$ and $\theta^j
	\partial_z$ is isomorphic to $\fg\fl(1|n)$. It follows from the definition of
	strongly conformal $N_W=N$ SUSY vertex algebras in
	\ref{defn:w.conformal.definition}, that
	we can
	exponentiate this representation of $\fg\fl (1|n)$ (the fact that the
	nilpotent part of the Lie superalgebra exponentiates follows easily from the OPE
	formula and the locality axiom). 
\end{nolabel}
\begin{nolabel} \label{no:rhozdefin}
	Let $X$ be an $N=n$ supercurve over a Grassman algebra $\Lambda$, let $x
	\in X$ and
	$\cO_x$ be the completion of the local super-ring at $x$. Let
	$Z=(z,\theta^i)$ be local coordinates at $x$ (recall that in the
	formal setting, $Z$ is an \'{e}tale map $X \rightarrow
	\mathbb{A}^{1|n}$). With such a choice of coordinates we get an
	isomorphism $\cO_x \equiv \Lambda[ [Z]]$, and the set of
	coordinates at $x$, $\mathrm{Aut}_x$, is an $\mathrm{Aut}
	\cO^{1|n}$-torsor. Let us work in the analytic setting
	first for the sake of simplicity as in \cite{frenkelzvi}. Let
	$D_x$ be a small disk around $x$. Let $p$ 
	be
	a $\Lambda$-point given in the local coordinates $Z=(z,\theta^i)$
	by  $Y=(y, \alpha^i)$. The coordinates $Z$ induce
	coordinates $Z-Y = (z-y, \theta^i - \alpha^i)$ at $p$. Now let $\rho
	\in \mathrm{Aut} \cO^{1|n}$ be a change of coordinates. Recall
	that this change of coordinates is given by power series
	  $(F(Z),
	\Psi^i(Z))$, where $F(Z) \in \Lambda[ [Z]]$ is even and $\Psi^i \in
	\Lambda[ [Z]]$ are odd.
	This change of coordinates induce new coordinates at $p$, given by:
	\begin{equation}
		\rho(Z)-\rho(Y) = (F(Z) - F(Y), \Psi^i (Z) -
		\Psi^i
		(Y)).
		\label{eq:newcoordatp}
	\end{equation}
	 The coordinates $Z-Y=(z-y, \theta^i -
	\alpha^i)$ and (\ref{eq:newcoordatp}) at $p$ are related by a
	change of coordinates $\rho_{Y}=(F_{Y},
	\Psi^i_{Y})$ satisfying: 
	\begin{equation*}
		\rho_Y (Z-Y) = \rho(Z) - \rho(Y).
	\end{equation*}
	Therefore, letting  $W = (w, \zeta^i) = Z-Y$, we get:
	\begin{equation}
		\rho_Y(W) = \rho (W+Y) - \rho(Y).
		\label{eq:rhoz}
	\end{equation}

	In the formal setting we can not consider a small disk, but
	given a point $x$ and coordinates $Z$ at $x$, we can
	still define $\rho_{Z} \in \mathrm{Aut} \cO^{1|n}$ for
	any $\rho \in \mathrm{Aut} \cO^{1|n}$ by formula (\ref{eq:rhoz}) with $Y$
	replaced by $Z$. 

	Let $V$ be a strongly conformal $N_W=n$ SUSY vertex algebra, so that
	$V$ is an $\mathrm{Aut} \cO^{1|n}$-module. We will call this
	representation $R$. 
\end{nolabel}
\begin{thm} \label{thm:changewn}
	Let $V$ be a strongly conformal $N_W = n$ SUSY vertex algebra.
	  Let
	$\rho = (F, \Psi^j) \in \mathrm{Aut} \cO^{1|n}$ and $a \in V$.
	The following
	\emph{change of
	coordinates} formula is true:
	\begin{equation}
		\ys(a, Z) = R(\rho) \ys \left( R(\rho_{Z
		})^{-1} a, \rho(Z)
		\right) R(\rho)^{-1}
		\label{eq:changewn}
	\end{equation}
	where by $\rho(Z)$ we understand the images of $z,
	\theta^j$ under $\rho$, namely $F(z, \theta^i), \, \Psi^j(z,
	\theta^i)$. 
\end{thm}
\begin{proof}
	The proof is similar to the analogous formula in the ordinary vertex
	algebra
	case. Namely, the state-field correspondence $\ys(\cdot, Z)$
	is an element in the vector space $\Hom(V, \cF(V))$,
	where $\cF(V)$ is the space of all $\End(V)$-valued superfields.
	For each $\rho \in \mathrm{Aut} \cO^{1|n}$ consider the linear
	operator in $\Hom(V, \cF(V))$ given by
	\begin{equation*}
		(T_\rho X)(a, Z) = R(\rho) X (R(\rho_{Z})^{-1} a,
		\rho(Z)) R(\rho)^{-1}. 
		\label{}
	\end{equation*}
	It is easy to check that $T_\rho X \in \Hom(V, \cF(V))$.
	Moreover, this action defines a representation of $\mathrm{Aut}
	\cO^{1|n}$ in $\Hom(V, \cF(V))$. Recall that the group structure
	in $\mathrm{Aut} \cO^{1|n}$ is given by composition, namely, if
	$\rho = (F, \Psi^j)$ and $\tau = (G, \Theta^j)$ then $\rho \star
	\tau$ is given by $H, \Sigma^j$ where
	\begin{equation*}
			H(z, \theta^j) = G(F(z, \theta^j)),
			\Psi^k(z,\theta^j), \qquad
			\Sigma^i(z, \theta^j) = \Theta^i(F(z,\theta^j),
			\Psi^k(z,\theta^j)).
	\end{equation*}
	It follows that 
	$	\rho_{Z} \star \tau_{\rho(Z)} =
		(\rho \star \tau)_{Z}$. 
	Indeed, the LHS, when evaluated in $W$ is given by
	\begin{equation*}
		\tau_{\rho(Z)} \left( \rho(W+Z) - \rho(Z) \right) = \tau\left(
		\rho(W+Z) - \rho(Z) + \rho(Z) \right) - \tau(\rho(Z)),
	\end{equation*}
	which is the RHS.

	It follows from this formula that $\rho \mapsto
	T_\rho$ defines a representation of $\mathrm{Aut} \cO^{1|n}$.
	In fact, we have:
	\begin{multline*}
		(T_{\rho \star \tau}X)(a,Z) = R(\rho \star \tau)
		X(R( (\rho \star \tau)_{Z})^{-1} a,
		\tau(\rho(Z))) R(\rho \star \tau)^{-1} = \\  
		R(\rho) R(\tau) X(R(\rho_{Z} \star
		\tau_{\rho(Z)})^{-1} a, \tau(\rho(Z)))
		R(\tau)^{-1} R(\rho)^{-1} = \\ 
		R(\rho) \left[R(\tau) X(R(\tau_{\rho(Z)})^{-1}
		R(\rho_{Z})^{-1} a, \tau(\rho(Z)))
		R(\tau)^{-1}  \right] R(\rho)^{-1}= 
		\bigl[T_\rho (T_\tau X)](a,Z).
		\label{}
	\end{multline*}
	We have reduced the proof of the theorem to show that
	$\ys(\cdot, Z)$ is fixed under this action.
	Since the exponential map $\exp: \mathrm{Der}_0 \cO^{1|n}
	\rightarrow \mathrm{Aut} \cO^{1|n}$ is surjective, we need only
	to show that $\ys(\cdot, Z)$ is stable under the
	induced
	infinitesimal action of $\mathrm{Der}_0 \cO^{1|n}$. For this we
	let $\rho =\exp(\varepsilon \mathbf{v})$, where $\mathbf{v} = v(Z) \partial_Z
	\in \mathrm{Der}_0 \cO^{1|n}$, $v(Z) = (f(Z), g^1(Z),
	\dots, g^n(Z))$
	with $f(Z)$ an even function and $g^i(Z)$ odd functions of $Z$. As
	before,
	$\partial_Z = (\partial_z, \partial_{\theta^1}, \dots, \partial_{\theta^{n}})$ and the product $v(Z)
	\partial_Z$ denotes the scalar product $f(Z) \partial_z + \sum_{i = 1}^n
	g^i(Z) \partial_{\theta^i}$. 
	We want to compute 
	$\rho_{Z}$.  For this we put $\rho_Z = \exp( \varepsilon \mathbf{u})$.
	Expanding $\rho_Z(W)$ in powers of $\varepsilon$, we get
	\begin{equation}
		\mathbf{u} = v(Z + W) \partial_W - v(Z)\partial_W =\left( e^{Z \partial_W}
		v(W)\right) \partial_W - v(Z) \partial_W.
		\label{eq:alpedo2}
	\end{equation}
	Noting that the operators corresponding to $\partial_W = (\partial_w,
	\partial_{\zeta^1}, \dots, \partial_{\zeta^n})$ are $- \nabla = (-T, -S^1,
	\dots, -S^n)$, we obtain:
	\begin{equation}
		R(\mathbf{u}) = e^{-Z \nabla} R(\mathbf{v}) e^{Z\nabla} + v(Z)
		\nabla. 
		\label{eq:lo_que_faltaba}
	\end{equation}
%
	
	The (infinitesimal) action of $T_\rho$ on
	$\ys(a, Z)$ is given by $\ys(a, Z)$ plus the
	linear term in $\varepsilon$, which in turn is:
	\begin{equation*}
		[R(\mathbf{v}), \ys(a, Z)
		] - \ys(R(\mathbf{u}) a, Z) + v(Z)
		\nabla_Z \ys(a,Z).
	\end{equation*}
	The first term comes from the adjoint action of $R(\rho)$, the
	second term is the $\varepsilon$-linear term in
	$R(\rho_{Z})^{-1}$, and the last term comes from the
	Taylor expansion of the change of coordinates. The result follows from 
	(\ref{eq:alpedo2}), (\ref{eq:lo_que_faltaba}), Proposition \ref{prop:commutator_w.1}  and Theorem
	\ref{thm:structure.2}.
\end{proof}
\begin{nolabel}
	Now we can define a vector bundle associated to an $N_W=n$ SUSY 
	vertex algebra over any $N=n$ supercurve. Moreover, we will
	define a canonical section of this bundle and a flat connection
	on it. First recall that from any finite dimensional
	$\mathrm{Aut} \cO^{1|n}$-module we can construct a vector bundle
	over an $N=n$ supercurve $X$ by twisting this $\mathrm{Aut} \cO^{1|n}$-module by the
	$\mathrm{Aut} \cO^{1|n}$-torsor $\mathrm{Aut}_X$ (see Remark
	\ref{rem:4.1.9}). Given a
	strongly conformal $N_W = n$ SUSY vertex algebra $V$, we have a
	filtration $V_{\leq i}$ by finite dimensional submodules, namely,
	$V_{\leq i}$ is the span of fields of conformal weight
	less or equal than $i$. By our assumptions, these are finite
	dimensional $\mathrm{Aut} \cO$-submodules of $V$. Let $\cV_{\leq
	i}$ be the corresponding $\mathrm{Aut}_X$ twist. These vector bundles come
	equipped with embeddings $\cV_{\leq i} \hookrightarrow \cV_{\leq
	i + 1}$. The limit of this directed system is a $\cO_X$-module
	$\cV_X$\footnote{When there is no possible confusion, we will denote this bundle simply
	by $\cV$.}:
	\begin{equation*}
		\cV_X = \varinjlim_{i \rightarrow \infty} \cV_{\leq i}.
		\label{}
	\end{equation*}
	This $\cO_X$-module is quasi-coherent by definition.
	
	On the other hand, the dual modules $V_{\leq i}^*$ come equipped
	with surjections $V_{\leq i+1}^* \twoheadrightarrow V_{\leq i}^*$
	therefore we get a projective system of $\cO_X$-modules
	$\cV_{\leq i+1}^* \twoheadrightarrow \cV^*_{\leq i}$. The inverse
	limit of this system is by definition $\cV^*_X$, namely:
	\begin{equation*}
		\cV^*_X = \varprojlim_{i \rightarrow \infty} V_{\leq i}^*.
		\label{}
	\end{equation*}

	Thus, we have defined $\cO_X$-modules associated with the SUSY
	vertex algebra $V$. We will call these modules the
	\emph{SUSY vertex algebra bundle} and its dual. By construction, the fiber of the
	bundle $\cV$ at a
	point $x \in X$ is isomorphic as a vector space, to $V$. 

	Similar constructions can be applied when $X$ is replaced by a
	formal superdisk near a point $x \in X$. Namely, let $D_x$ be such a
	formal superdisk, we have as before an $\mathrm{Aut} \cO^{1|n}$-torsor
	$\mathrm{Aut}_{D_x}$ over $D_x$. Then $\cV_{D_x}$ is the twist of
	$V$ by this torsor. It is easy to see that in this case we get
	$\cV_X|_{D_x} = \cV_{D_x}$.

	 Let $\mathrm{Aut}_x$ be the torsor of coordinates at
	$x$ as before. Then the fiber of $\cV$ at $x$ is given by:
	\begin{equation*}
		\cV_x = \mathrm{Aut}_x \overset{\mathrm{Aut} \cO}{\times}
		V.
		\label{}
	\end{equation*}
	Let $D_x^\times$ be the punctured disk at $x$, that is the formal
	completion:
	\begin{equation*}
		D_x^\times = \varinjlim_{i \rightarrow \infty}
		\Spec(\tilde{\cK}_x/\fm^{i+1}),
		\label{}
	\end{equation*}
	where $\tilde{\cK}_x$ is the ring of fractions of the local ring
	at $x$ and $\fm$ is the maximal ideal defining $x$. If
	$Z=(z,\theta^i)$ are coordinates at $x$, this is isomorphic to
	the formal spectrum of $\Lambda( (Z))$. 

	We will define an $\End{\cV_x}$-valued section of $\cV^*$ on
	$D_x^\times$. In order to define such a section it is enough to
	give its matrix coefficients, namely, for each $\varphi \in
	\cV^*_x$, $v \in \cV_x$ and $s$ a section of $\cV|_{D_x}$ we
	assign a function on $D_x^\times$, that is an element of $\cK_x$,
	the ring of fractions of $\cO_x$. This assignment is denoted by:
	\begin{equation*}
		\varphi, v, s \mapsto <\varphi, \cys_x(s) \cdot v>,
		\label{}
	\end{equation*}
	and should be linear in $v$ and $\varphi$ and $\cO_x$ linear in
	$s$.
	Let $Z=(z, \theta^i)$ be coordinates at $x$, we obtain a
	trivialization 
	$i_{Z}: V[ [Z]] \risom \Gamma(D_x, \cV)$, of $\cV|_{D_x}$.
	This induces isomorphisms $V \risom \cV_x$ and $V^* \risom
	\cV_x^*$, where $V^*$ is the \emph{restricted dual} of $V$. Let $v \in V$ and $\varphi
	\in V^*$. Denote their images in $\cV_x$ and $\cV_x^*$, under these
	isomorphisms   by
	$(Z, v)$ and $(Z, \varphi)$ respectively. Let $s \in V[ [Z]]$,
	 its image under the isomorphism $i_{Z}$ is a
	regular section of $\cV$ in $D_x$. By $\cO_x$ linearity, we may assume
	that $s = a \in V$. To this data, we assign the
	function:
	\begin{equation}
		<(Z,\varphi), \cys_x(i_Z(a)) \cdot (Z, v)> = <\varphi,
		\ys(a, Z) v>.
		\label{eq:ysdefinition}
	\end{equation}
\end{nolabel}
\begin{thm}
	The assignment (\ref{eq:ysdefinition}) is independent of the
	coordinates $Z=(z,\theta^1, \dots, \theta^n)$ chosen, i.e. $\cys_x$ is a well defined
	$\End(\cV_x)$-valued section of $\cV^*$ on $D_x^\times$.
	\label{thm:assignment.1}	
\end{thm}
\begin{proof}
	The proof follows the lines of the ordinary vertex algebra case in \cite{frenkelzvi}. Let
	$W=(w,\zeta^i)$ be another set of coordinates at $x$. Then $W$
	and $Z$ are related by $\rho \in \mathrm{Aut} \cO$, $\rho(Z) =
	W$. Given these new coordinates, we construct another assignment by
	the same formula (\ref{eq:ysdefinition}), namely
	\begin{equation*}
		<(W,\varphi), \tilde{\cys}(i_W(a)) \cdot(W,v)> =
		<\varphi, Y(a, W)v>.
		\label{}
	\end{equation*}
	We need to show that this assignment coincides with $\cys$. By
	the definition of the bundle $\cV$ we have
	\begin{equation*}
		(Z,v) = (\rho^{-1}(W), v) = (W, R(\rho)^{-1} v),
		\label{}
	\end{equation*}
	where $R(\cdot)$ is the representation of $\mathrm{Aut}
	\cO^{1|n}$ in $V$. Similarly $(Z,\varphi) = (W, \varphi
	R(\rho))$. We need to find how does the section $i_Z(a)$
	transform by this change of coordinates. Recall from
	\ref{no:rhozdefin} that in the analytic setting, if we trivialize
	$\cV|_{D_x}$ with the coordinates $Z$, we can use the coordinates
	$(Z-Y):=(z-y, \theta^i - \alpha^i)$ at $Y=(y, \alpha^i)$ to identify
	$\cV_y$ with $V$. We obtain:
	\begin{equation}
		(Z-Y, a) = (W - \rho(Y), R(\rho_{Y})^{-1} a),
		\label{eq:acabo.de.agregarte}
	\end{equation}
	therefore the section $i_Z(a)$ is 
	$i_W(R(\rho_{Z})^{-1}a)$ in the $W$-trivialization. 

	In the formal setting, we can replace the coordinates by their
	$n$-jets, but these in turn can be extended by definition to a
	small Zariski open neighborhood of $x$, in this case, the formula
	(\ref{eq:acabo.de.agregarte})
	is true as we have shown.
	
	We have reduced the problem to prove:
	\begin{equation*}
		<\varphi, R(\rho)\ys(R(\rho_{Z})^{-1}a, W)
		 R(\rho)^{-1} v> = <\varphi, Y(a, Z) v>,
		\label{}
	\end{equation*}
	thus, the theorem follows from Theorem \ref{thm:changewn}
\end{proof}
\begin{nolabel}
	In the superconformal case, the
	situation is slightly more
	complicated. Roughly, the only changes that we have to make in the
	above prescription are the induced coordinates at a
	$\Lambda$-point and consequently the definition of
	$\rho_{Z}$. 

	Like in the $N_K=n$ SUSY vertex algebra situation, given two set of coordinates $Z =
	(z,\theta^1, \dots,\theta^n )$ and $W = (w,\zeta^1, \dots, \zeta^n)$ we will write \[Z-W
	= \Bigl(z-w-\sum_{i=1}^n
	\theta^i \zeta^i, \theta^1 - \zeta^1, \dots, \theta^n - \zeta^n\Bigr).\]

	 Let $V$
	be a strongly conformal $N_K=n$ SUSY vertex algebra $(n \leq 4)$, hence  $V$
	is an $\mathrm{Aut}^\omega \cO^{1|n}$-module.
	Moreover, $V$ has a filtration by finite dimensional
	submodules $V_{\leq i}$ given by conformal weight as above. Let
	$X$ be an oriented superconformal $N=n$ supercurve over $\Lambda$. We
	constructed an $\mathrm{Aut}^\omega \cO$-torsor $\mathrm{Aut}^\omega_X$ over
	$X$ (see \ref{no:4.1.8}). As above we can define the vertex algebra bundles $\cV$
	and $\cV^*$. Similarly, we can define the $N_K=n$ SUSY vertex algebra
	bundles over the superconformal disks $D^\omega_x$. The fibers
	$\cV_x$ of
	these bundles are the $\mathrm{Aut}_x^\omega$-twists of $V$, where
	$\mathrm{Aut}_x^\omega$ is the torsor of coordinates at $x$, compatible
	with the superconformal structure (see Remark \ref{rem:4.1.9}). We define and
	$\End(\cV_x)$-valued section $\cys_x$ of $\cV^*$ on the punctured
	disk $D_x^\times$ by formula (\ref{eq:ysdefinition}).
	\label{no:k.bundle.def}
\end{nolabel}
\begin{thm}
	The assignment $\cys_x$ is independent of the coordinates $Z=(z, \theta^i)$
	 chosen as long as they are compatible with the
	superconformal structure on $X$.
	\label{thm:6.8}
\end{thm}
\begin{proof}
	Let us first work in the analytic setting. If $p$ is a
	$\Lambda$-point in $D_x$ (now a small analytic disk near $x \in
	X$) given by local parameters $Y=(y, \alpha^i)$, then $Z$
	induces local coordinates at $T=(t,\eta^i)= Z-Y$ near $p$. The
	coordinates $T$ are
	compatible with the superconformal structure. Indeed, we have:
	\begin{equation*}
			dt = dz + \sum_{i=1}^n \alpha^i d\theta^i, \qquad
			d\eta^i = d\theta^i,
		\label{}
	\end{equation*}
	therefore:
	\begin{equation*}
		dt + \sum_{i=1}^n \eta^i d\eta^i = d z + \sum_{i=1}^n \alpha^i
		d\theta^i + (\theta^i + \alpha^i) d\theta^i = dz + \sum_{i=1}^n
		\theta^i d \theta^i.
		\label{}
	\end{equation*}
	If $W=(w,\zeta^i)=\rho(Z)$ is another set of
	coordinates compatible with the superconformal structure at x,
	with $\rho = (F, \Psi^i) \in \mathrm{Aut}^\omega \cO^{1|n}$, then $W$
	induces another
	set of coordinates at $p$, namely 
	\begin{equation*}
		\rho(Z)-\rho(Y)=\left( F(z,\theta) - F(y, \alpha) - \sum_{i=1}^n \Psi^i(z,\theta)
		\Psi^i(y,\alpha), \Psi^j(z,\theta)-
		\Psi^j(y,\alpha)\right).
		\label{}
	\end{equation*}
	These are related with the coordinates $T$ by a change of
	coordinates 
	\begin{equation*}
		\rho_{Y} = (F_{Y},
	\Psi^i_{Y}) \in \mathrm{Aut}^\omega \cO^{1|n}.
	\end{equation*}
	We have:
	\begin{equation*}
		\rho_Y (T) = \rho(T + Y) - \rho(Y),
	\end{equation*}
	where, as in the $N_K=n$ SUSY vertex algebra case, we write $T + Y = T - (-Y)$.

	The theorem will follow if we
	prove formula (\ref{eq:changewn}) for $\rho \in \mathrm{Aut}^\omega
	\cO^{1|n}$. This is achieved as in the proof of Theorem
	\ref{thm:changewn} by first showing that the action $T_\rho$
	 is a representation of
	$\mathrm{Aut}^\omega \cO^{1|n}$ in $\Hom(V, \cF(V))$. For this we first note
	that $(\rho \star \tau)_{Z} =
	\rho_{Z} \star \tau_{\rho(Z)}$ in exactly the same way as in the $N_W=n$
	case. Again we just
	have to prove that $\ys$ is fixed under this action, and
	we check this at the level of Lie algebras. Denote $D_W = (\partial_w, 
	D^1_W, \dots, D^n_W)$, where $D^i_W = \partial_{\zeta^i} + \zeta^i \partial_w$. Similarly, denote $\bar{D}_W =
	(\partial_w, \bar{D}^1_W, \dots, \bar{D}^n_W)$, where $\bar{D}^i_W = \partial_{\zeta^i} - \zeta^i \partial_w$. Let $\rho =
	\exp( \varepsilon \mathbf{v})$ where $\mathbf{v} = v(W) \bar{D}_W \in
	\mathrm{Der}_0^\omega \cO^{1|n}$, and put $\rho_Z = \exp(
	\varepsilon \mathbf{u})$. Expanding $\rho_Z$ in powers of
	$\varepsilon$ we find:
	\begin{equation*}
		\mathbf{u} = v(W + Z) \bar{D}_W - v(Z) \bar{D}_W.
		\label{}
	\end{equation*}
	Note that in this context we have two
	\emph{different} Taylor expansions\footnote{Note that $Z+ W \neq W + Z$.}:
	\begin{equation*}
		e^{Z D_W} f(W) = f(Z + W), \qquad e^{Z \bar{D}_W} f(W) = f(W + Z),
		\label{}
	\end{equation*}
	using the second, we see that 
	\begin{equation*}
		\mathbf{u} = \left( e^{Z \bar{D}_W} v(W) \right) \bar{D}_W - v(Z) \bar{D}_W.
		\label{}
	\end{equation*}
	From this and the fact that the operators corresponding to $\bar{D}_W$
	are:
	\[-\nabla = (-T, -S^1, \dots, -S^n),\] we obtain:
	\begin{equation*}
		R(\mathbf{u}) = e^{-Z \nabla} R(\mathbf{v}) e^{Z \nabla} + v(Z) \nabla.
		\label{}
	\end{equation*}
	The theorem now follows as in the $N_W = n$ case.
\end{proof}
Now we construct connections on the vector bundles $\cV$ from the
previous paragraphs. 
\begin{thm}
	Let $X$ be a $(1|N)$ dimensional supercurve. Let $U \subset X$ be open
	and $Z$ be coordinates in $U$ defining the vector fields
	$\partial_z$ and $\partial_{\theta^i}$. Let $V$ be a strongly conformal
	$N_W=N$ SUSY vertex algebra and $\cV$ the associated bundle. Define
	the connection operators $\nabla_\chi : \cV_{|U} \rightarrow \cV_{|U}$
	for each  vector field $\chi$ in $U$ by
	\begin{equation*}
		\nabla_{\partial_z} = \partial_z + T,  \qquad
		\nabla_{\partial_{\theta^i}} = \partial_{\theta^i} + S^i.
	\end{equation*}
	Then $\nabla$ is a well defined (left) connection on $\cV$ (independent
	of the coordinates chosen). Moreover, this connection is flat. 
	\label{thm:connection.w}
\end{thm}
\begin{proof}
	The proof is verbatim the proof of the analogous statement in
	\cite[16.1]{frenkelzvi}. Indeed, strongly conformal SUSY vertex algebras are
	modules for the Harish Chandra pair ($\mathrm{Der} \cO^{1|N},\,
	\mathrm{Aut} \cO^{1|N}$) and this in turn acts simply transitively on the
	torsor $\mathrm{Aut}_X \rightarrow X$. The localization procedure of
	formal geometry applies without difficulties. 
\end{proof}
\begin{rem}
	Note that this connection endows $\cV$ with a structure of a left
	$\cD_X$-module for any supercurve $X$ and any strongly conformal
	$N_W=N$ 
	SUSY vertex algebra $V$. 

	Let $V$ be a strongly conformal $N_K=N$ SUSY
	vertex algebra, and let $\cV$ be the associated vector bundle over an oriented
	superconformal curve $X$.  For an open $U$ as before, and superconformal
	coordinates $Z$ in $U$ we will define the superconformal differential
	operators $\cD_X(U)$ to be the super ring of differential operators generated by
	all the $D^i_Z$. This defines a sheaf of algebras of superconformal
	differential operators
	$\cD_X$ over any (oriented) superconformal curve $X$. The assignment
	\begin{equation}
		D^i_Z \cdot f(Z) a = (D^i_Z f(Z)) a + (-1)^{f} f(Z) S^i a 
		\label{eq:rem.connection.1}
	\end{equation}
	gives $\cV$ the structure of a left $\cD_X$-module.
	\label{rem:connection}
\end{rem}
\subsection{Examples}\label{sec:7}
\begin{nolabel}
	\label{no:7.1} In this section we give the first non-trivial examples of
	the super vector bundles that arise with the construction of the previous
	sections. To simplify the notation, we will use the ordinary
	description of the involved vertex algebras. For example, when we analyze
	the \emph{boson-fermion} system (cf. example \ref{ex:7.2}) we will work
	with the fermion $\varphi$ and the boson $\alpha$ instead of the 
	superfields $\Psi$ and $S\Psi$. Note that the Grassman algebra $\Lambda$
	is a SUSY vertex algebra (either $N_W=N$ or $N_K=N$) with $T = S^i = 0$ and
	$\vac = 1$. In this section, given a SUSY vertex algebra $V$, we will
	consider the tensor product $W = \Lambda \otimes V$ (either of $N_W=N$ or
	$N_K=N$ SUSY vertex algebras), therefore we can view $W$ as a SUSY vertex
	algebra over $\Lambda$, namely, $W$ is a $\Lambda$-module and the vertex
	operators are $\Lambda$-linear.

	Let us start with $N_K=1$ bundles. For this let $X$ be a super
	conformal $N=1$ supercurve over $\Lambda$. Let $U_\alpha$ and $U_\beta$ be
	open in $X$ 
	and $p=(t,\zeta)$ a $\Lambda$-point in the intersection. Let $V$ be a 
	strongly conformal $N_K=1$ SUSY vertex algebra, so that $V$ carries a
	representation of $Der_0^\omega \cO^{1|1}$ that exponentiates to a
	representation of $\mathrm{Aut}^\omega \cO^{1|1}$.
	Suppose we have coordinates
	$(z_\alpha, \theta_\alpha)$ in $U_\alpha$ and $(z_\beta,\theta_\beta)$ in
	$U_\beta$ that are compatible with the superconformal structure.
            They are related by a change of coordinates $\rho_{\beta\alpha}
	= (F(z_\alpha, \theta_\alpha), \Psi(z_\alpha, \theta_\alpha))$ satisfying
	$DF = \Psi D \Psi$ where $D = \partial_{\theta_\alpha} + \theta_\alpha
	\partial_{z_\alpha}$. These
	coordinates define coordinates at the point $p$, therefore we obtain
	different trivializations of the bundle $\cV$. The transition functions for
	the structure sheaf give us transition functions for $\cV$, in particular,
	they act in the fiber at the point $p$ as $R(\rho_p)^{-1}$ (cf.
	\ref{eq:acabo.de.agregarte}). 

	In order to compute $R(\rho_p)$ we need to look only at the odd
	coordinate, namely, expand in Taylor series:
	\begin{equation}
		\begin{aligned}
			\Psi_{z,\theta}(t,\zeta) &= \Psi(t + z + \zeta \theta,
			\zeta + \theta) - \Psi(z,\theta) \\
			&= \zeta D \Psi + t D^2 \Psi + \zeta t D^3 \Psi +
			\frac{t^2}{2} D^4 \Psi + \dots \\ 
			&= \exp \left( - \sum_{i \geq 1} (v_i L_i + w_i G_{(i)})
			\right) A^{-2 L_0} \cdot \zeta,
			\end{aligned}
		\label{eq:7.1.1}
	\end{equation}
	where as in 
	(\ref{eq:neveu_shwarz_generators}) we have
	\begin{equation*}
			L_n = - \frac{n+1}{2} t^n \zeta \partial_\zeta
			- t^{n+1} \partial_t, \quad
			G_{(n+1/2)}=G_n = - t^{n+1/2} (\partial_\zeta - \zeta
			\partial_t),
	\end{equation*}
	and $v_i = v_i(z,\theta)$ are even functions and $w_i=w_i(z,\theta)$ are
	odd functions. Truncating the series in (\ref{eq:7.1.1}) at order $2$ we have:
	\begin{equation*}
		\Psi_{(z,\theta)}(t,\zeta) = A\left( \zeta + t w_1 + \zeta t v_1 +
		t^2 (w_2 + v_1 w_1) + \dots \right),
	\end{equation*}
	from where we get the equations:
	\begin{equation}
	\begin{aligned}
		A &= D\Psi, & w_1 A &= D^2\Psi, \\
		v_1 A &= D^3 \Psi, & (w_2 + v_1 w_1)A &= \frac{1}{2} D^4 \Psi.
	\end{aligned}
		\label{eq:7.1.4}
	\end{equation}
	We can solve this system to get:
	\begin{equation}
		\begin{gathered}
			v_1 = \frac{D^3 \Psi}{D \Psi}, \qquad
			w_1 = 	\frac{D^2 \Psi}{D \Psi}, \\
			w_2 = \frac{1}{2} \left( \frac{D^4 \Psi}{D \Psi} - 2
			\frac{D^3\Psi D^2 \Psi}{(D\Psi)^2} \right) =
			\frac{1}{2} \sigma \bigl(D\Psi\bigr),
		\end{gathered}
		\label{eq:7.1.5}
	\end{equation}
	where $\sigma$ is the \emph{$N=1$ super-schwarzian} defined in
	(\ref{eq:schwarzian_1_def}).  
\end{nolabel}
\begin{ex}[Free Fields]
	\label{ex:7.2} Recall the strongly conformal $N_K=1$ SUSY vertex algebra $B_1$ defined
	in Example \ref{ex:2.16}. We will
	denote this vertex algebra as $B(1)$. As an
	ordinary vertex algebra, it  is graded with respect to
	conformal weight. The fermion $\varphi$ is primary of conformal weight
	$1/2$ and the boson $\alpha$ has conformal weight $1$ but it is not primary
	unless $m = 0$. It follows easily from the non-commutative Wick formula,
	that the only
	non-trivial relations with the fermion $\varphi$ are given by:
	\begin{equation*}
		G_{(1)} \varphi = - m \vac, \qquad L_0 \varphi = \frac{1}{2}
		\varphi,
	\end{equation*}
	 therefore the subspace $B(1)_{\leq 1/2}$ of $B(1)$
	  spanned by
	$\{ \vac, \varphi\}$ 
	is an $\mathrm{Aut}^\omega \cO^{1|1}$ submodule. For a given change of
	coordinates $\rho = (F,\Psi)$ we can compute the action of
	$R(\rho_{(z,\theta)})^{-1}$. For this we write in the basis
	$\{\vac,\varphi\}$: 
	\begin{equation}
		\begin{aligned}
			R(\rho_{(z,\theta)})^{-1} &= A^{2 L_0} \cdot \exp \left(
			\sum_{i \geq
		1} (v_i
		L_i + w_i G_{(i)}) \right) \\
		&= \begin{pmatrix}
			1 & - m w_1 \\
			0 & A 
		\end{pmatrix}
		= \begin{pmatrix}
			1 & - m \frac{D^2 \Psi}{D \Psi} \\
			0 & D\Psi
		\end{pmatrix}.
	\end{aligned}
		\label{eq:7.2.2}
	\end{equation}
	Hence, if $\cB(1)$ is the vector bundle associated to the
	$\mathrm{Aut}^\omega
	\cO^{1|1}$ module $B(1)$ and $\cB(1)_{\leq 1/2}$ is the
	vector bundle corresponding to $B(1)_{\leq 1/2}$ we see that the transition functions
	that define $\cB(1)_{\leq 1/2}$ are given on the intersections $U_\alpha
	\cap U_\beta$ by the functions (\ref{eq:7.2.2}). 

	Dually, sections of the bundle $\cB(1)_{\leq 1/2}^*$ transform by (note
	that we use the \emph{super-transpose} instead of the transpose, as defined
	in \cite[$\S$ 3.1]{manin2}):
	\begin{equation}
		\begin{pmatrix}
			1 & 0 \\
			m \frac{D^2 \Psi}{D \Psi} & D\Psi 
		\end{pmatrix}.
		\label{eq:7.2.3}
	\end{equation}
	Note that we have a section $\cys$ of $\cB(1)^*$ which projects to a
	section of $\cB(1)_{\leq 1/2}^*$. In the basis $\{\vac, \varphi\}$ this
	section is given by:
	\begin{equation*}
		\begin{pmatrix}
			\mathrm{Id} \\
			\varphi(z, \theta)
		\end{pmatrix},
	\end{equation*}
	where, according to (\ref{eq:nimepreguntes}), we have 
	\begin{equation*} 
		\varphi(z,\theta) =
		Y(\varphi_{(-1)} \vac, z) + \theta Y(G_{-1/2}
		\varphi_{(-1)} \vac, z).
	\end{equation*}
	According to (\ref{eq:7.2.3}) and Theorem \ref{thm:6.8} we see that
	the field $\varphi(z,\theta)$ transforms as:
	\begin{equation}
		\varphi(z,\theta) = R(\rho) \varphi(\rho(z,\theta)) R(\rho)^{-1}
		D\Psi + m \frac{D^2 \Psi}{D \Psi} \mathrm{Id},
		\label{eq:7.2.5}
	\end{equation}
	where $\rho=(F,\Psi)$. In particular, since $X$ is a superconformal $N=1$
	curve we have:
	\begin{equation*}
		D\Psi = D\left( \frac{DF}{D\Psi} \right) = \mathrm{sdet}
		\begin{pmatrix}
			F_z & \Psi_z \\
			F_\theta & \Psi_\theta
		\end{pmatrix}.
	\end{equation*}
	Therefore when $m = 0$, $\varphi(z,\theta)[dz d\theta]$ transforms as
	an $\End \cB(1)_p$-valued section of the Berezinian bundle of $X$ on the
	punctured disk $D^\times_p$ for any $\Lambda$-point $p \in X$. When $m \neq 0$
	this bundle is
	not split and $\varphi(z,\theta)$ gives rise to an $\End \cB(1)_p$-valued
	section of $\cB(1)^*_{\leq 1/2}$ that projects onto the section $1 \otimes
	\mathrm{Id}$ of the quotient $\cO_X \otimes \End \cB(1)_p$ and transforms
	according to (\ref{eq:7.2.5}) with changes of coordinates.
	In other words, the bundle $\cB_{\leq 1/2}(1)^*$ is an
	extension:
	\begin{equation}
		0 \rightarrow \ber_X \rightarrow \cB_{\leq 1/2} (1)^* \rightarrow
		\cO_X \rightarrow 0
		\label{eq:7.2.6p}
	\end{equation}
	which is non-split unless $m = 0$. In the case when $m \neq 0$ the
	section $\cys$ projects into the constant section $1$ of $\cO_X$\footnote{From now one, we abuse notation and forget the fact that
	the section $\cys$ is  $\End(\cV_x)$-valued.}.

	In analogy to \cite{frenkelzvi} we want to understand the geometric meaning
	of these sections. Equivalently, we want to find the set of splittings of the
	extension (\ref{eq:7.2.6p}). This set, if non-empty, is a torsor over the
	space of even sections of $\mathrm{Ber}_X$. Recall also that the operator $D =
	\partial_\theta + \theta \partial_z$ takes values in 
	$\ber_X$ for a superconformal $N=1$ curve. 
\end{ex}
\begin{thm}
	\label{thm:7.3} The superfield $\varphi(z,\theta)$ transforms as an odd differential
	operator $\nabla : \ber_X \rightarrow \ber_X^{\otimes 2}$ locally of the
	form $\nabla = - m D_\alpha + g_\alpha(z_\alpha,\theta_\alpha)$, where on the
	open subset $U_\alpha$ with coordinates $(z_\alpha, \theta_\alpha)$ we have
	$D_\alpha = \partial_{\theta_\alpha} + \theta_\alpha \partial_{z_\alpha}$
	and $g_\alpha$ is an odd function.  
\end{thm}
\begin{proof}
	Recall that in a superconformal $N=1$ curve the generators $[dz_\alpha
	d\theta_\alpha]$ of the Berezinian bundle transform as:
	\begin{equation*}
		[dz_\beta d\theta_\beta] = (D_\alpha \Psi_{\beta,\alpha} )
		[dz_\alpha d\theta_\alpha],
	\end{equation*}
	where the change of coordinates is $\theta_\beta = \Psi_{\beta,
	\alpha}(z_\alpha, \theta_\alpha)$. 
	Since $\nabla :\ber_X \rightarrow \ber_X^{\otimes 2}$ we have:
	\begin{equation*}
		\nabla_\alpha f_\alpha = (D_\alpha \Psi_{\beta, \alpha})^2
		\nabla_\beta \left( (D_\alpha \Psi_{\beta,\alpha})^{-1} f_\alpha
		\right).
	\end{equation*}
	Therefore we get:
	\begin{equation*}
		\begin{aligned}
			\nabla_\alpha f_\alpha &= - m D_\alpha f_\alpha +
			g_\alpha f_\alpha = -m (D_\alpha \Psi_{\beta,\alpha})
			D_\beta f_\alpha + g_\alpha f_\alpha \\
			&= (D_\alpha \Psi_{\beta,\alpha})^2 \left( -m \left(
			D_\beta \left( D_\alpha \Psi_{\beta,\alpha} \right)^{-1}
			f_\alpha \right) + g_\beta \left( D_\alpha
			\Psi_{\beta,\alpha} \right)^{-1} f_\alpha   \right) \\
			&= - m \left( D_\alpha \Psi_{\beta,\alpha} \right)
			D_\beta f_\alpha - m \left( D_\alpha \Psi_{\beta,\alpha}
			\right)^2 D_\beta \left(  \left( D_\alpha
			\Psi_{\beta,\alpha} \right)^{-1} \right) f_\alpha +\\
			&\quad + g_\beta 
			\left( D_\alpha \Psi_{\beta,\alpha} \right) f_\alpha \\
			&= - m \left( D_\alpha \Psi_{\beta,\alpha} \right)
			D_\beta f_\alpha + m \left( D_\alpha \Psi_{\beta,\alpha}
			\right)^{-1} D_\alpha^2 \Psi_{\beta,\alpha} f_\alpha +
			g_\beta \left( D_\alpha \Psi_{\beta,\alpha} \right)
			f_\alpha,\\
			g_\alpha &= \left( D_\alpha \Psi_{\beta,\alpha} \right)
			g_\beta + m \frac{D_\alpha^2
			\Psi_{\beta,\alpha}}{D_\alpha \Psi_{\beta,\alpha}}.
		\end{aligned}
	\end{equation*}
	Hence we find: 
	\begin{equation*}
		\begin{pmatrix}
			1 \\ g_\alpha
		\end{pmatrix} = 
		\begin{pmatrix}
			1 & 0\\ m \frac{D^2 \Psi}{D \Psi} & D \Psi
		\end{pmatrix}
		\begin{pmatrix}
			1 \\ g_\beta
		\end{pmatrix}.
	\end{equation*}
	thus proving the theorem.
\end{proof}
\begin{nolabel}
	\label{no:7.4} Given that we can integrate a section of $\ber_{X}$ along a
	super contour as in \ref{no:3.16}, we can state \cite[7.1.9]{frenkelzvi} in
	this situation. We define an \emph{affine structure} on a superconformal
	$N=1$ curve to be a (equivalence class of) coordinate atlas $U_\alpha$ with
	coordinates 
	$(z_\alpha, \theta_\alpha)$ such that the transition functions on overlaps
	satisfy\footnote{These are SUSY changes of coordinates where the odd
	coordinate changes by affine transformations.}:
	\begin{equation}
		\begin{aligned}
		z_\beta &= F_{\beta,\alpha}(z_\alpha, \theta_\alpha) = a^2
		z_\alpha + \theta_\alpha \xi a + b, \\
		\theta_\beta &= \Psi_{\beta,\alpha}(z_\alpha,\theta_\alpha) =
		\theta_\alpha a + \xi,
	\end{aligned}
		\label{eq:7.4.1}
	\end{equation}
	where $a, b$ are even constants with $a$ invertible and $\xi$ is an odd
	constant (these constants may change with $\alpha$ and $\beta$). Given such
	an atlas, we can define $\nabla_\alpha = - m D_\alpha$ and we get from:
	\begin{equation*}
		- m D_\alpha = - m D_\alpha \Psi_{\beta,\alpha} D_\beta,
	\end{equation*}
	and the fact that $D^2 \Psi = 0$ for these transition functions, that
	$\nabla_{\alpha}$ is a well defined operator as in Theorem \ref{thm:7.3}.
	
	On the other hand, suppose we have such a differential operator
	$\nabla_\alpha = - m D_\alpha + g_\alpha$.
	Let $f_\alpha [dz_\alpha d\theta_\alpha]$ be a section of
	$\ber_X$ in $U_\alpha$ such that $f_\alpha$ is an even function and
	$\nabla_\alpha \cdot f_\alpha = 0$. Choose a $\Lambda$-point $P=(x,\pi)$
	of $U_\alpha$ and, for any other point $Q$ in $U_\alpha$, we define the
	function $\xi_\alpha$ to be
	\begin{equation*}
		\xi_\alpha (Q) = \int_P^Q f_\alpha.
	\end{equation*}
	From the definition of this integral we see that $\xi$ is an odd function,
	indeed, to compute this integral we need to solve $D \omega = f$ and then
	this integral becomes $\omega(Q) - \omega(P)$. By shrinking if necessary
	the open cover $U_\alpha$ we may assume that $f_\alpha$ does not vanish
	everywhere (it is an even function), it follows that $D \xi$ is
	invertible everywhere. We now solve the differential equation $D w = \xi
	D\xi$ (we may need to shrink $U_\alpha$ even more) and obtain thus a
	coordinate atlas $U_\alpha$ with new coordinates $(w_\alpha, \xi_\alpha)$.
	We claim that this atlas is indeed an affine structure on $X$. We have made
	some choices. One is the reference point $P$ which shifts the function
	$\xi_\alpha$ by an odd constant. 
 	The other choice was the solution $\xi$, which is
	unique up to an invertible even multiple (for this we can apply a version
	of Cauchy's theorem in super geometry). Therefore $\xi$ is well defined up
	to affine transformations of the form $\xi \mapsto a \xi + \zeta$. This
	forces $w$ to change to $\tilde{w}$ with 
	\begin{equation*}
		D \tilde{w} = (a \xi + \zeta) D(a \xi + \zeta) = a^2 \xi D\xi + a
		\zeta D\xi,
	\end{equation*}
	hence $\tilde w = a^2 w + w'$ with $D w' = D a \zeta \xi$. Finally we see
	that $w' = a \zeta \xi + w''$ where $D w'' = 0$, namely the choices made
	combine into changes of the form:
	\begin{equation*}
			\xi \mapsto a \xi + \zeta,  \qquad
			w \mapsto a^2 w + a \zeta \xi + b,
	\end{equation*}
	where $a, b$ are even constants ($a$ is invertible) and $\zeta$ is odd.
	Since these changes of coordinates are of the form
	(\ref{eq:7.4.1}),  we have proved:
\end{nolabel}
\begin{thm}
	\label{thm:7.5} Let $X$ be an $N=1$ superconformal curve. 
	For every $m \neq 0$ the set of differential operators
	$\nabla : \ber_X \rightarrow \ber_X^{\otimes 2}$ locally defined as
	$\nabla_\alpha = - m D_\alpha + g_\alpha$ for odd functions $g_\alpha$
	are in one to one correspondence with the set of affine structures on the
	curve $X$. These in turn are in one to one correspondence with the set of
	splittings of the extension (\ref{eq:7.2.6p}). 
\end{thm}
\begin{ex}
	\label{ex:7.6} \textbf{The Neveu-Schwarz vertex algebra} Recall the strongly
	conformal $N_K=1$ SUSY vertex algebra $K_1$ defined in Example
	\ref{ex:ns.example} (see also Example \ref{ex:k.n.series}). Denote this vertex
	algebra by $K(1)$.
	Note that the sub-vector space spanned by the primary elements of
	conformal weight less or equal to $3/2$, namely the vacuum vector and the
	$N=1$ superconformal vector $\tau$, is $\mathrm{Aut}^\omega \cO^{1|1}$-invariant. 
	 In order to
	compute the transition functions we use the relevant
	relations in this case:
	\begin{equation*}
		L_{(1)}\tau = \frac{3}{2} \tau,  \qquad G_{(2)} \tau =
		\frac{2}{3} c.
	\end{equation*}
	Therefore we can compute in the basis $\{\vac, \tau\}$:
	\begin{equation*}
		\begin{aligned}
			R(\rho_{(z,\theta}))^{-1} \begin{pmatrix} \vac \\ \tau
			\end{pmatrix} &= A^{2 L_0} \exp \left(  \sum_{i \geq 1}
		v_i L_i + w_i G_{(i)} \right) \begin{pmatrix}
			\vac \\ \tau
		\end{pmatrix}\\
		&= \begin{pmatrix}
			1 & \frac{2}{3} c w_2 \\ 0 & A^3 	
		\end{pmatrix}\begin{pmatrix}
			\vac \\ \tau
		\end{pmatrix}.
	\end{aligned}
	\end{equation*}
	It follows from (\ref{eq:7.1.5}) and (\ref{eq:7.1.4}) that the transition
	functions for the
	corresponding bundle $\cK_{\leq 3/2}(1)$ are given by:
	\begin{equation}
		R(\rho_{(z,\theta)})^{-1} =\begin{pmatrix}
			1 & \frac{c}{3} \sigma \bigl(D\Psi\bigr) \\ 0 & (D\Psi)^3
		\end{pmatrix},
		\label{eq:7.6.3}
	\end{equation}
	where, as before,  $\sigma \bigl(D \Psi\bigr)$ is the \emph{super
	shwarzian} derivative. Dualizing, we obtain an
	extension:
	\begin{equation}
		0 \rightarrow \ber_X^{\otimes 3} \rightarrow \cK_{\leq 3/2}(1)^*
		\rightarrow \cO_{X} \rightarrow 0.
		\label{eq:7.6.4}
	\end{equation}
	This extension is not split if $c \neq 0$ and, as for the free fields, we
	see that the section $\cys$ of $\cK_{\leq 3/2}(1)^*$ projects onto the
	section $1 \in \cO_X$
	in this case. Denote by $\tau(z,\theta) = G(z) + 2 \theta L(z)$ the
	superfield $\ys(\tau, z,\theta)$. By taking the super transpose of
	(\ref{eq:7.6.3}) we find that $\tau(z,\theta)$ transforms as:
	\begin{equation*}
		\tau(z,\theta) = R(\rho) \tau(\rho(z,\theta)) R(\rho)^{-1}
		(D\Psi)^3 - \frac{c}{3} \sigma \bigl(D\Psi\bigr)  
	\end{equation*}
	which in turn implies, according to Proposition \ref{prop:3.22} the
	following:
\end{ex}
\begin{thm}
	 The
	set of splittings of (\ref{eq:7.6.4}) is in one to one correspondence with
	the set of superprojective structures in $X$.
	\label{thm:agregadoeste}
\end{thm}
\begin{nolabel}
	\label{no:7.7} Now we turn our attention to the oriented superconformal
	$N=2$ case.  We will use the coordinates
	$(z, \theta^\pm = \theta^1 \pm i \theta^2)$  and the change of coordinates $\rho =
	(F, \Psi^\pm = \Psi^1 \pm i\Psi^2)$ (cf.
	\ref{no:n=2supergroup}). It follows that:
	\begin{multline*}
		\Psi^\pm_{(z,\theta^+,\theta^-)}(t,\zeta^+,\zeta^-) = \Psi^\pm
		\left( t + z + \frac{1}{2} (\zeta^+ \theta^- + \zeta^- \theta^+),
		\zeta^+ + \theta^+, \zeta^- + \theta^- \right) - \\ -
		\Psi^\pm(z,\theta^+,\theta^-),
	\end{multline*}
	which we want to expand in Taylor series (here $\Psi$ denotes
	either $\Psi^+$ or $\Psi^-$):
	\begin{equation*}
		\begin{aligned}
		\Psi_{(z,\theta^\pm)} &= \left( 1 
			 + \frac{1}{2} (\zeta^+ \theta^- + \zeta^-
			 \theta^+)\partial_z
			  + \frac{1}{8}(\zeta^+ \theta^- + \zeta^-
			 \theta^+)^2 \partial^2_z \right) \cdot \\ & \qquad \qquad
			 \cdot \Psi(t +
			 z, \zeta^+ +
			 \theta^+, \zeta^- + \theta^-) - \Psi(z,\theta^+,\theta^-)
			 \\
			 &= \left( 1 
			 + \frac{1}{2} (\zeta^+ \theta^- + \zeta^-
			 \theta^+)\partial_z
			  + \frac{1}{4} \zeta^+ \zeta^- \theta^+ \theta^-
			  \partial_z^2
			 \right) \cdot \\ & \qquad \qquad
			 \cdot \Psi(t +
			 z, \zeta^+ +
			 \theta^+, \zeta^- + \theta^-) - \Psi(z,\theta^+,\theta^-)\\
			 &= \biggl[ \zeta^+ \left( \partial_{\theta^+} +
			 \frac{1}{2} \theta^- \partial_z \right) + \zeta^-
			 \left( \partial_{\theta^-} + \frac{1}{2} \theta^+
			 \partial_z
			 \right)  + \\
			 & \quad + t \partial_z  + \zeta^+ t \left(
			 \partial_{\theta^+} + \frac{1}{2} \theta^- \partial_z
			 \right) \partial_z +  \zeta^- t \left( \partial_{\theta^-} +
			 \frac{1}{2} \theta^+ \partial_z \right)\partial_z + \\ &\quad +
			 \zeta^+ \zeta^- \left( \partial_{\theta^-, \theta^+} -
			 \frac{1}{2} \theta^- \partial_{z,\theta^-} +
			 \frac{1}{2} \theta^+ \partial_{z,\theta^+} +
			 \frac{1}{4} \theta^+ \theta^- \partial_z^2 \right) +
			 \frac{1}{2} t^2 \partial^2_z \biggr] \Psi +
			 \dots \\
			 &= \Bigl( \zeta^+ D^-  + \zeta^- D^+  + t \partial_z 
			 + \zeta^+ t D^- \partial_z + \zeta^- t D^+ \partial_z + \\
			 & \qquad \qquad +\zeta^+
			 \zeta^- (D^+ D^-
			 - \frac{1}{2} \partial_z) + \frac{1}{2} t^2 \partial^2_z
			 \Bigr) \Psi + \dots, 
		 \end{aligned}	 
	\end{equation*}
	where $D^\pm = \partial_{\theta^\mp} + \tfrac{1}{2} \theta^\pm \partial_z$. 
	Since the curve is oriented (cf. (\ref{eq:7.7.8.a})), this reduces to:
	\begin{equation*}
		\begin{aligned}
			\Psi^+_{(z,\theta^\pm)} &= \Bigl( \zeta^+ D^-  
			+ t D^+D^- 
			 + \zeta^+ t D^- D^+ D^-  + \\
			 & \qquad \qquad +\zeta^+
			 \zeta^- \frac{1}{2}(D^+ D^-
				) + \frac{1}{2} t^2 (D^+ D^-)^2
			 \Bigr) \Psi^+ + \dots \\	
			 \Psi^-_{(z,\theta^\pm)} &= \Bigl( \zeta^- D^+  
			+ t D^-D^+ 
			 + \zeta^- t D^+ D^- D^+  + \\
			 & \qquad \qquad +\zeta^-
			 \zeta^+ \frac{1}{2}(D^- D^+
				) + \frac{1}{2} t^2 (D^- D^+)^2
			 \Bigr) \Psi^- + \dots \\
		\end{aligned}
	\end{equation*}
	We want to express these as the exponential of a vector field. For this we
	compute:
	\begin{multline*}
		\exp\left( -\sum_{i \geq 1} v_i L_i + u_i J_i + w^\pm G^\pm_{(i)}
		\right) B^{-J_0} A^{-2 L_0} \cdot \zeta^\pm = B^{\pm 1} A
		\biggl[\zeta^\pm + t w_1^\pm + \\ +\zeta^\pm t (v_1 \pm u_1 +
		\frac{1}{2} w_1^\mp w_1^\pm)  + 
		t^2 (w_2^\pm + \frac{1}{2}w_1^\pm(2
		v_1 \pm u_1)) + \frac{1}{2} \zeta^\pm \zeta^\mp w_1^\pm \biggr] +
		\dots,
	\end{multline*}
	where we have used (\ref{eq:7.7.9.a}). 
	We get the equations:
	\begin{equation*}
		\begin{gathered}
			B^{\pm 1} A = D^\mp \Psi^\pm, \qquad
			w_1^\pm = \frac{D^\pm D^\mp \Psi^\pm}{D^\mp \Psi^\pm} =
			\frac{\Psi^\pm_z}{D^\mp \Psi^\pm} =
			\frac{\Psi^\pm_z}{\Psi^\pm_{\theta^\pm}}, \\
			v_1 \pm u_1 + \frac{1}{2} w_1^\pm w_1^\mp =
			\frac{D^\mp D^\pm D^\mp \Psi^\pm}{D^\mp \Psi^\pm}, \quad
			w_2^\pm + \frac{1}{2} w_1^\pm (2 v_1 \pm u_1) =
			\frac{1}{2} \frac{(D^\pm D^\mp)^2 \Psi^\pm}{D^\mp
			\Psi^\pm}.
		\end{gathered}
	\end{equation*}
	We can solve this system to get:
	\begin{equation*}
		\begin{gathered}
			v_1 = \frac{1}{2} \left( \frac{D^- \Psi^+_z}{D^- \Psi^+} +
			\frac{D^+ \Psi^-_z}{D^+ \Psi^-}\right), \quad
			w_2^\pm = \frac{1}{2 D^\mp \Psi^\pm} \left( \Psi^\pm_{z,z}
			- \frac{1}{2} \left( \frac{D^\pm \Psi^\mp_z}{D^\pm
			\Psi^\mp} + 3 \frac{D^\mp \Psi^\pm_z}{D^\mp \Psi^\pm} \right)
			\right),\\
			u_1 = \frac{1}{2}\left( \frac{D^- \Psi^+_z}{D^- \Psi^+} -
			\frac{D^+ \Psi^-_z}{D^+ \Psi^-}\right) -
			\frac{1}{2} \frac{\Psi^+_z \Psi^-_z}{D^-\Psi^+
			D^+\Psi^-} = - \sigma_2 (\Psi^+, \Psi^-),
		\end{gathered}
	\end{equation*}
	where $\sigma_2$ is the $N=2$ \emph{schwarzian derivative} (cf.
	\cite{cohn}). 
\end{nolabel}
\begin{ex}\textbf{Free Fields.}
	\label{ex:7.8} With the results of
	the previous sections we can compute now
	explicitly some vector bundles over oriented superconformal $N=2$ curves.
	Let $Y$ be such a curve and let $B(2)$ be the strongly conformal $N_K=2$ SUSY vertex
	algebra described in Example \ref{ex:2.23}. Let $\cB(2)$ be the associated
	vector bundle over $Y$. The vector subspace spanned by the vacuum vectors
	and the two fermions (namely the fields with conformal weight less or equal
	to $1/2$) is an $\mathrm{Aut}^\omega \cO^{1|2}$-submodule. Let us denote these
	vectors, as in \ref{ex:2.23}, by $\{\vac, \varphi^\pm\}$ respectively, and let
	$\cB_{\leq
	1/2}(2)$ be the associated rank $1|2$ vector bundle over $Y$. In order to
	compute its transition functions explicitly we note that the only
	nontrivial relations (for our purposes) are:
	\begin{equation*}
		G^\pm_{(1)} \varphi^\mp = \mp m \vac, \qquad J_0 \varphi^\pm = \pm
		\varphi^\pm, \qquad L_0 \varphi^\pm = \frac{1}{2} \varphi^\pm.
	\end{equation*}
	We compute the transition functions as:
	\begin{equation}
		\begin{aligned}
			R(\rho)^{-1}\begin{pmatrix}
				\vac \\ \varphi^+ \\ \varphi^-
			\end{pmatrix} &= A^{2 L0} B^{J_0} \exp \left( \sum_{i \geq
			1} v_1 L_i + u_i J_i + w^\pm_i G^\pm_{(i)}
			\right) \begin{pmatrix}
				\vac \\ \varphi^+ \\ \varphi^-
			\end{pmatrix} \\
			&=\begin{pmatrix}
				1 &  m w_1^- & - m w_1^+ \\
				0 & B A & 0 \\
				0 & 0 & B^{-1} A
			\end{pmatrix}\begin{pmatrix}
				\vac \\ \varphi^+ \\ \varphi^-
			\end{pmatrix} \\
			&=\begin{pmatrix}
				1 &  m \frac{\Psi^-_z}{D^+ \Psi^-} & - m
				\frac{\Psi^+_z}{D^- \Psi^+} \\
				0 & D^- \Psi^+ & 0 \\
				0 & 0 & D^+ \Psi^-
			\end{pmatrix}\begin{pmatrix}
				\vac \\ \varphi^+ \\ \varphi^-
			\end{pmatrix}.
		\end{aligned}
		\label{eq:7.8.2}
	\end{equation}
	Recall now that an oriented superconformal  $N=2$ curve projects onto two
	$N=1$ supercurves: $X$ and its dual $\hat{X}$ (cf.
	\ref{no:n=2superconformal}). Using the coordinates (cf.
	\ref{no:n=2supergroup}): \[\left(u = z +
	\frac{1}{2} \theta^+ \theta^-, \theta^+, \theta^- \right),\] we obtain from
	(\ref{eq:7.8.4.a}) and (\ref{eq:7.8.3'.a}) that 
	\begin{multline}
			D^+ \Psi^- = D^+ \left( \frac{D^- G}{D^- \Psi^+ } \right)
			= D^+ \left(
			\frac{1}{(\Psi^+_{\theta^+})^2}(\Psi^+_{\theta^+} -
			\theta^- \Psi^+_u)(G_{\theta^+} + \theta^- G_u \right)  
			 = \\ = \frac{\Psi^+_{\theta^+} G_u - \Psi^+_u
			G_{\theta^+}}{(\Psi^+_{\theta^+})^2} 
			= \mathrm{sdet } \begin{pmatrix}
				G_u & \Psi^+_u \\ G_{\theta^+} & \Psi^+_{\theta^+}
			\end{pmatrix},
		\label{eq:7.8.5}
	\end{multline}
	where $G = F + \tfrac{1}{2}\Psi^+ \Psi^-$ as in \ref{no:n=2supergroup}.
	Similarly, we find
	\begin{equation}
		D^- \Psi^+ = \sdet \begin{pmatrix}
			G'_{u'} & \Psi^-_{u'} \\
			G'_{\theta^-} & \Psi^-_{\theta^-}
		\end{pmatrix}.
		\label{eq:7.8.6}
	\end{equation}
	Let us call $\pi$ and $\hat{\pi}$ the projections from $Y$ onto $X$ and
	$\hat{X}$ respectively. 
	We see from (\ref{eq:7.8.6}) and (\ref{eq:7.8.5}) that taking the
	super-transpose in (\ref{eq:7.8.2}) we obtain an extension (of sheaves of
	$\cO_Y$-modules):
	\begin{equation}
		0 \rightarrow \pi^* \ber_X \oplus \hat{\pi}^* \ber_{\hat{X}}
		\rightarrow \cB(2)_{\leq 1/2}^* \rightarrow \cO_{Y} \rightarrow 0.
		\label{eq:7.8.7}
	\end{equation}
	As in the $B(1)$ case, this extension is not split unless $m$ vanishes.
	It follows in the same way as in the $N=1$ case that  the set of splittings
	of this extension corresponds to \emph{affine 
	structures} on the $N=2$ superconformal curve $Y$. Indeed, we see in the
	same way as in Theorem \ref{thm:7.3}, that the
	pair of fields $(\varphi^{+},\, \varphi^{-})$ transforms as 
	a differential operator $\nabla : \ber_{\hat{X}} \oplus \ber_X \rightarrow
	\ber_{\hat{X}}^{\otimes 2} \oplus \ber_{X}^{\otimes 2}$ which is locally
	of the form $(m D^+ + g^+, - m D^- + g^-)$ for $g^\pm$ odd functions of
	$(u,\theta^+)$ and $(u', \theta^-)$ respectively. We note that according to
	\ref{no:3.16} sections of $\ber_X \oplus \ber_{\hat{X}}$ can be integrated
	in $Y$ up to an additive constant. The argument in the proof of Theorem
	\ref{thm:7.5} generalizes to this setting without difficulty. 

	We will return to this example below (cf. \ref{ex:7.11}).
\end{ex}
\begin{ex}\textbf{The $N=2$ vertex algebra.}
	\label{ex:7.9} Let $K(2):= K_2$ be the strongly conformal $N_K=2$ SUSY vertex algebra described in
	Example  \ref{ex:n=2} (see also Example \ref{ex:k.n.series}), and let $\cK(2)$ be the associated vector bundle
	over an oriented superconformal $N=2$ curve $Y$. The vector subspace spanned by primary fields of conformal weight
	$0$ and $1$ is an $\mathrm{Aut}^\omega \cO^{1|2}$ submodule. Let us denote
	these
	vectors as above by $\{\vac, J\}$ respectively, and let $\cK(2)_{\leq 1}$
	be the associated rank $2|0$ vector bundle over $Y$. To compute the
	transition functions we note that the only non-trivial 
	relations we need are:
	\begin{equation}
		L_0 J = J, \qquad J_1 J = \frac{c}{3} \vac.
		\label{eq:7.9.1}
	\end{equation}
	Therefore the transition functions are given by:
	\begin{equation*}
		\begin{aligned}
			R(\rho)^{-1}\begin{pmatrix}
				\vac \\ J
			\end{pmatrix} &=\begin{pmatrix}
				1 & \frac{c}{3} u_1 \\ 
				0 & A^2
			\end{pmatrix}\begin{pmatrix}
				\vac \\ J
			\end{pmatrix} \\ &=
			\begin{pmatrix}
				1 & - \frac{c}{3} \sigma_2( \Psi^+, \Psi^-) \\
				0 & D^+\Psi^- D^-\Psi^+
			\end{pmatrix}\begin{pmatrix}
				\vac \\ J
			\end{pmatrix}.
		\end{aligned}
	\end{equation*}
	It follows as before, by taking the super-transpose, that when $c=0$, the superfield
	$J(z,\theta^+, \theta^-)$ transforms as a section of $\pi^*\ber_X \otimes
	\hat{\pi}^*\ber_{\hat{X}}$,
	namely in this case we get an extension:
	\begin{equation}
		0 \rightarrow \pi^*\ber_X \otimes \hat{\pi}^*\ber_{\hat{X}}
		\rightarrow \cK(2)_{\leq 1}^* \rightarrow \cO_Y \rightarrow 0 ,
		\label{eq:7.9.3}
	\end{equation}
	which is split if and only if $c = 0$. When $c \neq 0$, 
	 the superfield $J(z,\theta^+, \theta^-)$ transforms
	as:
	\begin{equation*}
		J(z,\theta^+,\theta^-) = (D^+ \Psi^-)(D^- \Psi^+)
		J(\rho(z,\theta^+,\theta^-)) + \frac{c}{3} \sigma_2 (\Psi^+,
		\Psi^-). 
	\end{equation*}
	We see that the
	section $\cys$ is an even section projecting onto $1 \in \cO_Y$, therefore
	giving a splitting of (\ref{eq:7.9.3}). The set of such splittings, if
	non-empty, is a torsor for the even part of $\pi^*\ber_X \otimes
	\hat{\pi}^*\ber_{\hat{X}}$.

	Analyzing this algebra further, we can consider the space $K(2)_{\leq 3/2}$
	spanned by vectors
	of conformal weight less than or equal to $3/2$.
	This space is spanned by $\{\vac, J, G^\pm\}$. In addition to (\ref{eq:7.9.1})
	we have the following relations:
	\begin{equation*}
		\begin{aligned}
			L_0 G^- &= \frac{3}{2} G^-, & J_0 G^- &= - G^-, & G^+_{(1)}
			G^- &= J, &  G^+_{(2)} G^- &= \frac{c}{3} \vac, \\
			L_0 G^+ &= \frac{3}{2} G^+, & J_0 G^+ &= G^+, & G^-_{(1)} G^+
			&= - J, & G^-_{(2)}G^+ &= \frac{c}{3} \vac.
		\end{aligned}
	\end{equation*}
	With these we can compute the transition functions in the basis
	$\{\vac, J, G^-, G^+\}$ explicitly:
	\begin{equation}
		\begin{aligned}
		R(\rho)^{-1} &=\begin{pmatrix}
			1 & \frac{c}{3} u_1 & \frac{c}{3} w_2^+ & \frac{c}{3} w_2^-
			\\
			0 & A^2 & A^2 w_1^+ & - A^2 w_1^- \\
			0 & 0 & A^3 B^{-1} & 0  \\
			0 & 0 & 0 & A^3 B 
		\end{pmatrix},
	\end{aligned}
		\label{eq:7.9.6}
	\end{equation}
	the first three by three block being:
	\begin{equation*}
		 \begin{pmatrix}
			1 & - \frac{c}{3} \sigma_2( \Psi^+, \Psi^-) &
			\frac{c}{6 D^- \Psi^+} \left( \Psi^+_{z,z}
			- \frac{1}{2} \left( \frac{D^+ \Psi^-_z}{D^+
			\Psi^-} + 3 \frac{D^- \Psi^+_z}{D^- \Psi^+} \right)
			\right) \\
			0 & (D^+\Psi^-)(D^-\Psi^+) & (D^+\Psi^-) \Psi^+_z \\
			0 & 0 & (D^- \Psi^+)(D^+ \Psi^-)^2 
		\end{pmatrix},
	\end{equation*}
	and the $4,4$ entry in (\ref{eq:7.9.6}) is $(D^+ \Psi^-)(D^-\Psi^+)^2$.
	Taking the super-transpose of (\ref{eq:7.9.6}) it follows that
	$\cK(2)_{Y,\leq 3/2}^*$ fits in a short exact sequence of the form:
	\begin{equation*}
		0 \rightarrow \pi^* \ber_X \otimes \left(\hat{\pi}^*
		\ber_{\hat{X}}\right)^{\otimes 2} \rightarrow \cK(2)_{Y, \leq
		3/2}^*
		\rightarrow \cN^* \rightarrow 0 .
	\end{equation*}
	The bundle $\cN$ in turn fits in the exact sequence:
	\begin{equation*}
		0 \rightarrow \left( \pi^* \ber_X \right)^{\otimes 2} \otimes
		\hat{\pi}^*\ber_{\hat{X}} \rightarrow \cN^*
		\rightarrow \cK(2)_{Y,\leq 1}^* \rightarrow 0.
	\end{equation*}
	In a more ``symmetric'' fashion, if we look at the lower two by two
	block in (\ref{eq:7.9.6}), we see that we have an extension:
	\begin{multline*}
		0 \rightarrow \pi^* \ber_{X} \otimes \left( \hat{\pi}^*
		\ber_{\hat{X}} \right)^{\otimes 2} \oplus \left( \pi^* \ber_X 
		\right)^{\otimes 2} \otimes \hat{\pi}^* \ber_{\hat{X}} \rightarrow \\ 
		\rightarrow \cK(2)^*_{\leq 3/2}  \rightarrow \cK(2)_{\leq 1}^*
		\rightarrow 0.
	\end{multline*}
\end{ex}
\begin{nolabel}
	We turn our attention now to the $N_W=1$ case. For this let $X$ be a
	general $N=1$ supercurve. As before, given a change of coordinates
	$\rho=(F,\Psi)$, we expand in Taylor series:
	\begin{equation*}
		\begin{aligned}
			F_{(z,\theta)}(t,\zeta) &= F(t + z, \zeta + \theta) - F(z,
			\theta)\\ &= t F_z + \zeta F_\theta + \zeta t F_{\theta,z}
			+ \frac{t^2}{2} F_{z,z} + \dots \\
			\Psi_{(z,\theta)}(t,\zeta) &= \Psi(t+z,\zeta+\theta) -
			\Psi(z,\theta) \\
			&= t \Psi_z + \zeta \Psi_\theta + \zeta t \Psi_{z,\theta} +
			\frac{t^2}{2} \Psi_{z,z}.
		\end{aligned}
	\end{equation*}
	We need to express these as:
	\begin{multline}
		 \begin{pmatrix} F_{(z,\theta)} \\ \Psi_{(z,\theta)} \end{pmatrix}
			 = \exp \left( - \sum_{i \geq 1} v_i T_i + u_i J_i + q_i
			 Q_{i} + h_i H_i \right) \times \\ \times \exp (-q_0 Q_0)
			 \exp(- h_0 H_0)
			 B^{-J_0} A^{-T_0}\begin{pmatrix}
				 t \\ \zeta
			 \end{pmatrix},
		\label{eq:7.10.2}
	\end{multline}
	where, as in (\ref{eq:4.2.7}), we have:
	\begin{subequations}
	\begin{xalignat*}{2}
		T_n &=  -
		t^{n+1} \partial_t - (n+1) t^n \zeta \partial_\zeta, &
		J_n &= -t^n \zeta \partial_\zeta, \\
		Q_n &= -t^{n + 1} \partial_\zeta, & H_n &=
		t^{n} \zeta \partial_t.
	\end{xalignat*}
	\end{subequations}
	Expanding (\ref{eq:7.10.2}) up to second order, we find:
	\begin{equation*}
		\begin{aligned}
			F_{(z,\theta)} &= t A(1+ q_0 h_0) + \zeta A h_0 + t^2  (v_1
			(A + A q_0 h_0) + A q_1 h_0)
			+ \\
			& \quad + \zeta t (A(1+q_0 h_0) h_1 + 2 A v_1 h_0 + A u_1 h_0) 
			+ \dots \\
			\Psi_{(z,\theta)} &= \zeta B A + t q_0 B A + t\zeta B A(2
			v_1 + u_1 + h_1 q_0) + t^2 BA (q_1 + v_1 q_0) + \dots,
		\end{aligned}
	\end{equation*}
	and we get the equations:
	\begin{xalignat*}{2}
		A(1+q_0 h_0) &= F_z, & BA  &= \Psi_{\theta}, \\
		A h_0 &= F_\theta, & q_0 B A &= \Psi_z, \\
		v_1 F_z + q_1 F_\theta &= \frac{1}{2} F_{z,z}, & h_1 \Psi_z + (2 v_1
		+ u_1) \Psi_\theta &= \Psi_{\theta,z}, \\
		h_1F_z + (2 v_1 + u_1) F_\theta &= F_{z,\theta}, & v_1 \Psi_z + q_1
		\Psi_\theta&=
		\frac{1}{2}\Psi_{z,z}.
	\end{xalignat*}
	From this we find:
	\begin{equation}
	\begin{aligned}
		A &= \frac{F_z \Psi_\theta - \Psi_z F_\theta}{\Psi_\theta}, & B &=
		\frac{\Psi_\theta^2}{F_z \Psi_\theta - \Psi_z F_\theta},\\  
		h_0 &= \frac{F_\theta \Psi_\theta}{F_z \Psi_\theta - \Psi_z
		F_\theta}, &
		q_0 &= \frac{\Psi_z}{\Psi_\theta}, \\
		v_1 &= \frac{1}{2} \frac{F_{z,z} \Psi_\theta - \Psi_{z,z}
		F_\theta}{F_z \Psi_\theta - \Psi_z F_\theta}, & q_1 &=
		\frac{1}{2} \frac{F_{z,z} \Psi_z - \Psi_{z,z} F_z}{F_\theta \Psi_z
		- \Psi_\theta F_z}, \\
		h_1  &= \frac{F_{z,\theta} \Psi_\theta
		- \Psi_{z,\theta} 
		F_\theta}{F_z \Psi_\theta - \Psi_z F_\theta},  & u_1 &=
		\frac{F_{z,\theta}\Psi_z - \Psi_{z,\theta} F_z}{F_\theta \Psi_z -
		\Psi_\theta F_z} + \frac{\Psi_{z,z} F_\theta - F_{z,z}
		\Psi_\theta}{F_z \Psi_\theta - \Psi_z F_\theta}.
	\end{aligned}
		\label{eq:7.10.6}
\end{equation}
	\label{no:7.10}
\end{nolabel}
\begin{ex}
	\textbf{Free Fields} Consider the vertex algebra $B(2)$ as in example
	\ref{ex:7.8} but as a $N_W=1$ SUSY vertex algebra. As such, for each
	$N=1$ supercurve $X$ we obtain a vector bundle $\cB(2)_X$. Recall that 
	with respect to the Virasoro field $\tilde{L}$, the vector $\varphi^-$ has conformal
	weight $0$. Therefore the vector space spanned by $\vac$ and $\varphi^-$ is
	an $\mathrm{Aut} \cO^{1|1}$-submodule. We obtain then a rank $1|1$ vector
	bundle over $X$, to be denoted $\cB(2)_{X,\leq 0}$. Let us compute
	explicitly the transition functions for this bundle. The relevant
	relations are in this case:
	\begin{equation*}
		J_0 \varphi^- = - \varphi^-, \qquad Q_0 \varphi^- = -m \vac.
	\end{equation*}
	Hence we obtain:
	\begin{equation*}
		R(\rho)^{-1} \begin{pmatrix} \vac \\ \varphi^- \end{pmatrix} = A^{
			T_0} B^{J_0} \exp (h_0 H_0) \exp(q_0 Q_0)\begin{pmatrix}
				\vac \\ \varphi^-
			\end{pmatrix} = 
			\begin{pmatrix}
				1 & - m q_0 \\
				0 & B^{-1}
			\end{pmatrix}
			\begin{pmatrix}
				\vac \\ \varphi^-
			\end{pmatrix},
	\end{equation*}
	which implies:
	\begin{equation}
		R(\rho)^{-1}=\begin{pmatrix}
			1 & - m \frac{\Psi_z}{\Psi_\theta} \\
			0 & \frac{F_z \Psi_\theta - \Psi_z F_\theta}{\Psi_\theta^2}
		\end{pmatrix}.
		\label{eq:7.11.3}
	\end{equation}
	Noting that 
	\begin{equation*}
		 \sdet \begin{pmatrix}
			F_z & \Psi_z \\ F_\theta & \Psi_\theta 
		\end{pmatrix} = \frac{F_z \Psi_\theta- \Psi_z
		F_\theta}{\Psi_\theta^2},
	\end{equation*}
	we see that by taking the super-transpose in (\ref{eq:7.11.3}) we obtain an
	extension
	\begin{equation}
		0 \rightarrow \ber_X \rightarrow \cB(2)_{X,\leq 0}^* \rightarrow
		\cO_X \rightarrow 0.
		\label{eq:7.11.5}
	\end{equation}
	This short exact sequence is split if and only if $m = 0$. In that case,
	we see that $\varphi^-(z,\theta)[dzd\theta]$ transforms as a section of
	$\ber_X$. On the other hand, when $m \neq 0$, (\ref{eq:7.11.5}) is not
	split and  the section $\cys$ projects into $1 \in \cO_X$, giving a
	splitting of (\ref{eq:7.11.5}). In order to analyze the splittings of these
	sequences, recall from \ref{no:n=2superconformal} that to the $N=1$
	supercurve $X$ we have associated another ``dual'' curve $\hat{X}$ and
	an oriented superconformal $N=2$ curve $Y$. Introduce maps of sheaves
	on $Y$,  
	 $\nabla : \ber_{\hat X} \rightarrow
	\ber_{X} \otimes \ber_{\hat{X}}$ which are locally of the form
	$\nabla_\alpha = - m D^+_\alpha + g_\alpha$, for an odd function $g_\alpha
	= g_\alpha(u, \theta^+)$. Here we consider $X$ with
	coordinates $u, \theta^+$ and $\hat{X}$ with coordinates $u', \theta^-$ as
	in \ref{no:n=2supergroup}. We will write $\hat{f}$ to
	denote a function of $u', \theta^-$. It follows from (\ref{eq:7.8.6}),
	(\ref{eq:7.8.7}) and the fact that $\nabla$ maps $\ber_{\hat{X}}
	\rightarrow \ber_X \otimes \ber_{\hat{X}}$ that on overlaps we must have:
	\begin{equation*}
		\nabla_\alpha \hat{f}_\alpha = (D^+ \Psi^-) (D^- \Psi^+)
		\nabla_\beta \left( (D^- \Psi^+)^{-1} \hat{f}_\alpha \right).
	\end{equation*}
	Replacing $\nabla$ in both sides by its local form and using
	(\ref{eq:7.7.7.a}) and (\ref{eq:7.7.8.a}),
	 we get:
	\begin{multline*}
		- m D^+_\alpha \hat{f}_\alpha + g_\alpha \hat{f}_\alpha = - m
		D^+\alpha \hat{f}_\alpha + m (D^- \Psi^+)^{-1} D^+_\alpha
		D^-_\alpha \Psi^+ \hat{f}_\alpha + (D^+\Psi^-) g_\beta f_\alpha.
	\end{multline*}
	Now noting that $D^+ D^- \Psi^+ = \Psi^+_u$ and that:
	\begin{equation*}
		\frac{\Psi^+_u}{D^- \Psi^+} = \frac{\Psi^+_u}{\Psi^+_{\theta^+} +
		\theta^- \Psi^+_u} = \frac{\Psi^+_u}{\Psi^+_{\theta^+}},
	\end{equation*}
	we get
	\begin{equation}
		g_\alpha = \sdet \begin{pmatrix}
			G_u & \Psi^+_u \\ G_{\theta^+} & \Psi^+_{\theta^+}
		\end{pmatrix}g_\beta + m \frac{\Psi^+_u}{\Psi^+_{\theta^+}}
		\label{eq:7.11.5.4}
	\end{equation}
	therefore proving the following
	
	\textbf{Theorem.} The set of splittings of (\ref{eq:7.11.5}) for $m \neq
	0$ is in one to
	one correspondence with operators $\nabla : \ber_{\hat{X}} \rightarrow
	\ber_X \otimes \ber_{\hat{X}}$ locally of the form $-m D^+_\alpha +
	g_\alpha$.

	Let $\nabla$ be such an operator, and let $0 \neq \psi_\alpha \in 
	\ber_{\hat{X}}(U_\alpha)$ be a flat even section, namely $\nabla_\alpha
	\psi_\alpha = 0$. As a section of $\ber_{\hat{X}}$, it can be integrated
	along any contour in $X$ (cf. \ref{no:3.16}), namely, let $P$ be a reference $\Lambda$-point in
	$U_\alpha$, then for any other $\Lambda$-point in $U_\alpha$ we put:
	\begin{equation*}
		\zeta_\alpha(Q) = \int_P^Q \psi_\alpha.
	\end{equation*}
	The solution $\zeta_\alpha$ is unique up to an even multiplicative constant,
	while changing the reference point $P$ changes $\zeta_\alpha$ by an
	additive odd
	constant, shrinking $U_\alpha$ we may assume that $D_\alpha \zeta$ is
	invertible. Choosing any other even function $t_\alpha$ with invertible
	differential, we obtain charts $U_\alpha, (t_\alpha, \zeta_\alpha)$. The
	transition functions between these charts are clearly affine functions for
	the odd coordinates, namely $\zeta_\beta = a_{\beta,\alpha} \zeta_\alpha +
	\varepsilon_{\beta,\alpha}$ for some even constants $a$ and odd constants
	$\varepsilon$. Conversely, given such a covering of $X$, we define
	$\nabla_\alpha = - m D^+_\alpha$, where we take $\zeta$ instead of
	$\theta^+$ and $t$ instead of $u$ in the definition of $D^+$. It follows
	from (\ref{eq:7.11.5.4}) that $\nabla$ is well defined globally since the
	second term in the right hand side of (\ref{eq:7.11.5.4}) vanishes. 

	Combining the above paragraph with the previous theorem we have

	\textbf{Theorem.} The set of splittings of (\ref{eq:7.11.5}) for $m \neq 0$
	is in one to one correspondence with (equivalence classes of) atlases
	$U_\alpha, z_\alpha, \theta_\alpha$, such that the transition functions
	are affine in the odd coordinate, namely $\theta_\beta = a \theta_\alpha +
	\varepsilon$ for some even constant $a$ and some odd constant
	$\varepsilon$.

	 Note that
	from (\ref{eq:7.11.3}) and (\ref{eq:7.8.2}) it follows that the
	following sequences of $\cO_Y$-modules are exact:
	\begin{equation*}
		\begin{aligned}
			& 0 \rightarrow \hat{\pi}^* \ber_{\hat{X}} \rightarrow
			\cB(2)_{Y, \leq 1/2}^* \rightarrow \pi^* \cB(2)_{X,\leq
			0}^* \rightarrow 0, \\
			& 0 \rightarrow \pi^* \ber_X \rightarrow \cB(2)_{Y, \leq
			1/2}^* \rightarrow \hat{\pi}^*\cB(2)_{\hat{X}, \leq 0}
			\rightarrow 0.
		\end{aligned}
	\end{equation*}
	The bundle $\cB(2)_Y$ is the corresponding bundle constructed in Example
	\ref{ex:7.8} from this vertex algebra, but viewed as an $N_K=2$ SUSY vertex
	algebra. These two extensions show how the different vector bundles
	constructed from the same vertex algebras in these three different curves
	($X$, $\hat{X}$ and $Y$) are related. 	
	\label{ex:7.11}
\end{ex}
\begin{ex}
	\textbf{The $N=2$ vertex algebra}. Let, as before $K(2)$ be the $N=2$ super
	vertex algebra defined in Example \ref{ex:n=2}, but considered as an $N_W=1$ SUSY vertex algebra. Let $X$ be an
	$N=1$ supercurve. The vector space spanned by the vacuum vector, the
	current $J$, and the fermion $H$,
	is an $\mathrm{Aut} \cO^{1|1}$-submodule. Indeed, with respect to the
	Virasoro field $\tilde{L}$,  
	the fermion $H$ has conformal weight $1$. Denote the
	corresponding rank $2|1$ vector bundle over $X$ by $\cK(2)_{X,\leq 1}$. It
	follows from 
	the general considerations in appendix \ref{ap:gl1_1_rep},
	that the dual of this vector bundle fits in a short exact
	sequence of the form:
	\begin{equation}
		0 \rightarrow \Omega^1_X \otimes \ber_X \rightarrow \cK(2)_{X,\leq
		1}^* \rightarrow \cO_X \rightarrow 0.
		\label{eq:7.12.1}
	\end{equation}
	Indeed, the relevant relations are in this case:
	\begin{equation*}
		\begin{aligned}
			T_0 J &= J, & T_1 J &= \frac{-c}{3}\vac, & J_1 J &=
			\frac{c}{3} \vac, & H_0 J &= H, \\
			T_0 H &= H, & Q_0 H &= J,  & J_0 H &= - H,  & Q_1
			 H &= \frac{c}{3} \vac,
		\end{aligned}
	\end{equation*}
	therefore the vector space $K(2)_{1}$ spanned by $\{J,H\}$ is isomorphic
	(as a $\fg\fl(1|1)$-module) to $\pi_+(1,0)$ (cf. Appendix
	\ref{ap:gl1_1_rep}),
	and its dual module is then 
	$\pi_-(-1,0) \equiv \pi_{+}(-1,0) \otimes \pi_-(1)$. Also  we know that the
	$\mathrm{Aut} \cO^{1|1}$-twist of $\pi_+(-1,0)$ (resp. $\pi_-(1)$) is
	$\Omega^1_X$ (resp. $\ber_X$). 
	
	We can actually compute these transition functions explicitly as before by
	exponentiating vector fields:
	\begin{equation*}
		\begin{aligned}
			R(\rho)^{-1}\begin{pmatrix}
				\vac \\ J \\ H
			\end{pmatrix} &= A^{T_0} B^{J_0} \exp(h_0 H_0) \exp(q_0
			Q_0) \times \\
			& \quad \times \exp\left( \sum_{i \geq 1} v_i L_i + u_i J_i
			+ q_i Q_i + h_i H_i\right) \cdot\begin{pmatrix}
				\vac \\ J \\ H
			\end{pmatrix} \\
			&=  A^{T_0} B^{J_0} \exp(h_0 H_0) \exp(q_0
			Q_0) \cdot\begin{pmatrix}
				\vac \\ J + (u_1 - v_1) \frac{c}{3} \vac \\ H + q_1
				\frac{c}{3} \vac
			\end{pmatrix} \\
			&=\begin{pmatrix}
				1 & \frac{c}{3} (v_1 - u_1) & \frac{c}{3} q_1 \\
				0 & A & A q_0 \\
				0 & B^{-1} A h_0 & B^{-1} A (1 - h_0 q_0)
			\end{pmatrix}
			\begin{pmatrix}
				\vac \\ J \\H 
			\end{pmatrix},
		\end{aligned}
	\end{equation*}
	which, according to (\ref{eq:7.10.6}), implies:
	\begin{equation}
		R(\rho)^{-1} =\begin{pmatrix}
			1 & \frac{c}{3} \left(\frac{F_{z,\theta}\Psi_z - \Psi_{z,\theta} F_z}{F_\theta \Psi_z -
			\Psi_\theta F_z} + \frac{3}{2} \frac{\Psi_{z,z} F_\theta - F_{z,z}
		\Psi_\theta}{F_z \Psi_\theta - \Psi_z F_\theta}  \right) &
		\frac{c}{6} \frac{F_{z,z} \Psi_z - \Psi_{z,z} F_z}{F_\theta \Psi_z
		- \Psi_\theta F_z} \\
		0 & \frac{F_z \Psi_\theta - \Psi_z F_\theta}{\Psi_\theta} &
		\frac{F_z \Psi_z}{\Psi_\theta} \\
		0 & \frac{F_\theta}{\Psi_\theta} & \frac{F_z^2 \Psi_\theta - \Psi_z
		F_\theta F_z}{\Psi_\theta^2}
		\end{pmatrix}.
		\label{eq:7.12.4}
	\end{equation}
	Taking the super-transpose of the lower two by two block we easily see
	that this block corresponds to the transition functions in $\ber_X \otimes
	\Omega^1$, proving thus that \linebreak $\cK(2)_{X,\leq 1}^*$ is given by an extension
	as in (\ref{eq:7.12.1}). This extension is non-split unless $c = 0$, in
	which case the pair of fields $\{J(z,\theta), H(z,\theta)\}$ transforms as
	a section of $\ber_X \otimes \Omega^1_X$. In order to study the
	splittings of this extension we need to understand the
	differential operators appearing in the first row of (\ref{eq:7.12.4}).
	We leave this to the reader. 
	\label{ex:7.12}
\end{ex}
\section{Chiral algebras on supercurves}\label{sec:chiral}
\begin{nolabel} In this section we follow closely the treatment in chapter 18 of
\cite{frenkelzvi}. We note that most definitions carry over to the ``super'' case
with minor technical changes. In particular we give a sheaf theoretical
interpretation of
the OPE formula (\ref{eq:ope_w.1}).  We
define the superconformal blocks in section \ref{sub:conformalblocks}

We will restrict our analysis to the $(1|1)$ dimensional case for simplicity. All
the results in this section can be generalized to arbitrary odd dimensions
without difficulty. 

For the definitions of chiral algebras over ordinary curves, the reader is
referred to \cite{beilinsondrinfeld} and \cite{frenkelzvi}. For the theory of
$\cD$-modules, we refer to 
\cite{bernstein}, and \cite{penkov} in the supermanifold
case.
\end{nolabel}
\subsection{Chiral algebras}
\begin{nolabel}
	When trying to define chiral algebras on supercurves the first problem
	that we encounter is that given a $(1|N)$ dimensional supercurve $X$
	over $S$,
	the diagonal embedding $\Delta \hookrightarrow X \times_S X$ has relative
	codimension $(1|N)$. In particular, the diagonal is not a divisor in $X
	\times_S X$ unless $N=0$. 
	

	The situation is much simpler in the superconformal case (corresponding
	to $N_K=N$ SUSY vertex algebras). In this case, we can define
	canonically a divisor in $X \times_S X$. Basically, all the arguments in
	the classical case work without change in the superconformal case, given
	that we have replaced the diagonal by a \emph{super diagonal}.

	Since we can carry explicitly the
	computations in the $N=1$ case, without introducing extra notation, we
	will assume that this is the case in the following.
	\label{no:ch.1}
\end{nolabel}
\begin{lem}[6.3 \cite{manin3}] (cf. \ref{no:ch.6} below)
	Let $X$ be a superconformal $N=1$ supercurve. Let $J$ be the ideal
	defining the diagonal $i:\Delta \hookrightarrow X \times_S X$. In local
	coordinates $J$ is defined by $(z-w, \theta - \zeta)$. Let $\Delta^{(1)}$
	be defined by $J^2$. Let $I$ be the kernel of the natural map
	$\Omega^1_{X/S} \rightarrow \mathrm{Ber}_{X/S}$. Finally we define $\Delta^s$ by:
	\begin{equation*}
		\cO_{\Delta^s} = \cO_{\Delta^{(1)}}/i_* (I).
	\end{equation*}
	Then $\Delta^s$ is a $(1|0)$ codimensional divisor in $X \times_S X$,
	locally defined by the equation 
	\begin{equation*}
		0 = z - w - \theta \zeta .
	\end{equation*}
	\label{lem:sdelta}

	This divisor will be called the \emph{super diagonal} and we will simply
	call it the \emph{diagonal} when no confusion should arise. 
\end{lem}
\begin{nolabel}
	Given an $\cO_X$-module $\cM$, we define two extensions of $\cM$ along the
	super diagonal: \emph{extension by principal parts in the transversal
	direction} and \emph{extension by
	delta functions in the transversal direction}. The former is given by
	\begin{equation*}
		\Delta_+^s \cM := \frac{\cO \boxtimes \cM (\infty \Delta^s)}{\cO
		\boxtimes \cM},
	\end{equation*}
	and the latter by
	\begin{equation*}
		\Delta^s_! \cM := \frac{\omega \boxtimes \cM (\infty
		\Delta^s)}{\omega \boxtimes \cM}, 
	\end{equation*}
	where $\omega$ is the Berezinian bundle of $X$ defined in
	\ref{no:berezinian.def}. 
	\label{no:ch.3}
\end{nolabel}
\begin{nolabel}
	As in the non-super case, we have a sheaf-theoretical interpretation of
	the OPE formula. For this we let $X$ be a superconformal $N=1$ curve
	over $\Lambda$.
	 Let $V$ be a strongly conformal
	$N_K=1$ SUSY vertex algebra over $\Lambda$ and let $\cV$ be the associated vector bundle
	over $X$ (cf. \ref{no:k.bundle.def}). Recall that, given any
	$\Lambda$-point $x$ in $X$, we have defined a local section $\cys_x$ (cf.
	\ref{thm:6.8}). Choose local coordinates $Z$ at $x$ compatible with the
	superconformal structure. Using this coordinates we trivialize the
	bundle $\cV$ in the formal superdisk $D_x$ around $x$, namely we have an
	isomorphism $i_Z:V[ [Z]] \rightarrow \Gamma(D_x, \cV)$. Let $W$ be another copy
	of $Z$, so that $D_x^2$ is identified with $\Spec \Lambda[ [Z, W] ]$. The
	bundle $\cV \boxtimes \cV (\infty \Delta^s)$, when restricted to $D_x^2$,
	is the sheaf associated to the $\Lambda[ [ Z, W]] $-module $V \otimes V[
	[Z, W]][ (z-w-\theta\zeta)^{-1}]$. Similarly, the restriction of the
	sheaf $\Delta_+^s \cV$ to $D_x^2$ is associated to the $\Lambda[ [ Z,
	W]]$-module $V[ [Z,W]][(z-w-\theta \zeta)]/V[ [Z,W]]$.
	\label{no:ch.4}
\end{nolabel}
\begin{thm}
	Define a map of $\cO_{D_x^2}$-modules $\cys_{2,x}:\cV \boxtimes \cV(\infty
	\Delta^s) \rightarrow \Delta_+^s \cV$ by the formula
	\begin{equation*}
		\cys_{2,x} ( f(Z,W)a \boxtimes b) = f(Z,W) \ys(a,Z-W)b \mod V[ [
		Z,W]].
	\end{equation*}
	Then $\cys_{2,x}$ is independent of the choice of the coordinates $Z$ as
	long as they are compatible with the superconformal structure induced in
	$D_x$ from that of $X$. 
	\label{thm:ch.5}
\end{thm}
\begin{proof}
	Exactly as in the non-super case, we reduce the proof of this theorem to
	the identity:
	\begin{equation*}
		\ys(a,Z-W) = R(\mu_W) \ys\left( R(\mu_Z)^{-1}a, \mu(Z) - \mu(W)
		\right) R(\mu_W)^{-1}, \qquad a \in V,
	\end{equation*}
	for any $\mu \in \mathrm{Aut}^\omega \cO^{1|1}$. This identity is equivalent
	to (\ref{eq:changewn}) by substituting $Z-W$ instead of $Z$ and $\mu_W(Z-W)
	= \mu(Z) - \mu(W)$ instead of $\rho(Z)$. Recall that in this case we have
	\begin{equation*}
		Z-W = (z - w - \theta\zeta, \theta - \zeta).
	\end{equation*}
\end{proof}
\begin{rem}
	In order to prove a similar statement for a general $N=1$ supercurve $X$
	over $\Lambda$,
	we could define a ``super-diagonal'' as follows. Recall that
	any such curve $X$ gives rise to an oriented superconformal $N=2$ super
	curve $Y$ (cf. \ref{no:n=2superconformal}). Recall also that the
	curve $Y$ comes equipped with two maps $\pi : Y \rightarrow X$ and
	$\hat{\pi}: Y \rightarrow \hat{X}$, where $\hat{X}$ is the \emph{dual} curve.
	In local coordinates these maps are described by (cf. \ref{ex:7.9})
	\begin{equation*}
		\begin{aligned}
			(z, \theta^+, \theta^-) &\xrightarrow{\pi} \left( z +
			\frac{1}{2} \theta^+ \theta^-, \theta^+
			\right), \\
			(z,\theta^+, \theta^-) &\xrightarrow{\hat{\pi}} \left(
			z - \frac{1}{2} \theta^+\theta^-,\theta^- \right)
		\end{aligned}
	\end{equation*}
	It is easy to show that $Y$ embeds as a $(1|0)$ codimensional divisor in
	$X \times_\Lambda \hat{X}$. Indeed, for a $\Lambda$-point $x$ in $X$
	given by local parameters $Z = (z,\theta)$ the preimage in $Y$ is given
	by local parameters $(z - \tfrac{1}{2} \theta \zeta, \theta, \zeta)$.
	Similarly, for a point $W = (w,\zeta)$ in $\hat{X}$ we have its preimage
	in $Y$ given by local parameters $(w + \tfrac{1}{2} \theta \zeta, \theta,
	\zeta)$. Then the point $(Z,W)$ in $X \times_\Lambda \hat{X}$ is in the
	image of $Y$ if and only if $z - w - \theta \zeta = 0$. Note in
	particular that when $X$ is superconformal, namely $X \equiv \hat{X}$
	this ``diagonal'' $Y \hookrightarrow X \times_\Lambda \hat{X}$ agrees with
	Manin's super-diagonal given in Lemma
	\ref{lem:sdelta}. 
	
	We could try to repeat the argument given above for superconformal
	curves, but the operation $\cys_{2}$ turns out to be
	coordinate-dependent\footnote{It will be nice to find a way of describing
	the vertex algebra multiplication as an expression when a point $x \in X$
	``collides'' with a point $\hat{x} \in \hat{X}$ along the ``diagonal'' $\Delta^s
	\subset X \times \hat{X}$.}.
\end{rem}
\begin{nolabel}
	Instead of using the approach in the previous remark, note that we can
	define the push-forwards $\Delta_+$ and $\Delta_!$ even when $\Delta$ is
	not a divisor. In our case these are easy to describe. Let $\Delta$ be
	the diagonal $\Delta \hookrightarrow X \times_S X$. Even though $\Delta$
	is not a divisor in $X \times_S X$, its reduction $|\Delta|$ is a divisor
	in $|X \times_{S} X| = |X| \times_{|S|} |X|$. We have then an open
	immersion $j : X \times X \setminus \Delta \hookrightarrow X \times X$, where
	$X \times X \setminus \Delta$ is $U = |X| \times |X| \setminus |\Delta|$ as a
	topological space and the structure sheaf is the restriction of
	$\cO_{X^2}$ to $U$. We can now define the corresponding push-forwards of an
	$\cO_X$-module $\cM$ as:
	\begin{equation*}
			\Delta_+ \cM = \frac{j_* j^* (\cO_X \boxtimes \cM)}{\cO_X
			\boxtimes \cM}, \qquad
			\Delta_! \cM = \frac{j_* j^* (\omega \boxtimes
			\cM)}{\omega \boxtimes \cM}.
	\end{equation*}
	When no confusion can arise, for any sheaf $\cF$, we will denote by $\cF(\infty \Delta)$ the
	sheaf $j_* j^* \cF$. 
	\label{no:ch.6}
\end{nolabel}
\begin{rem}
	As in the non-super case, these pushforwards are in fact the push forward
	of left (resp. right) $\cD_X$-modules along the diagonal, where in the
	superconformal case we understand for a $\cD_X$ module, a module over
	the ring of superconformal differential operators as in Remark
	\ref{rem:connection} (see also Remark \ref{rem:10b}).
	\label{rem:ch.6b}
\end{rem}
\begin{nolabel}
	We construct now a
	morphism of $\cO_{D^2_x}$-modules $\cys_{2,x} : j_* j^* (\cV_X
	\boxtimes \cV_{X})
	\rightarrow  \Delta^s_+ \cV_{X}$ by the formula:
	\begin{equation}
		\cys_{2,x}( f(Z,W) a \boxtimes b) = f(Z,W) \ys(a, Z-W)b
		\mod V[ [Z,W]].
		\label{eq:ch.7.1}
	\end{equation}
	As in \ref{thm:ch.5} we have
	\label{no:ch.7}
\end{nolabel}
\begin{thm}
	The map $\cys_{2,x}$ defined by (\ref{eq:ch.7.1}) does not depend
	on the coordinates $Z$ chosen. 
	\label{thm:ch.8}
\end{thm}
\begin{nolabel}
	We can now generalize all the results in \cite[chapter 18]{frenkelzvi} on
	chiral algebras 
	without difficulty. For simplicity let us assume that $X$ is a general
	$1|N$-dimensional supercurve. Suppose that the sheaf $\cM$ on $X$
	carries a (left) action of the sheaf of differential operators $\cD_X$. 
	Let $\sigma_{12}:X^2 \rightarrow X^2$ be the transposition of the two
	factors.  
	We obtain a canonical isomorphism of sheaves $\Delta_+ \cM \simeq
	\sigma_{12}^* \Delta_+ \cM$, given in local coordinates by the formula:
	\begin{equation*}
		\frac{1 \otimes \psi}{(Z-W)^{k|K}} \mapsto e^{(Z-W)\nabla} \cdot
		\frac{\psi \otimes 1}{(Z-W)^{k|K}}\mod \cM \boxtimes \cO_X,
	\end{equation*}
	Where $\psi$ is a local section of $\cM$ and $\nabla$ is the connection
	that we obtain from the $\cD$-module structure in $\cM$.
	When $\cM$ carries a
	\emph{right} action of $\cD_X$, we obtain similarly an isomorphism $\Delta_!
	\cM \simeq \sigma_{12}^* \Delta_! \cM$. Note that the Berezinian bundle
	is of rank $(0|1)$ if $N$ is odd, hence in the above formula we need to
	multiply by $(-1)^{\psi N}$ in this case. 

	Similarly, if $X$ is a superconformal curve and $\cM$ carries a (left)
	action of the sheaf of superconformal differential operators $\cD_X$ (cf.
	\ref{rem:connection}), the above formula defines isomorphisms as in
	the general case. 
	\label{no:ch.9}
\end{nolabel}
\begin{nolabel}
	The Berezinian bundle $\omega_X$ is a right
	$D_X$-module, the action given by the Lie derivative \cite{deligne2}.
	Therefore for any left $\cD_X$-module $\cF$ we obtain a right
	$\cD_X$-module $\cF^r := \omega \otimes \cF$. This operation establishes
	an equivalence of categories between left and right $\cD_X$-modules
	\cite{penkov}. 
	The same results hold for $\cD_X$-modules over superconformal curves in
	the sense of \ref{rem:connection}.

	Let $X$ be a supercurve, the sheaf $\omega_X \boxtimes \omega_X$ on
	$X^2$ is isomorphic to $\omega_{X^2}$. The natural map is expressed in
	local coordinates as:
	\begin{equation}
		dZ \boxtimes dW \mapsto [dZdW],
		\label{eq:ch.10.1}
	\end{equation}
	where as before $dZ$ denotes the section $[dz d\theta^1 \dots d\theta^N]$
	of $\omega_X$ and $[dZdW]$ denotes the section $[dz dw d\theta^1 d\zeta^1
	\dots d\theta^N d\zeta^N]$ of $\omega_{X^2}$. We note the skew-symmetry
	in (\ref{eq:ch.10.1}) since (recall the definition of the Berezinian in
	\ref{no:berezinian.def})
	\begin{equation}
		dZ \boxtimes dW \mapsto - (-1)^N [dWdZ].
		\label{eq:ch.10.2}
	\end{equation}
	We obtain thus $\Delta_! \omega_X \simeq \omega_{X^2}(\infty \Delta) /
	\omega_{X^2}$. Let $\mu_\omega$ denote the composition of the
	identification $\omega \boxtimes \omega (\infty \Delta) \simeq
	\omega_{X^2}(\infty \Delta)$ with the projection onto $\Delta_!
	\omega_X$. This map is a morphism of right $\cD_{X^2}$-modules
	satisfying the skew-symmetry condition:
	\begin{equation}
		\mu_\omega \circ \sigma_{12} = - 
		\mu_\omega.
		\label{eq:ch.10.3}
	\end{equation}
	Note that this formula differs from (\ref{eq:ch.10.2}) by a factor
	$(-1)^N$. Indeed this factor appears when applying $\sigma_{12}$, namely
	the composition in the LHS of (\ref{eq:ch.10.3}) is given by:
	\begin{equation*}
		dZ \boxtimes dW \xrightarrow{\sigma_{12}} (-1)^N dW \boxtimes dZ
		\xrightarrow{\mu_\omega} (-1)^N [dWdZ] = - [dZdW] = - \mu_\omega
		dZ \boxtimes dW.
	\end{equation*}
	\label{no:ch.10}
\end{nolabel}
\begin{rem}
	Let $X$ be a supercurve and $Z \hookrightarrow X$ a closed embedding, We
	define the functor $\underline{\Gamma}_Z$ from the category of sheaves on
	$X$ to itself by letting sections of $\underline{\Gamma}_Z(\cF)$ be
	sections of $\cF$ supported on $Z$. This functor is left exact. Let
	$\cH^i_Z$ be the higher derived functors. In this sense the basic
	definitions of local cohomologies in \cite{hartshornelocal} extend in a
	straightforward way to the super case. Similarly we can define the
	\emph{relative local cohomologies} as the higher derived functors of
	$\underline{\Gamma}_{Z/Z'}$ where $Z' \hookrightarrow Z$ is another
	closed embedding and $\underline{\Gamma}_{Z/Z'}$ is defined in the usual
	way as the quotient of sections supported in $Z$ modulo those supported
	in $Z'$ \cite{hartshornelocal}. 
	From the exact sequence
	\begin{equation*}
		0 \rightarrow \underline{\Gamma}_Z (\cF) \rightarrow \cF
		\rightarrow j_* ( \cF|_U) \rightarrow \cH^1_Z (\cF) \rightarrow 0,
	\end{equation*}
	where $U = X \setminus Z$ and $j : U \hookrightarrow X$ is the open
	immersion, we obtain:
	\begin{equation*}
		\Delta_! \omega_X = \cH^1_\Delta(\omega_{X^2}).
	\end{equation*}
	This identification of sheaves extended by delta functions on the
	diagonals with local cohomology sheaves shows that indeed these are
	push-forwards of $\cD_X$-modules.
	\label{rem:10b}
\end{rem}
\begin{nolabel}
	We have also a dictionary between $\cD_X$-modules and delta functions.
	The space $\mathbb{C}[ [ Z^{\pm 1}, W^{\pm 1} ] ] $ carries a structure
	of a module over the algebra of differential operators $\mathbb{C}[ [ Z, W]][
	\nabla_Z, \nabla_W]$ (here $\nabla_Z = (\partial_z, \partial_{\theta^i})$
	in the general case and $\nabla_Z = (\partial_z, D^i_Z)$ in the
	superconformal case). The formal delta-function $\delta(Z,W)$ satisfies
	the relations:
	\begin{equation*}
			(Z - W)^{1|0} \delta(Z,W) = 0, \quad
			(Z - W)^{0|e_i} \delta(Z,W) = 0, \quad
			(\nabla_Z + \nabla_W) \cdot \delta(Z,W) = 0.
	\end{equation*}
	Therefore the $\mathbb{C}[ [Z, W]][ \nabla_Z, \nabla_W]$-submodule of
	$\mathbb{C}[ [ Z^{\pm 1}, W^{\pm 1}]]$ generated by $\delta(Z, W)$ is
	spanned by $\nabla_W^{j|K} \delta(Z,W)$ with $j \geq 0$. This module
	gives rise to a $\cD$-module on the disk $D^2 = \Spec \mathbb{C}[ [ Z,
	W]]$ supported on $z=w$ (note that this is also the case in the
	superconformal case, where the poles are in $z - w - \sum \theta^i \zeta^i$). 
	The assignment
	\begin{equation*}
		(Z-W)^{-1-j|N\setminus J} dW \mapsto \sigma(J)
		\partial^{(j|J)}_W \delta(Z,W), 
	\end{equation*}
	induces an isomorphism of left $\cD$-modules on $D^2$ between $\Delta_+
	\omega$ and the left $\cD$-module generated by $\delta(Z,W)$. Similarly,
	tensoring with $\omega$ we obtain an isomorphism of right $\cD$-modules.
	In the superconformal case the situation is analogous, the proof follows
	from (\ref{eq:k.deriv.1}). 
	\label{no:ch.11}
\end{nolabel}
\begin{nolabel}
	Recall that from Theorem \ref{thm:connection.w} and
	(\ref{eq:rem.connection.1}), we have a natural (left) action of differential operators
	on $\cV$. It follows then that the push-forward $\Delta_+ \cV$ is also a
	(left) $\cD$-module. Indeed, the action of vector fields locally is given
	by ($a \in V$):
	\begin{gather*}
		\partial_z: f(Z,W)a \rightarrow (\partial_z f(Z,W))a, \\
		\partial_w: f(Z,W)a \rightarrow (\partial_w f(Z,W))a + f(Z,W)
		(Ta) \\
		\partial_{\theta^i} : f(Z,W)a \rightarrow (\partial_{\theta^i}
		f(Z,W))a \\ \partial_{\zeta^i}:  f(Z,W)a \rightarrow
		(\partial_{\zeta^i} f(Z,W))a + (-1)^f f(Z,W) S^i a,
	\end{gather*}
	and similarly in the superconformal case, using $D^i_Z$ (resp. $D^i_W$)
	instead of $\partial_{\theta^i}$ (resp. $\partial_{\zeta^i}$).
	Also, we obtain a $\cD$-module structure on the sheaves
	$\cV \boxtimes \cV (\infty \Delta)$ where $\partial_{\theta^i}$ acts as
	$\partial_{\theta^i} + S^i$ and $\partial_{\zeta^i}$ acts as
	$\partial_{\zeta^i} + S^i$. Similarly, in the superconformal case,
	$D^i_Z$ acts as $D^i_Z + S^i$ and $D^i_W$ acts as $D^i_W + S^i$. 
	\label{no:ch.12}
\end{nolabel}
\begin{prop}
	The map $\cY_{2,x}$ commutes with the action of differential operators on
	$D_x^2$, making this map a morphism of $\cD$-modules.
	\label{prop:ch.13}
\end{prop}
\begin{proof}
	For a general supercurve $X$ the proof is the same as in the non-super case.  
	We sketch the proof in the superconformal case where a subtlety arises.
	Let $X = (x,\eta^1, \dots, \eta^N)$. 
	The identity
	\begin{equation*}
		\ys(S^ia, Z-W)b = D_X^i \ys(a, X)b|_{X = Z-W} = D^i_Z \ys(a,
		Z-W)b
	\end{equation*}
	translates into:
	\begin{equation*}
		\cY_{2,x} (D_Z^i \cdot f(Z,W) a \boxtimes b) = D_Z^i \cdot \cY_{2,x}(f
		(Z,W) a \boxtimes b).
	\end{equation*}
	On the other hand, consider translation invariance:
	\begin{multline*}
			{[}S^i, \ys(a,Z-W)]b = (\partial_{\eta^i} - \eta^i
			\partial_x) \ys(a, X)b|_{X=Z-W} =
			\\= (- \partial_{\zeta^i} + \theta^i \partial_x - \eta^i
			\partial_x) \ys(a, Z-W)b|_{X=Z-W} \\
			= (- \partial_{\zeta^i} - \zeta^i \partial_w) \ys(a,
			Z-W)b 
			= - D_W^i \ys(a, Z-W)b.
	\end{multline*}
	From where we obtain:
	\begin{equation*}
		\ys(a, Z-W) S^i b = (-1)^a S^i \ys(a, Z-W)b + (-1)^a D_W^i \ys(a,
		Z-W)b,
	\end{equation*}
	and this translates into:
	\begin{equation*}
		\cY_{2,x}( D_W^i \cdot f(Z,W) a \boxtimes b) = D_W^i \cdot
		\cY_{2,x}(f(Z,W) a \boxtimes b).
	\end{equation*}
\end{proof}
\begin{rem}
	Since $\Delta_+ \cV$ is supported on the diagonal, we obtain a global
	version $\cY^2$ of $\cY_{2,x}$ by gluing these morphisms in the diagonal
	with the zero morphism outside of the diagonal. By the previous
	proposition, this morphism is a map of $\cD$-modules on $X^2$.
	\label{rem:ch.14}
\end{rem}
\begin{prop}
	The map $\cY^2: \cV \boxtimes \cV(\infty \Delta) \rightarrow \Delta_+
	\cV$ satisfies $\cY^2 = \sigma_{12} \circ \cY^2$ under the canonical
	identification $\Delta_+  \simeq \sigma_{12}^* \Delta_+\cV$. 
	\label{prop:ch.15}
\end{prop}
\begin{proof}
	From the skew-symmetry property of SUSY vertex algebras
	(\ref{eq:skew-symmetry.w.1}) it follows:
	\begin{equation*}
		\ys(a,Z-W)b = (-1)^{ab} e^{(Z-W)\nabla} \ys(b, W-Z)a.
	\end{equation*}
	The sign cancels when applying $\sigma_{12}$ 
	and the exponential $e^{(Z-W)\nabla}$ is the coordinate expression for
	the parallel translation, using the $\cD$-module structure on $\cV$, from
	$W$ to $Z$ (see \ref{no:ch.9}).
\end{proof}
\begin{nolabel}
	In order to define \emph{chiral algebras} over supercurves, we need to
	understand the \emph{composition} of morphisms like $\cY^2$. For this we
	need to understand $\Delta_{123!} \cA$ for any right $\cD$-module $\cA$
	over $X$, where $\Delta_{123}$ is the small diagonal in $X^3$ where the
	three points collide. As in the non-super case, we can write this as a
	composition 
	\begin{equation}
		\Delta_{123!} \cA \simeq \Delta_{23!} \Delta_! \cA.
		\label{eq:ch.16.1}
	\end{equation}
	This identity follows from the fact that the push-forward of right
	$\cD$-modules is exact for closed embeddings (cf. \cite{bernstein}). 

	Now let $\mu: \cA \boxtimes \cA (\infty \Delta) \rightarrow \Delta_! \cA$ be a
	morphism of $\cD$-modules on $X^2$. We define a \emph{composition} of
	$\mu$:
	\begin{equation*}
		\mu_{1\{23\}}: j_* \cA \boxtimes \cA \boxtimes \cA|_U \rightarrow
		\Delta_{123!} \cA,
	\end{equation*}
	where $U = X^3 \setminus \cup \Delta_{ij}$ and $j : U \rightarrow X^3$ is
	the open immersion. In order to define such a composition we first apply
	$\mu$ to the second and third argument, and then we apply $\mu$ to the
	first argument and the result (cf. \cite[18.3.1]{frenkelzvi}). We define
	other compositions of $\mu$ by changing the order in which we group the
	points. As in \cite{frenkelzvi} we denote these compositions in the
	following way: given local sections $a, b$ and $c$ of $\cA$ and a
	meromorphic function $f(X,Y,Z)$ with poles along the diagonals, we have:
	\begin{equation*}
		\begin{aligned}
			\mu_{1\{23\}}(f(X,Y,Z) a \boxtimes b \boxtimes c) &= \mu(
			f(X,Y,Z) a \boxtimes \mu(b \boxtimes c)) \\
			\mu_{\{12\}3} (f(X,Y,Z) a \boxtimes b \boxtimes c) &=
			\mu(\mu(f(X,Y,Z)a \boxtimes b) \boxtimes c) \\
			\mu_{2\{13\}} (f(X,Y,Z) a \boxtimes b \boxtimes c) &=
			\sigma_{12} \circ \mu(f(X,Y,Z) b \boxtimes \mu(a
			\boxtimes c)).
		\end{aligned}
	\end{equation*}
	With these compositions defined, we can now define a \emph{chiral
	algebra} in the usual way:
	\label{no:ch.16}
\end{nolabel}
\begin{defn}
	A \emph{chiral algebra} on a $1|N$ dimensional supercurve $X$ is a right $\cD$-module
	$\cA$ equipped with a morphism of $\cD$-modules:
		$\mu: \cA \boxtimes \cA (\infty \Delta) \rightarrow \Delta_! \cA$, 
	satisfying the following conditions:
	\begin{itemize}
		\item (skew-symmetry) $\mu = - \mu \circ
			\sigma_{12}$.
		\item (Jacobi identity) $\mu_{1\{23\}} = \mu_{\{12\}3} +
			\mu_{2\{13\}}$. 
		\item (Unit) We are given a canonical embedding $\omega_X
			\hookrightarrow \cA$ of the Berezinian bundle compatible
			with the homomorphism $\mu_\omega$ defined in
			\ref{no:ch.10}. 
	\end{itemize}
	\label{defn:ch.17}
\end{defn}
\begin{rem}
	Note that this definition is exactly the same as in the non-super case,
	namely, the signs appearing when anticommuting odd-elements are taken
	care by the symmetric structure of the category of modules over
	super-rings. That is, given a super-ring $R$ and two $R$-modules $M$ and
	$N$, the isomorphism $\sigma: M \otimes N \simeq N \otimes M$ is given
	by:
	\begin{equation*}
		\sigma: \; m \otimes n \mapsto (-1)^{mn} n \otimes m.
	\end{equation*}
	Indeed the only difference with the non-super case is the fact that the
	unit $\omega$ is a rank $(0|1)$-bundle when $N$ is odd. From
	the SUSY vertex algebra point of view, this is translated into the fact that
	the $\Lambda$-bracket has parity $N \mod 2$. 

	In the superconformal case there is a subtlety. We note that the
	intersection of two different diagonals in the sense of \ref{lem:sdelta}
	depends on the diagonals chosen, namely:
	\begin{equation*}
		\Delta^s_{12} \cap \Delta^s_{23} \neq \Delta^s_{13} \cap
		\Delta^s_{23}.
	\end{equation*}
	But despite this fact, the pushforward $\Delta_{123!}$ is still
	well defined, independent of the composition chosen as in
	(\ref{eq:ch.16.1}). 
	\label{rem:18}
\end{rem}
Using the equivalence between left $\cD$-modules and right $\cD$-modules, we
obtain a right $\cD$-module $\cV^r = \omega_X \otimes \cV$ from any
strongly conformal SUSY vertex algebra. Similarly, this sheaf carries a
multiplication $\mu = (\cY^2)^r$ obtained from $\cY^2$. 
\begin{thm}
	The pair $(\cV^r, \mu)$ carries a structure of a chiral algebra over $X$. 
	\label{thm:ch.19}
\end{thm}
\begin{proof}
	The proof of this fact is the same as the proof in the non-super case
	\cite[Thm 18.3.3]{frenkelzvi}. This follows by considering the
	\emph{Cousin resolution} of the Berezinian bundle in $X^3$ and the
	corresponding Cousin property of SUSY vertex algebras 
	\ref{cor:cousin.w} proved in \cite{heluani3}. 
\end{proof}
\subsection{Conformal blocks} \label{sub:conformalblocks}
In this section we define the sheaves of coinvariants of SUSY vertex algebras.
The treatment follows \cite{frenkelzvi}. In fact, most results carry over without
change to our situation. We only mention the major differences.
\begin{nolabel} 
	Recall that the polar part of a SUSY vertex algebra is naturally a SUSY Lie
	conformal
	algebra (cf. \cite{heluani3}). We can consider then the operator $\cY_{x,-}$ which is the polar
	part of $\cY_x$. The notion of $\mathrm{Lie}^*$ algebra over a super
	curve is generalized in a straightforward manner from the non-super case. 
	\label{no:cb.1}
\end{nolabel}
\begin{nolabel}
	Let $\cA$ be a right $\cD$-module, the de Rham sequence of $\cA$ is the
	sequence:
	\begin{equation*}
		0 \rightarrow \cA \otimes \cT \rightarrow \cA \rightarrow 0
	\end{equation*}
	placed in cohomological degrees $0$ and $-1$, where $\cT$ is the tangent
	sheaf of $X$. In the superconformal case, we do not have an action of
	the entire tangent sheaf, but we can act by the subsheaf $\cT^s$ generated by the
	derivations $D^i_Z$ (i.e. the subsheaf $\cT^1$ of remark \ref{rem:3.14}
	in the $1|1$ dimensional case, and the sheaf $\cT' \oplus \cT''$ in the
	$1|2$ dimensional case). We define then \emph{the de Rham sheaf} $h(\cA)$
	of $\cA$ as
	\begin{equation*}
		h(\cA) = \cA / (\cA \cdot \cT).
	\end{equation*}
	whereas in the superconformal case we put $h(\cA) = \cA/(\cA \cdot
	\cT^s)$.
	\label{no:cb.2}
\end{nolabel}
\begin{prop}
	Let $(\cA, \mu)$ be a chiral algebra. Then 
	\begin{enumerate}
		\item $h(\cA)(D_x^\times)$ and $h(\cA)(\Sigma)$, for any open
			$x \notin \Sigma \subset X$ are Lie superalgebras, and there is a
			natural homomorphism of Lie superalgebras $h(\cA)(\Sigma) 
			\rightarrow h(\cA)(D_x^\times)$. 
		\item $h(\cA)(D_x^\times)$ acts on the fiber $\cA_x$. 
		\item If $(\cA, \mu)$ is associated to a SUSY vertex algebra $V$, then
			there is a canonical isomorphism $h(\cA)(D_x^\times)
			\simeq \lie'(V)$ (see Theorem \ref{thm:lie_fourier.vertex.1}
			for the definition of $\lie'(V)$). 
	\end{enumerate}
	\label{prop:cb.3}
\end{prop}
\begin{proof}
	We can think of $\cA \simeq \omega \otimes \cA^l$, where $\cA^l$ is a
	left $\cD$-module. Since we can integrate sections of the Berezinian
	bundle, we see
	immediately that we have $h(\Delta_! \cA) = \Delta_* h(\cA)$. On the
	other hand the map $\mu : \cA \boxtimes \cA(\infty \Delta) \rightarrow
	\Delta_! \cA$ induces
	\begin{equation*}
		h(\mu) : h(\cA) \boxtimes h(\cA) (\infty \Delta) \rightarrow
		h(\Delta_! \cA).
	\end{equation*}
	Restricting to regular sections and pulling back along the
	diagonal we obtain:
	\begin{equation*}
		[\,,\,] : h(\cA) \otimes h(\cA) \rightarrow h(\cA).
	\end{equation*}
	The fact that $[\,,\,]$ satisfies the axioms of a Lie superalgebra
	follows from the skew-symmetry and Jacobi identity of
	chiral algebras.  The rest of the theorem is proved in the
	same way as \cite[prop 18.4.12]{frenkelzvi}. 
	
	(3)
	follows from the definitions, in formulas (\ref{eq:acabodeagregarte}).
	Indeed, these formulas are the equivalent of the corresponding formulas
	for the action of vector fields on $\cA^l$ as defined in
	Theorem \ref{thm:connection.w} and in (\ref{eq:rem.connection.1}).
\end{proof}
\begin{rem}
	As in the non-super case, for a strongly conformal SUSY vertex algebra $V$,
	we have a natural map 
	\begin{equation*}
		\cys^\vee_x : \cV^r (D_x^\times) \rightarrow \End(\cV_x) \simeq
		\End \cV^r_x,
	\end{equation*}
	on $D_x^\times$. Namely, given a section $s \in \cV^r(D_x^\times)$ we
	obtain the endomorphism $\cys^\vee_x (s) = \res_X <\cys_x, s>$ on $\cV_x$.
	If $s$ is a total derivative, this residue vanishes and the map
	$\cys^\vee_x$ factors through $h(\cV^r)(D^\times_x)$. The resulting Lie
	superalgebra homomorphism $h(\cV^r)(D^\times_x) \rightarrow \End (\cV^r_x)$
	coincides with the homomorphism of Proposition \ref{prop:cb.3} (2) and with the
	homomorphism $\varphi'$ of Theorem \ref{thm:lie_fourier.vertex.1}.
\end{rem}
\begin{nolabel}
	We can now define the \emph{spaces of coinvariants} for a SUSY
	vertex algebra. For this let $X$ be a supercurve and $x \in X$ a point.
	We have a Lie superalgebra $U_\Sigma = h(\cV^r)(\Sigma)$, where $\Sigma =
	X \setminus
	\{x\}$ and this Lie superalgebra acts in $\cV_x$. 
	\label{no:cb.4}
\end{nolabel}
\begin{defn}
	The \emph{space of coinvariants} associated to $(V, X, x)$ is
	\begin{equation*}
		H(V,X,x) = \cV_x / (U_\Sigma \cdot \cV_x).
	\end{equation*}
	\label{defn:cb.5}
\end{defn}
\begin{rem}
	The extension of this definition to the multiple point case with
	arbitrary module insertions is straightforward and we leave it for the
	reader.
	\label{rem:cb.6}
\end{rem}
     Fix $N \geq 0$. Let $\fg$ be the Lie superalgebra of vector fields on the $1|N$
     dimensional punctured superdisk $D^\times$, namely $\fg$ is the completion
     of the Lie superalgebra $W(1|N)$. Let $\fg^\omega$ be the Lie subalgebra of
     $\fg$ consisting of vector fields preserving the form $\omega = dt + \sum
     \zeta^i d\zeta^i$, namely $\fg^\omega$ is the completion of the Lie
     superalgebra $K(1|N)$. Let
     $\cM_{g,1}$ be the moduli space of smooth $1|N$ dimensional genus $g$, pointed
     supercurves (here the genus of a supercurve $X$ is the genus of
     $X_{\mathrm{rd}}$). Let $\hat{\cM}_{g,1}$ be the moduli space of triples $(X, x,
     Z)$, where $(X,x) \in \cM_{g,1}$ and $Z$ is a coordinate system at $x$. Let
     $\cM_{g,1}^\omega$ and $\hat{\cM}_{g,1}^\omega$ be the superconformal analogous.
     \begin{thm}[{\cite{vaintrob2}}] The Lie algebra $\fg$ (resp. $\fg^\omega$) acts
       (infinitesimally) transitively on $\hat{\cM}_{g,1}$ (resp.
       $\hat{\cM}_{g,1}^\omega$). This action preserves the fibers of the projection
       $\hat{\cM}_{g,1} \rightarrow \cM_{g,1}$ (resp. $\hat{\cM}_{g,1}^\omega \rightarrow
       \cM_{g,1}^\omega$). 
\end{thm}

It follows from this theorem, by repeating the \emph{localization} construction in
\cite[ch. 16]{frenkelzvi} that, given a strongly conformal $N_W=n$ SUSY vertex
algebra (resp. a strongly conformal 
$N_K=n$ SUSY vertex algebra) $V$, we obtain a left $\cD$-module $\Delta(V)$ on
$\cM_{g,1}$ (resp. $\cM^\omega_{g,1}$), whose fiber at $(X,x)$ is the space of coinvariants
$H(X,x,V)$.

\appendix

\section{Representations of $\fg\fl(1|1)$} \label{ap:gl1_1_rep}
Let us pick a basis of $\mathfrak{g}\mathfrak{l}(1|1)$ such that
\begin{equation*}
		T = \begin{pmatrix}
			1 & 0\\
			0& 1
		\end{pmatrix}, \quad
		J =\begin{pmatrix}
			0 & 0 \\ 0 & 1
		\end{pmatrix}, \quad
		Q =\begin{pmatrix}
			0 & 1 \\
			0 & 0 
		\end{pmatrix}, \quad
		H =\begin{pmatrix}
			0 &0 \\ 1 & 0 
		\end{pmatrix}. 
\end{equation*}
Then the irreducible representations such that $T$ and $J$ act diagonally are
classified by
\begin{itemize}
	\item $1 | 0$ or $0|1$ dimensional: these are representation on
		$\mathbb{C}^{1|0}$ or $\mathbb{C}^{0|1}$ generated by an even
		(resp. odd) vector $\bar{1} \in \mathbb{C}$ such that in this basis
		we have $T = Q = H = 0$ and $J = j$ we call these representations
		$\pi_{\pm}(j)$. 
	\item $1|1$ dimensional: for each numbers $t, j \in \mathbb{C}$ there are
		two irreducible representations of dimension $1|1$. These are
		either of highest $\pi_+(t,j)$, or lowest weight $\pi_-(t,j)$:
		\begin{xalignat*}{4}
			T &=\begin{pmatrix}
				t & 0 \\ 0 & t
			\end{pmatrix}, & J &=\begin{pmatrix}
				j & 0 \\ 0 & j-1
			\end{pmatrix}, & Q &=\begin{pmatrix}
				0 & t \\ 0 & 0
			\end{pmatrix}, & H &=\begin{pmatrix}
				0 & 0 \\ 1 & 0
			\end{pmatrix}, \\
				T &=\begin{pmatrix}
				t & 0 \\ 0 & t
			\end{pmatrix}, & J &=\begin{pmatrix}
				j & 0 \\ 0 & j+1
			\end{pmatrix}, & Q &=\begin{pmatrix}
				0 & 0 \\ 1 & 0
			\end{pmatrix}, & H &=\begin{pmatrix}
				0 & t \\ 0 & 0
			\end{pmatrix}. 
		\end{xalignat*}
		
\end{itemize}
We note that by taking minus the super transpose we get that the duals of these
representations are given (in the dual basis $\{v^*, \omega^* \}$) by 
\begin{xalignat*}{4}
			T &=\begin{pmatrix}
				-t & 0 \\ 0 & -t
			\end{pmatrix}, & J &=\begin{pmatrix}
				-j & 0 \\ 0 & -j+1
			\end{pmatrix}, & Q &=\begin{pmatrix}
				0 & 0 \\ -t & 0
			\end{pmatrix}, & H &=\begin{pmatrix}
				0 & 1 \\ 0 & 0
			\end{pmatrix}, \\
				T &=\begin{pmatrix}
				-t & 0 \\ 0 & -t
			\end{pmatrix}, & J &=\begin{pmatrix}
				-j & 0 \\ 0 & -j-1
			\end{pmatrix}, & Q &=\begin{pmatrix}
				0 & 1 \\ 0 & 0
			\end{pmatrix}, & H &=\begin{pmatrix}
				0 & 0 \\ -
				t & 0
			\end{pmatrix}, 
\end{xalignat*}
which in the basis $\{-t^{-1} v,  \omega\}$ show that $\pi_\pm (t, j)^\vee \equiv
\pi_\mp (-t,-j)$

Finally we note that the parity changed modules are $\Pi \pi_\pm(t,j) = \pi_\mp
(t,j \mp 1)$. 

On the formal $1|1$ dimensional superdisk with coordinates $(z,\theta)$ we have the following
realization of these representations. Consider the basis for this
Lie algebra $-T = z \partial_z + \theta \partial_\theta$, $J = - \theta
\partial_\theta$, $Q = - z \partial_\theta$ and $H = \theta \partial_z$ acting on
sections of a vector bundle by the Lie derivative. By analyzing the action of
these derivations on the fibers of the corresponding bundles we obtain: 
\begin{equation*}
	\begin{aligned}
		\wedge^m \Omega^1 &= \mathrm{Aut_D} \overset{Aut \cO}{\times}
		\pi_+(-m, -m + 1) \qquad m \equiv 1 (2), \\
		\wedge^m \Omega^1 &= \mathrm{Aut_D} \overset{Aut \cO}{\times}
		\pi_-(-m,-m) \qquad m \equiv 0 (2), \\
		S^m \Omega^1 &= \mathrm{Aut_D} \overset{Aut \cO}{\times} \pi_+ (-m,0)\\ 
		\mathrm{Ber}_D &= \mathrm{Aut_D} \overset{Aut \cO}{\times} \pi_-
		(1). \\
	\end{aligned}
\end{equation*}

\bibliographystyle{alpha}
\bibliography{refs}
\end{document}